\documentclass[10pt]{amsart}
\usepackage{amsmath}
\usepackage{amsxtra}
\usepackage{amscd}
\usepackage{amsthm}
\usepackage{amsfonts}
\usepackage{amssymb}
\usepackage{eucal}

\newtheorem{cor}[subsection]{Corollary}
\newtheorem{lem}[subsection]{Lemma}
\newtheorem{prop}[subsection]{Proposition}

\newtheorem{conj}[subsection]{Conjecture}
\newtheorem{thm}[subsection]{Theorem}

\theoremstyle{definition}

\theoremstyle{remark}

\numberwithin{equation}{section}

\newcommand{\thmref}[1]{Theorem~\ref{#1}}
\newcommand{\secref}[1]{Sect.~\ref{#1}}
\newcommand{\lemref}[1]{Lemma~\ref{#1}}
\newcommand{\propref}[1]{Proposition~\ref{#1}}
\newcommand{\corref}[1]{Corollary~\ref{#1}}

\newcommand{\nc}{\newcommand}
\nc{\renc}{\renewcommand}
\nc{\ssec}{\subsection}
\nc{\sssec}{\subsubsection}
\nc{\on}{\operatorname}

\nc\ol{\overline}
\nc\wt{\widetilde}
\nc\tboxtimes{\wt{\boxtimes}}
\nc{\alp}{\alpha}

\nc{\ZZ}{{\mathbb Z}}
\nc{\NN}{{\mathbb N}}
\nc{\CC}{{\mathbb C}}
\nc{\OO}{{\mathbb O}}
\renc{\SS}{{\mathbb S}}
\nc{\DD}{{\mathbb D}}
\nc{\GG}{{\mathbb G}}

\nc{\Fq}{{\mathbb F}_q}
\nc{\Fqb}{\ol{{\mathbb F}_q}}
\nc{\Ql}{\ol{{\mathbb Q}_\ell}}
\nc{\id}{\text{id}}
\nc\X{\mathcal X}

\nc{\Hom}{\on{Hom}}
\nc{\Lie}{\on{Lie}}
\nc{\Loc}{\on{Loc}}
\nc{\Pic}{\on{Pic}}
\nc{\Bun}{\on{Bun}}
\nc{\IC}{\on{IC}}
\nc{\Aut}{\on{Aut}}
\nc{\rk}{\on{rk}}
\nc{\Sh}{\on{Sh}}
\nc{\Perv}{\on{Perv}}
\nc{\pos}{{\on{pos}}}
\nc{\Conv}{\on{Conv}}
\nc{\Sph}{\on{Sph}}
\nc{\Sym}{\on{Sym}}
\nc{\BunBb}{\overline{\Bun}_B}
\nc{\Buno}{\overset{o}{\Bun}}
\nc{\BunPb}{{\overline{\Bun}_P}}
\nc{\BunBM}{\overline{\Bun}_{B(M)}}
\nc{\BunPbw}{{\widetilde{\Bun}_P}}
\nc{\BunBP}{\widetilde{\Bun}_{B,P}}
\nc{\GUb}{\overline{G/U}}
\nc{\GUPb}{\overline{G/U(P)}}

\nc{\Hhom}{\underline{\on{Hom}}}
\nc\syminfty{\on{Sym}^{\infty}}
\nc\lal{\ol{\lambda}}
\nc\xl{\ol{x}}
\nc\thl{\ol{\theta}}
\nc\nul{\ol{\nu}}
\nc\mul{\ol{\mu}}
\nc\Sum\Sigma
\nc{\oX}{\overset{o}{X}{}}
\nc{\hl}{\overset{\leftarrow}h}
\nc{\hr}{\overset{\rightarrow}h}
\nc{\M}{{\mathcal M}}
\nc{\N}{{\mathcal N}}
\nc{\F}{{\mathcal F}}
\nc{\D}{{\mathcal D}}
\nc{\Q}{{\mathcal Q}}
\nc{\Y}{{\mathcal Y}}
\nc{\G}{{\mathcal G}}
\nc{\E}{{\mathcal E}}
\nc{\CalC}{{\mathcal C}}
\nc\Dh{\widehat{\D}}

\nc{\C}{{\mathcal C}}
\nc{\K}{{\mathcal K}}
\renewcommand{\H}{{\mathcal H}}

\nc{\T}{{\mathcal T}}
\nc{\V}{{\mathcal V}}
\renc{\P}{{\mathcal P}}
\nc{\A}{{\mathcal A}}
\nc{\B}{{\mathcal B}}
\nc{\U}{{\mathcal U}}

\nc{\Gr}{\on{Gr}}

\nc{\frn}{{\check{\mathfrak u}(P)}}
\nc{\p}{\mathfrak p}
\nc{\q}{\mathfrak q}
\nc\f{{\mathfrak f}}

\nc{\qo}{{\mathfrak q}}
\nc{\po}{{\mathfrak p}}
\nc{\s}{{\mathfrak s}}
\nc\w{\text{w}}

\nc\mathi\iota
\nc\Spec{\on{Spec}}
\nc\Mod{\on{Mod}}
\nc{\tw}{\widetilde{\mathfrak t}}
\nc{\pw}{\widetilde{\mathfrak p}}
\nc{\qw}{\widetilde{\mathfrak q}}
\nc{\jw}{\widetilde j}

\nc{\grb}{\overline{\Gr}}
\nc{\I}{\mathcal I}

\nc{\lambdach}{{\check\lambda}}
\nc{\Lambdach}{{\check\Lambda}{}}
\nc{\much}{{\check\mu}}
\nc{\omegach}{{\check\omega}}
\nc{\nuch}{{\check\nu}}
\nc{\etach}{{\check\eta}}
\nc{\alphach}{{\check\alpha}}
\nc{\betach}{{\check\beta}}
\nc{\rhoch}{{\check\rho}}
\nc{\ch}{{\check h}}

\nc{\Hb}{\overline{\H}}

\nc{\BA}{{\mathbb{A}}}
\nc{\BC}{{\mathbb{C}}}
\nc{\BG}{{\mathbb{G}}}
\nc{\BM}{{\mathbb{M}}}
\nc{\BN}{{\mathbb{N}}}
\nc{\BP}{{\mathbb{P}}}
\nc{\BR}{{\mathbb{R}}}
\nc{\BZ}{{\mathbb{Z}}}
\nc{\BV}{{\mathbb{V}}}
\nc{\BW}{{\mathbb{W}}}
\nc{\BS}{{\mathbb{S}}}
\nc{\BQ}{{\mathbb{Q}}}

\nc{\CA}{{\mathcal{A}}}
\nc{\CB}{{\mathcal{B}}}

\nc{\CE}{{\mathcal{E}}}
\nc{\CF}{{\mathcal{F}}}
\nc{\CG}{{\mathcal{G}}}
\nc{\CL}{{\mathcal{L}}}
\nc{\CM}{{\mathcal{M}}}
\nc{\CN}{{\mathcal{N}}}
\nc{\CK}{{\mathcal{K}}}
\nc{\CO}{{\mathcal{O}}}
\nc{\CP}{{\mathcal{P}}}
\nc{\CJ}{{\mathcal{J}}}
\nc{\CI}{{\mathcal{I}}}
\nc{\CQ}{{\mathcal{Q}}}
\nc{\CR}{{\mathcal{R}}}
\nc{\CS}{{\mathcal{S}}}
\nc{\CT}{{\mathcal{T}}}
\nc{\CU}{{\mathcal{U}}}
\nc{\CV}{{\mathcal{V}}}
\nc{\CW}{{\mathcal{W}}}
\nc{\CZ}{{\mathcal{Z}}}

\nc{\cM}{{\check{\mathcal M}}{}}
\nc{\csM}{{\check{\mathcal A}}{}}
\nc{\oM}{{\overset{\circ}{\mathcal M}}{}}
\nc{\obM}{{\overset{\circ}{\mathbf M}}{}}
\nc{\oCA}{{\overset{\circ}{\mathcal A}}{}}
\nc{\obA}{{\overset{\circ}{\mathbf A}}{}}
\nc{\ooM}{{\overset{\circ}{M}}{}}
\nc{\osM}{{\overset{\circ}{\mathsf M}}{}}
\nc{\vM}{{\overset{\bullet}{\mathcal M}}{}}
\nc{\nM}{{\underset{\bullet}{\mathcal M}}{}}
\nc{\oD}{{\overset{\circ}{\mathcal D}}{}}
\nc{\obD}{{\overset{\circ}{\mathbf D}}{}}
\nc{\oA}{{\overset{\circ}{\mathbb A}}{}}
\nc{\op}{{\overset{\bullet}{\mathbf p}}{}}
\nc{\cp}{{\overset{\circ}{\mathbf p}}{}}
\nc{\oU}{{\overset{\bullet}{\mathcal U}}{}}
\nc{\oZ}{{\overset{\circ}{\mathcal Z}}{}}
\nc{\ofZ}{{\overset{\circ}{\mathfrak Z}}{}}
\nc{\oF}{{\overset{\circ}{\fF}}}

\nc{\fa}{{\mathfrak{a}}}
\nc{\fb}{{\mathfrak{b}}}
\nc{\fg}{{\mathfrak{g}}}
\nc{\fgl}{{\mathfrak{gl}}}
\nc{\fh}{{\mathfrak{h}}}
\nc{\fj}{{\mathfrak{j}}}
\nc{\fm}{{\mathfrak{m}}}
\nc{\fl}{{\mathfrak{l}}}
\nc{\fn}{{\mathfrak{n}}}
\nc{\fu}{{\mathfrak{u}}}
\nc{\fp}{{\mathfrak{p}}}
\nc{\fr}{{\mathfrak{r}}}
\nc{\fs}{{\mathfrak{s}}}
\nc{\fsl}{{\mathfrak{sl}}}
\nc{\hsl}{{\widehat{\mathfrak{sl}}}}
\nc{\hgl}{{\widehat{\mathfrak{gl}}}}
\nc{\hg}{{\widehat{\mathfrak{g}}}}
\nc{\chg}{{\widehat{\mathfrak{g}}}{}^\vee}
\nc{\hn}{{\widehat{\mathfrak{n}}}}
\nc{\chn}{{\widehat{\mathfrak{n}}}{}^\vee}

\nc{\fA}{{\mathfrak{A}}}
\nc{\fB}{{\mathfrak{B}}}
\nc{\fO}{{\mathfrak{O}}}
\nc{\fD}{{\mathfrak{D}}}
\nc{\fE}{{\mathfrak{E}}}
\nc{\fF}{{\mathfrak{F}}}
\nc{\fG}{{\mathfrak{G}}}
\nc{\fK}{{\mathfrak{K}}}
\nc{\fL}{{\mathfrak{L}}}
\nc{\fM}{{\mathfrak{M}}}
\nc{\fN}{{\mathfrak{N}}}
\nc{\fP}{{\mathfrak{P}}}
\nc{\fU}{{\mathfrak{U}}}
\nc{\fV}{{\mathfrak{V}}}
\nc{\fZ}{{\mathfrak{Z}}}
\nc{\fz}{{\mathfrak{z}}}

\nc{\bb}{{\mathbf{b}}}
\nc{\bc}{{\mathbf{c}}}
\nc{\bd}{{\mathbf{d}}}
\nc{\be}{{\mathbf{e}}}
\nc{\bj}{{\mathbf{j}}}
\nc{\bn}{{\mathbf{n}}}
\nc{\bp}{{\mathbf{p}}}
\nc{\bq}{{\mathbf{q}}}
\nc{\bu}{{\mathbf{u}}}
\nc{\bv}{{\mathbf{v}}}
\nc{\bx}{{\mathbf{x}}}
\nc{\bs}{{\mathbf{s}}}
\nc{\by}{{\mathbf{y}}}
\nc{\bw}{{\mathbf{w}}}
\nc{\bA}{{\mathbf{A}}}
\nc{\bK}{{\mathbf{K}}}
\nc{\bB}{{\mathbf{B}}}
\nc{\bC}{{\mathbf{C}}}
\nc{\bD}{{\mathbf{D}}}
\nc{\bH}{{\mathbf{H}}}
\nc{\bM}{{\mathbf{M}}}
\nc{\bN}{{\mathbf{N}}}
\nc{\bV}{{\mathbf{V}}}
\nc{\bW}{{\mathbf{W}}}
\nc{\bX}{{\mathbf{X}}}
\nc{\bZ}{{\mathbf{Z}}}
\nc{\bS}{{\mathbf{S}}}

\nc{\sA}{{\mathsf{A}}}
\nc{\sB}{{\mathsf{B}}}
\nc{\sC}{{\mathsf{C}}}
\nc{\sD}{{\mathsf{D}}}
\nc{\sF}{{\mathsf{F}}}
\nc{\sG}{{\mathsf{G}}}
\nc{\sK}{{\mathsf{K}}}
\nc{\sM}{{\mathsf{M}}}
\nc{\sO}{{\mathsf{O}}}
\nc{\sQ}{{\mathsf{Q}}}
\nc{\sP}{{\mathsf{P}}}
\nc{\sZ}{{\mathsf{Z}}}
\nc{\sfp}{{\mathsf{p}}}
\nc{\sr}{{\mathsf{r}}}
\nc{\sg}{{\mathsf{g}}}
\nc{\ssf}{{\mathsf{f}}}
\nc{\sfb}{{\mathsf{b}}}
\nc{\sfc}{{\mathsf{c}}}
\nc{\sd}{{\mathsf{d}}}

\nc{\BK}{{\bar{K}}}

\nc{\tA}{{\widetilde{\mathbf{A}}}}
\nc{\tB}{{\widetilde{\mathcal{B}}}}
\nc{\tg}{{\widetilde{\mathfrak{g}}}}
\nc{\tG}{{\widetilde{G}}}
\nc{\TM}{{\widetilde{\mathbb{M}}}{}}
\nc{\tO}{{\widetilde{\mathsf{O}}}{}}
\nc{\tU}{{\widetilde{\mathfrak{U}}}{}}
\nc{\TZ}{{\tilde{Z}}}
\nc{\tx}{{\tilde{x}}}
\nc{\tbv}{{\tilde{\bv}}}
\nc{\tfP}{{\widetilde{\mathfrak{P}}}{}}
\nc{\tz}{{\tilde{\zeta}}}
\nc{\tmu}{{\tilde{\mu}}}

\nc{\urho}{\underline{\rho}}
\nc{\uB}{\underline{B}}
\nc{\uC}{{\underline{\mathbb{C}}}}
\nc{\ui}{\underline{i}}
\nc{\uj}{\underline{j}}
\nc{\ofP}{{\overline{\mathfrak{P}}}}
\nc{\oB}{{\overline{\mathcal{B}}}}
\nc{\og}{{\overline{\mathfrak{g}}}}
\nc{\oI}{{\overline{I}}}

\nc{\eps}{\varepsilon}
\nc{\hrho}{{\hat{\rho}}}

\nc{\one}{{\mathbf{1}}}
\nc{\two}{{\mathbf{t}}}

\nc{\Rep}{{\mathop{\operatorname{\rm Rep}}}}
\nc{\Tot}{{\mathop{\operatorname{\rm Tot}}}}
\nc{\Ker}{{\mathop{\operatorname{\rm Ker}}}}
\nc{\Hilb}{{\mathop{\operatorname{\rm Hilb}}}}
\nc{\End}{{\mathop{\operatorname{\rm End}}}}
\nc{\Ext}{{\mathop{\operatorname{\rm Ext}}}}
\nc{\CHom}{{\mathop{\operatorname{{\mathcal{H}}\it om}}}}
\nc{\GL}{{\mathop{\operatorname{\rm GL}}}}
\nc{\gr}{{\mathop{\operatorname{\rm gr}}}}
\nc{\Id}{{\mathop{\operatorname{\rm Id}}}}
\nc{\de}{{\mathop{\operatorname{\rm def}}}}
\nc{\length}{{\mathop{\operatorname{\rm length}}}}
\nc{\supp}{{\mathop{\operatorname{\rm supp}}}}

\nc{\Cliff}{{\mathsf{Cliff}}}
\nc{\Fl}{\on{Fl}}
\nc{\Fib}{{\mathsf{Fib}}}
\nc{\Coh}{{\mathsf{Coh}}}
\nc{\FCoh}{{\mathsf{FCoh}}}

\nc{\cplus}{{\mathbf{C}_+}}
\nc{\cminus}{{\mathbf{C}_-}}
\nc{\cthree}{{\mathbf{C}_*}}
\nc{\Qbar}{{\bar{Q}}}

\nc{\bh}{{\bar{h}}}

\nc{\seq}[1]{\stackrel{#1}{\sim}}

\nc{\wh}{\widehat}

\nc{\mc}{\mathcal}
\nc{\crit}{\on{crit}}
\nc{\reg}{\on{reg}}
\nc{\nil}{\wt{\on{reg}}}
\nc{\ka}{\kappa}
\nc{\g}{\fg}
\nc{\mb}{\mathbf}
\nc{\ren}{ren}
\nc{\la}{\lambda}
\nc{\FZ}{{\mathfrak Z}}
\nc{\Z}{{\mathbb Z}}

\begin{document}

\title[D-modules on the affine Grassmannian]
{D--modules on the affine Grassmannian and representations
of affine Kac-Moody algebras}

\dedicatory{Dedicated to Victor Kac on his 60th birthday}

\author{Edward Frenkel}

\address{Department of Mathematics, University of California,
  Berkeley, CA 94720, USA}

\author{Dennis Gaitsgory}

\address{Department of Mathematics, The University of Chicago,
  Chicago, IL 60637, USA}

\date{November 2003}

\maketitle

\section{Introduction}

\ssec{}    \label{negative}

Let $\fg$ be a simple Lie algebra over $\BC$, and $G$ the
corresponding algebraic group of adjoint type. Given an invariant
inner product $\kappa$ on $\g$, let $\hat\fg_\kappa$ denote the
corresponding central extension of the formal loop algebra $\fg
\otimes \BC((t))$, called the affine Kac-Moody algebra $\hat\fg_\kappa$,
$$0\to \BC {\mb 1}\to \hat\fg_\kappa\to \fg \otimes \BC((t))\to 0,$$
with the two-cocycle defined by the formula $$x \otimes f(t),y \otimes
g(t) \mapsto -\kappa(x,y)\cdot \on{Res}_{t=0} f dg.$$ Denote by
$\hat\fg_\kappa\text{--}\on{mod}$ the category of $\hat\fg_\kappa$-modules
which are discrete, i.e., any vector is annihilated by the Lie
subalgebra $\g \otimes t^N \BC[[t]]$ for sufficiently large $N\geq 0$,
and on which ${\mb 1} \in \BC\subset \hat\fg_\kappa$ acts as the
identity. We will refer to objects of these category as modules at
level $\kappa$.

Let $\Gr_G = G((t))/G[[t]]$ be the affine Grassmannian of $G$. For
each $\kappa$ there is a category $\on{D}_\kappa(\Gr_G)\text{--}\on{mod}$
of $\kappa$-twisted right D-modules on $\Gr_G$ (see \cite{BD1}). We
have the functor of global sections $$\Gamma:
\on{D}_\kappa(\Gr_G)\text{--}\on{mod}\to \hat\fg_\kappa\text{--}\on{mod}, 
\qquad \F \mapsto \Gamma(\Gr_G,\F).$$

Let $\kappa_{Kil}$ be the Killing form, 
$\kappa_{Kil}(x,y) = \on{Tr}(\on{ad}_\fg(x)\circ \on{ad}_\fg(y))$. 
The level $\kappa_{crit} = -\frac{1}{2}\kappa_{Kil}$ is called
critical. A level $\kappa$ is called positive (resp., negative, irrational)
if $\kappa=c\cdot \kappa_{Kil}$ and $c+\frac{1}{2}\in \BQ^{>0}$
(resp., $c+\frac{1}{2}\in \BQ^{<0}$, $c\notin \BQ$).

It is known that the functor of global sections cannot be exact when
$\kappa$ is positive. In contrast, when $\kappa$ is negative or irrational, 
the functor $\Gamma$ is exact and faithful, as shown by A. Beilinson and
V. Drinfeld in \cite{BD1}, Theorem 7.15.8.  This statement is a
generalization for affine algebras of the famous theorem of
A. Beilinson and J. Bernstein, see \cite{BB}, that the functor of global
sections from the category of $\lambda$-twisted D-modules on the flag
variety $G/B$ is exact when $\lambda-\rho$ is anti-dominant and it is
faithful if $\lambda-\rho$ is, moreover, regular.

The purpose of this paper is to consider the functor of global
sections in the case of the critical level
$\kappa_{crit}$. (In what follows we will slightly abuse the notation 
and replace the subscript ${}_{\kappa_{crit}}$ simply by ${}_{crit}$.)
Unfortunately, it appears that the approach of
\cite{BD1} does not extend to the critical level case, so we have to
use other methods to analyze it. Our main result is that the
functor of global sections remains exact at the critical level:

\begin{thm}  \label{main}
The functor $\Gamma:\on{D}_{crit}(\Gr_G)\text{--}\on{mod}\to
\hat\fg_{crit}\text{--}\on{mod}$ is exact.
\end{thm}

In other words, we obtain that for any object $\F$ of
$\on{D}_{crit}(\Gr_G)\text{--}\on{mod}$ we have $\on{H}^i(\Gr_G,\F)=0$ for
$i>0$. Moreover, we will show that if $\F\neq 0$, then
$\on{H}^0(\Gr_G,\F)=\Gamma(\Gr_G,\F)\neq 0$. This property is
sometimes referred to as ``D-affineness'' of $\Gr_G$.

In fact, we will prove a stronger result. Namely, we note
after \cite{BD1}, that for a critically twisted D-module $\F$ on 
$\Gr_G$, the action of
$\hat\fg_{crit}$ on $\Gamma(\Gr_G,\F)$ extends to an action of the
{\it renormalized enveloping algebra} $U^{\ren}(\hat\fg_{crit})$ of
Sect. 5.6 of {\it loc.cit.} Following a conjecture and suggestion of 
Beilinson, we show that the resulting functor from
$\on{D}_{crit}(\Gr_G)\text{--}\on{mod}$ to the category of
$U^{\ren}(\hat\fg_{crit})$-modules is fully-faithful.

\ssec{}

Our method of proof of \thmref{main} uses the chiral algebra of
differential operators $\fD_{G,\kappa}$ introduced in \cite{AG}.
Modules over $\fD_{G,\kappa}$ should be viewed as (twisted) D-modules
on the loop group $G((t))$. In particular, the category of
$\kappa$-twisted D-modules on $\Gr_G$ is equivalent to the
subcategory in $\fD_{G,\kappa}\text{--}\on{mod}$, consisting of modules,
which are integrable with respect to the {\it right} action of
$\fg[[t]]$ (see \thmref{AG}). The functor of global sections on 
$\Gr_G$ corresponds, under this equivalence, to the functor of
$\fg[[t]]$-invariants. Therefore, we need to prove that this 
functor of invariants is exact.

This approach may be applied both when the level $\kappa$ is negative 
(or irrational) and critical. In the case of the negative or irrational
level the argument is considerably simpler, and so we obtain a proof 
of the exactness of $\Gamma$, which is different from that of \cite{BD1} 
(see \secref{negative level}).

The argument that we use for affine Kac-Moody algebras yields also a
different proof of the exactness statement from \cite{BB}. Let us
sketch this proof. For a weight $\lambda$, let
$\on{D}^\lambda(G/B)\text{--}\on{mod}$ be the category of left
$\lambda$-twisted D-modules on $G/B$ (here for an integral $\lambda$,
by a $\lambda$-twisted D-module on $G/B$ we understand a module over
the sheaf of differential operators acting on the line bundle $G
\times_B \la$). Let $\pi$ denote the natural projection $G\to G/B$,
and observe that the pull-back functor (in the sense of quasicoherent
sheaves) lifts to a functor $\pi^*:\on{D}^\lambda(G/B)\text{--}\on{mod}\to
\on{D}(G)\text{--}\on{mod}$. Furthermore, for a D-module $\F'$ on $G$, the
space of its global sections $\Gamma(G,\F')$ is naturally a bimodule
over $\fg$ due to the action of $G$ on itself by left and the right
translations.

For $\F\in \on{D}^\lambda(G/B)\text{--}\on{mod}$ we have
$$\Gamma(G/B,\F)\simeq
\on{Hom}_{\fb}\left(\BC^{-\lambda},\Gamma(G,\pi^*(\F))\right),$$ where
$\fb$ is the Borel subalgebra of $\fg$, $\BC^{-\lambda}$ its
one-dimensional representation corresponding to weight $-\lambda$, and
$\Gamma(G,\pi^*(\F))$ is a $\fb$-module via $\fb\hookrightarrow \fg$
and the right action of $\g$. But the $\fg$-module $\Gamma(G,\F')$,
where $\F' = \pi^*(\F)$ (with respect the right $\g$-action), belongs
to the category $\CO$.  Thus, we obtain a functor
$$\Gamma':\on{D}^\lambda(G/B)-\on{mod} \to \CO, \qquad \F\mapsto
\Gamma(G,\pi^*(\F)),$$ and
$$\Gamma(G/B,\F)\simeq
\on{Hom}_{\CO}\left(M(-\lambda),\Gamma'(\F)\right),$$
where $M(-\lambda)$ is the Verma module with highest weight $-\lambda$.

The functor $\Gamma'$ is exact because $G$ is affine, and it is
well-known that $M(\mu)$ is a projective object of $\CO$
precisely when $\mu+\rho$ is dominant. Hence $\Gamma$ is the
composition of two exact functors and, therefore, is itself exact.

\medskip

This reproves the Beilinson-Bernstein exactness statement. Note,
however, that the methods described above do not give the
non-vanishing assertion of \cite{BB}.

\ssec{}

The proof of the exactness result in the negative (or irrational) level case is
essentially a word for word repetition of the above argument, once we
are able to make sense of the category of D-modules on $G((t))$ as the
category of $\fD_{G,\kappa}$-modules. The key fact that we will use will be
the same: that the corresponding vacuum Weyl module $\BV_{\fg,\kappa'}$
is projective in the appropriate category $\CO$ if $\kappa'$ is
positive or irrational.

This argument does not work at the critical level, because in
this case the corresponding Weyl module $\BV_{\fg,crit}$ is far from
being projective in the category $\hat\fg_{crit}-\on{mod}$. Roughly,
the picture is as follows. Modules over $\hat\fg_{crit}$ give rise to
quasicoherent sheaves over the ind-scheme
$\on{Spec}(\FZ_{\fg,x})$, where $\FZ_{\fg,x}$ is the center of
the completed universal enveloping algebra of $\hat\fg_{crit}$ (this
is the ind-scheme of $^L \fg$-opers on the punctured disc, where $^L
\fg$ is the Langlands dual Lie algebra to $\g$). The ind-scheme
$\on{Spec}(\FZ_{\fg,x})$ contains a closed subscheme
$\on{Spec}(\fz_{\fg,x})$ (this is the scheme of $^L
\fg$-opers on the disc). The module $\BV_{\fg,crit}$ is supported on
$\on{Spec}(\fz_{\fg,x})$ and is projective in the category
of $\hat\fg_{crit}$-modules, which are supported on
$\on{Spec}(\fz_{\fg,x})$ and are
$G(\hat\CO_x)$-integrable. 

The problem is, however, that the
$\hat\fg_{crit}$-modules of the form
$\Gamma\left(G((t)),\pi^*(\F)\right)$, where $\pi$ is the projection
$G((t))\to G((t))/G[[t]]\simeq \Gr_G$, are never supported on
$\on{Spec}(\fz_{\fg,x})$. Therefore we need to show that
the functor of taking the maximal submodule of
$\Gamma\left(G((t)),\pi^*(\F)\right)$, which is supported on
$\on{Spec}(\fz_{\fg,x})$, is exact. We do that by showing
that the action of $\hat\fg_{crit}$ on
$\Gamma\left(G((t)),\pi^*(\F)\right)$ automatically extends to the
action of the renormalized chiral algebra $\CA^{\ren,\tau}_\fg$, which
is closely related to the renormalized enveloping algebra
$U^{\ren}(\hat\fg_{crit})$, mentioned above.

Consider the following analogy. Let $X$ be a
smooth variety and $Y$ its smooth closed subvariety. Then we have a
natural functor, denoted $i^!$, from the category of ${\mathcal O}_X$-modules, 
set-theoretically supported on $Y$, to the category of 
${\mathcal O}_Y$-modules: this functor takes an 
${\mathcal O}_X$-module $\F$ to its maximal submodule supported 
scheme-theoretically on $Y$. This is not an exact functor.
But the corresponding functor from the category
of right D-modules on $X$, also set-theoretically supported on $Y$,
to the category of right D-modules on $Y$ is exact, according to 
a basic theorem due to Kashiwara. 

In our situation the role of
the category of ${\mathcal O}_X$-modules is played by the category
$\hat\fg_{crit}-\on{mod}$, and the role of the category of D-modules
is played by the category of modules over the chiral algebra
$\CA^{\ren,\tau}_\fg$. We show that the above functor of taking the
maximal submodule of $\Gamma\left(G((t)),\pi^*(\F)\right)$, which is
supported on $\on{Spec}(\fz_{\fg,x})$, factors through the latter
category, and this allows us to prove the required exactness.

\ssec{Contents} Let us briefly describe how this paper is organized.
In \secref{negative level} we treat the negative level case. In
\secref{sect center} we recall some facts about commutative
D-algebras and the description of the center of the Kac-Moody 
chiral algebra at the critical level. 
In \secref{renormalized algebra} we discuss several versions of the
renormalized universal enveloping algebra at the critical level in the
setting of chiral algebras. In \secref{diff op} we study the chiral
algebra of differential operators $\fD_{G,\kappa}$ when
$\kappa=\kappa_{crit}$. In \secref{sections on Gr} we derive our main
\thmref{main} from two other statements, Theorems \ref{support on
tangent bundle} and \ref{Kashiwara}.  In \secref{proof of Kashiwara}
we prove \thmref{Kashiwara}, generalizing Kashiwara's theorem
about D-modules supported on a subvariety. In
\secref{other results} we prove \thmref{support on tangent bundle} and
describe the category of $\hat\fg_{crit}$-modules, which are supported on
$\on{Spec}(\fz_{\fg,x})$ and are $G(\hat\CO_x)$-integrable. Finally, in
\secref{faithfulness} we prove that the functor $\Gamma$ is faithful.  

\ssec{Conventions}

Our basic tool in this paper is the theory of chiral algebras. The
foundational work \cite{BD} on this subject will soon be published
(in our references we use the most recent version; a previous
one is currently available on the Web). In addition, an abridged 
summary of the results of \cite{BD} that are used in this paper may be found in
\cite{AG}. We wish to remark that all chiral algebras considered in
this paper are universal in the sense that they come from 
quasi-conformal vertex
algebras by a construction explained in \cite{FB}, Ch. 18. Therefore
all results of this paper may be easily rephrased in the language of
vertex algebras. We have chosen the language of chiral algebras in
order to be consistent with the language used in \cite{AG}.

We also use some the results from \cite{BD1}, which is still
unpublished, but available on the Web.

\medskip

The notation in this paper mainly follows that of \cite{AG}. Throughout
the paper, $X$ will be a fixed smooth curve; we will denote by $\CO_X$
(resp., $\omega_X$, $T_X$ and $D_X$) its structure sheaf (resp., the
sheaf of differentials, the tangent sheaf and the sheaf of
differential operators). 

We will work with D-modules on $X$, and in our notation we will not
distinguish between left and right D-modules, i.e., we will denote by
the same symbol a left D-module $\CM$ and the corresponding right
D-module $\CM\otimes \omega_X$.  The operations of tensor product,
taking symmetric algebra, and restriction to a subvariety must be
understood accordingly.

We will denote by $\Delta$ the diagonal embedding $X\to X\times X$,
and by $j$ the embedding of its complement $X\times X-\Delta(X)\to
X\times X$.  If $x\in X$ is a point, we will often consider D-modules
supported at $x$. In this case, our notation will not distinguish
between such a D-module and the underlying vector space.

We will use the notation $A \underset{B}\times C$ for a fiber
product of $A$ and $C$ over $B$, and the notation ${\mathcal P}
\times_G V$ for the twist of a $G$-module $V$ by a $G$-torsor
${\mathcal P}$.

Finally, if $\C$ is a category and $C$ is an object of 
$\C$, we will often write $C\in \C$.

\ssec{Acknowledgments}

D.G. would like to express his deep gratitude to A.~Beilinson 
for explaining to him the theory of chiral algebras, as well
as for numerous conversations related to this paper. He would
also like to thank S.~Arkhipov, J.~Bernstein for stimulating and 
helpful discussions.

In addition, both authors would like to thank A.~Beilinson for helpful 
remarks and suggestions and B.~Feigin for valuable discussions.

The research of E.F. was supported by grants from the Packard
foundation and the NSF. D.G. is a long-term prize fellow of the Clay
Mathematics Institute.

\section{The case of affine algebras at the negative and irrational levels}
\label{negative level}

\ssec{}
In this section we will show that the functor of global sections
$$\Gamma:\on{D}_{\kappa}(\Gr_G)-\text{mod}\to
\hat\fg_{\kappa}-\text{mod}$$ is exact when $\kappa$ is negative 
or irrational. A similar result has been proved by Beilinson and Drinfeld in \cite{BD1},
Theorem 7.15.8, by other methods. The setting of
\cite{BD1} is slightly different: they consider twisted D-modules on
the affine flag variety $\on{Fl}_G=G((t))/I$ instead of
$\Gr_G=G((t))/G[[t]]$, where $I\subset G[[t]]$ is the Iwahori
subgroup, i.e., the preimage of a fixed Borel subgroup $B \subset G$ under
the projection $G[[t]] \to G$. Here is the precise statement of
their theorem:

Recall that for any affine weight $\hat{\lambda}=(\lambda,2\check h\cdot c)$
(where $\lambda$ is a weight of $\fg$, $c\in \BC$ and $\check h$ is the dual
Coxeter number), we can consider
the corresponding category $\on{D}_{\hat{\lambda}}(\on{Fl}_G)\text{--}\on{mod}$,
of right $\hat{\lambda}$-twisted D-modules on $\on{Fl}_G$. A weight
$\hat{\lambda}$ is called anti-dominant if the corresponding Verma module
$M(\hat{\lambda})$ over $\hat\fg_\kappa$ (where $\kappa=c\cdot \kappa_{Kil}$) is
irreducible. According to a theorem of Kac and Kazhdan (see \cite{KK}),
this condition can be  expressed combinatorially as 
$\langle \hat{\lambda}+\rho_{aff},\check \alpha_{aff}\rangle\notin \BZ^{> 0}$,
where $\alpha_{aff}$ runs over the set of all positive affine coroots.
We have:
\begin{thm}   \label{on flags}
If $\hat{\lambda}$ is anti-dominant, then the functor
of sections $\Gamma:\on{D}_{\hat{\lambda}}(\on{Fl}_G)\text{--}\on{mod}\to
\hat\fg_\kappa\text{--}\on{mod}$ is exact.
\end{thm}

\thmref{on flags} formally implies the exactness statement on $\Gr_G$ 
(i.e., \thmref{exactness on negative} below)
only for $\kappa=c\cdot \kappa_{Kil}$ with $c$ either irrational, or
$c+\frac{1}{2}< -1+\frac{1}{2\check h}$; so our exactness result
is slightly sharper than that of \cite{BD1}. The proof 
of \thmref{exactness on negative} given below can be extended in a rather 
straightforward way to reprove \thmref{on flags}. In contrast, in the case of the critical level, 
it is essential that we consider D-modules on $\Gr_G$ and not on
$\on{Fl}_G$; in the latter case the naive analogue of the exactness
statement is not true.

Finally, note that Theorem 7.15.8 of \cite{BD1} contains also the
assertion that for $0\neq \F\in \on{D}_{\wt{\lambda}}(\on{Fl}_G)\text{--}\on{mod}$,
then the space of sections $\Gamma(\Gr_G, \F)$ is non-zero, implying
a similar statement for $\F\in \on{D}_{\kappa}(\Gr_G)$.
In \secref{faithfulness}, we will reprove this fact as well, by a
different method. This proof is the same in the negative and the
critical level cases.

\ssec{}

Thus, our goal in this section is to prove the following theorem:

\begin{thm} \label{exactness on negative}
The functor $\Gamma:\on{D}_{\kappa}(\Gr_G)\text{--}\on{mod}\to
\hat\fg_{\kappa}\text{--}\on{mod}$ is exact when $\kappa$ is
negative or irrational.
\end{thm}

The starting point of our proof is the following. Recall the chiral
algebra $\fD_{G,\kappa}$ (on our curve $X$), introduced in \cite{AG}.  
Let $\fD_{G,\kappa}\text{--}\on{mod}$ denote the category of chiral
$\fD_{G,\kappa}$-modules concentrated at a point $x\in X$. In
\cite{AG} it was shown that $\fD_{G,\kappa}\text{--}\on{mod}$ is a
substitute for the category of twisted D-modules on the loop group
$G((t))$, where $t$ is a formal coordinate on $X$ near $x$.

In particular, we have the forgetful functor 
$$\fD_{G,\kappa}\text{--}\on{mod}\to (\hat\fg_{\kappa}\times
\hat\fg_{2\kappa_{crit}-\kappa})\text{--}\on{mod},$$ where
$\hat\fg_{\kappa}\text{--}\on{mod}$ (resp.,
$\hat\fg_{2\kappa_{crit}-\kappa}\text{--}\on{mod}$) is the category
of representations of the affine algebra at the level $\kappa$ (resp.,
$2\kappa_{crit}-\kappa$).  This functor corresponds to the action of
the Lie algebra $\fg((t))$ on $G((t))$ by left and right
translations. In what follows, for a module $\CM\in 
\fD_{G,\kappa}\text{--}\on{mod}$, we will refer to the corresponding
actions of $\hat\fg_{\kappa}$ and $\hat\fg_{2\kappa_{crit}-\kappa}$
on it as ``left'' and ``right'', respectively.

Let $\wh{\CO}_x \simeq \BC[[t]]$ be the completed local ring at $x$.
Consider the subalgebra $\fg(\wh{\CO}_x)\subset
\hat\fg_{2\kappa_{crit}-\kappa}$. Let
$\hat\fg_{2\kappa_{crit}-\kappa}\text{--}\on{mod}^{G(\wh{\CO}_x)}$ be
the subcategory of $\hat\fg_{2\kappa_{crit}-\kappa}\text{--}\on{mod}$
whose objects are the $\hat\fg_{2\kappa_{crit}-\kappa}$-modules, on
which the action of $\fg(\wh{\CO}_x)$ may be exponentiated to an
action of the corresponding group $G(\wh{\CO}_x)$. Let
$\fD_{G,\kappa}\text{--}\on{mod}^{G(\wh{\CO}_x)}$ denote the full
subcategory of $\fD_{G,\kappa}\text{--}\on{mod}$ whose objects belong to
$\hat\fg_{2\kappa_{crit}-\kappa}\text{--}\on{mod}^{G(\wh{\CO}_x)}$
under the right action of
$\hat\fg_{2\kappa_{crit}-\kappa}\text{--}\on{mod}$.

The following result has been established in \cite{AG}:

\begin{thm} \label{AG}
There exists a canonical equivalence of categories
$$\on{D}_{\kappa}(\Gr_G)\text{--}\on{mod}\simeq
\fD_{G,\kappa}\text{--}\on{mod}^{G(\wh{\CO}_x)}.$$ If $\F$ is an
object of $\on{D}_{\kappa}(\Gr_G)\text{--}\on{mod}$, and $\CM_\F$ the
corresponding object of
$\fD_{G,\kappa}\text{--}\on{mod}^{G(\wh{\CO}_x)}$, then the
$\hat\fg_{\kappa}$-module $\Gamma(\Gr_G,\F)$ identifies with
$(\CM_\F)^{\fg(\wh{\CO}_x)}$, the space of invariants in $\CM_\F$
with respect to the Lie subalgebra $\fg(\wh{\CO}_x)\subset
\hat\fg_{2\kappa_{crit}-\kappa}$ under the right action.
\end{thm}

\ssec{}    \label{Weyl}

To prove the exactness of the functor
$\Gamma:\on{D}_{\kappa}(\Gr_G)\text{--}\on{mod}\to
\hat\fg_{\kappa}\text{--}\on{mod}$, for negative or irrational $\kappa$, 
we compose it with the tautological forgetful functor
$\hat\fg_{\kappa}\text{--}\on{mod}\to \on{Vect}$. By \thmref{AG}, this
composition can be rewritten as
$$\fD_{G,\kappa}\text{--}\on{mod}^{G(\wh{\CO}_x)}\to
\hat\fg_{2\kappa_{crit}-\kappa}\text{--}\on{mod}^{G(\wh{\CO}_x)} \to
\on{Vect},$$ where the first arrow is the forgetful functor, and the
second arrow is $\CM\mapsto \CM^{\fg(\wh{\CO}_x)}$.

For an arbitrary level $\kappa'$, let $\BV_{\fg,\kappa'}$ be the
vacuum Weyl module, i.e., $$\BV_{\fg,\kappa'}\simeq
\on{Ind}^{\hat\fg_{\kappa'}}_{\fg(\wh{\CO}_x)\oplus \BC {\mb
1}}(\BC),$$ where $\fg(\wh{\CO}_x)$ acts on $\BC$ by zero and ${\mb
1}$ acts as the identity. Tautologically, for any $\CM\in
\hat\fg_{\kappa'}\text{--}\on{mod}$, we have:
\begin{equation}    \label{hom}
\on{Hom}_{\hat\fg_{\kappa'}}(\BV_{\fg,\kappa'},\CM)\simeq
\CM^{\fg(\wh{\CO}_x)}.
\end{equation}
Moreover, $\BV_{\fg,\kappa'}$ is $G(\wh{\CO}_x)$-integrable,
i.e., belongs to $\hat\fg_{\kappa'}\text{--}\on{mod}^{G(\wh{\CO}_x)}$.

Observe that the condition that $\kappa$ is negative or irrational
is equivalent to $\kappa':=2\kappa_{crit}-\kappa$ being positive or irrational.
Therefore, to prove the exactness of $\Gamma$, it is enough to 
establish the following:

\begin{prop}    \label{noncrit}
If $\kappa'$ is positive or irrational, the module
$\BV_{\fg,\kappa'}$ is projective in
$\hat\fg_{\kappa'}\text{--}\on{mod}^{G(\wh{\CO}_x)}$.
\end{prop}

This proposition is well-known, and the proof is based on considering
eigenvalues of the Segal-Sugawara operator $L_0$. We include the proof for
completeness.

\begin{proof}
Recall that for every non-critical value of $\kappa'$, the vector space
underlying every object $\CM\in
\hat\fg_{\kappa'}\text{--}\on{mod}$ carries a canonical
endomorphism $L_0$ obtained via the Segal-Sugawara construction, such 
that the action
of $\hat\fg_{\kappa'}$ commutes with $L_0$ in the following way:
\begin{equation}    \label{comm rel}
[L_0,x\otimes t^n]=-n\cdot x\otimes t^n, \qquad x\in\fg, n \in \BZ.
\end{equation}
Explicitly,
let $\{x^a, x_a \}$ be bases in $\g$, dual with respect to
$\kappa_{Kil}$. The operator
\begin{equation}    \label{S0}
S_0 = \sum_a x^a\cdot x_a + 2\sum_a \sum_{n>0} 
x^a \otimes t^{-n}\cdot x_a\otimes t^n
\end{equation}
is well-defined on every object of $\hat\fg_{\kappa'}\text{--}\on{mod}$,
and it has the following commutation relation with elements of 
$\hat\fg_{\kappa'}$:
\begin{equation}   
[S_0,x \otimes t^n]=-(2c'+1)\cdot n\cdot x\otimes t^n, \qquad x\in\fg, 
n \in \BZ,
\end{equation}
where $c'$ is such that $\kappa'=c'\cdot \kappa_{Kil}$. Therefore, for
$c'\neq -\frac{1}{2}$, the operator
$L_0:=\frac{1}{2c'+1}\cdot S_0$ has the required properties.

\medskip

For an integral dominant weight $\la$ of $\fg$, let $V^\la$ be the
finite-dimensional irreducible $\fg$-module with highest weight $\la$
and $\BV^\lambda_{\fg,\kappa'}$ the corresponding Weyl module over
$\hat\fg_{\ka'}$,
$$\BV^\lambda_{\fg,\kappa'}=
\on{Ind}^{\hat\fg_{\kappa'}}_{\fg(\wh{\CO}_x)\oplus \BC {\mb
1}}(V^\lambda),$$ where $\fg(\wh{\CO}_x)$ acts on $V^\la$ through the
homomorphism $\fg(\wh{\CO}_x) \to \g$ and ${\mb 1}$ acts as the
identity. Then we find from formula \eqref{S0} that $L_0$ acts on the
subspace $V^\lambda \subset \BV^\lambda_{\fg,\kappa'}$ by the scalar
$\frac{C_\fg(\lambda)}{2c'+1}$, where $C_{\fg}(\lambda)$ is the
scalar by which the Casimir element $\sum_a x^a \cdot
x_a$ of $U(\fg)$ acts on $V^\lambda$. Note that $C_{\fg}(\lambda)$
is a non-negative rational number for any dominant integral weight $\lambda$,
and $C_{\fg}(\lambda)\neq 0$ if $\lambda\neq 0$.

Since $\BV^\lambda_{\fg,\kappa'}$ is generated from $V^\lambda$ by
the elements $x\otimes t^n\in \hat\fg_{\kappa'}$, $n<0$, we obtain that 
the action of $L_0$ on 
$\BV^\lambda_{\fg,\kappa'}$ is semi-simple. Moreover, since every
object $\CM\in\hat\fg_{\kappa'}\text{--}\on{mod}^{G(\wh{\CO}_x)}$ has a
filtration whose subquotients are quotients of the
$\BV^\lambda_{\fg,\kappa'}$'s, the action of $L_0$ on any such
$\CM$ is locally-finite.

\medskip

Suppose now that we have an extension
\begin{equation}    \label{ext}
0\to \CM\to \wt{\CM}\to \BV_{\fg,\kappa'}\to 0
\end{equation}
in $\hat\fg_{\kappa'}\text{--}\on{mod}^{G(\wh{\CO}_x)}$. Let
$\wt{v}^0\in \wt{\CM}$ be a lift to $\wt{\CM}$ of the generating vector
$v^0\in \BV_{\fg,\kappa'}$. Without loss of generality we may assume
that $\wt{v}^0$ has the same generalized eigenvalue as $v^0$, i.e.,
$0$, with respect to the action of $L_0$. It is sufficient to show
that then $\wt{v}^0$ belongs to $(\wt{\CM})^{\fg \otimes
t\BC[[t]]}$. Indeed, if this is so, then $\wt{v}^0$ is annihilated by
the entire Lie subalgebra $\fg(\wh{\CO}_x)$, due to the eigenvalue
condition, which would mean that the extension \eqref{ext} splits.

Suppose that this is not the case, i.e., that $\wt{v}^0$ is not
annihilated by $\fg \otimes t\BC[[t]]$. Then we can find a
sequence of elements $x_i\otimes t^{n_i} \in \fg \otimes t\BC[[t]]$, which we
can assume to be homogeneous, automatically of negative degrees 
with respect to $L_0$, such that the vector 
$$w=x_1\otimes t^{n_1}\cdot ...\cdot x_k\otimes t^{n_k}\cdot \wt{v}^0\in\CM$$
is non-zero and is annihilated by $\fg \otimes t\BC[[t]]$.
But then, on the one hand, the eigenvalue of $L_0$ on $w$ 
is $\on{deg} (x_1\otimes t^{n_1}) +...+\on{deg} (x_k\otimes t^{n_k})=
-(n_1+...+n_k) \in
\BZ^{<0}$, but on the other hand, it must be of the form
$\frac{C_\fg(\lambda)}{c'+\frac{1}{2}}$, which is not in $\BQ^{<0}$, by
our assumption. 

\end{proof}

\section{Center of the Kac-Moody chiral algebra at the critical level}  \label{sect center}

\ssec{}

Let $\CA$ be a unital chiral algebra on $X$. In what follows we will work with
a fixed point $x\in X$ and denote by $\CA\text{--}\on{mod}$ the category
of chiral $\CA$-modules, supported at $x$.

Recall that the center of $\CA$, denoted by $\fz(\CA)$, is by definition the 
maximal D-submodule of $\CA$ for, which the Lie-* bracket 
$\fz(\CA)\boxtimes \CA\to \Delta_!(\CA)$ 
vanishes. It is easy to see that $\fz(\CA)$ is a commutative chiral
subalgebra of $\CA$. For example, the unit $\omega_X\hookrightarrow \CA$
is always contained in $\fz(\CA)$.

\medskip

Let $\CA_{\fg,\kappa}$ be the chiral universal enveloping
algebra of the Lie-* algebra 
$L_{\fg,\kappa}:=\fg\otimes D_X\oplus \omega_X$ at the level $\kappa$
(modulo the relation equating the two embeddings of $\omega_X$).
We have the basic equivalence of categories:
$$\CA_{\fg,\kappa}\text{--}\on{mod}\simeq \hat\fg_{\kappa}\text{--}\on{mod}.$$

It is well-known that when $\kappa\neq \kappa_{crit}$, the inclusion
$\omega_X\to \fz(\CA_{\fg,\kappa})$ is an isomorphism. Let us denote by
$\fz_\fg$ the commutative chiral algebra $\fz(\CA_{\fg,crit})$. In
\thmref{FF isom} below we will recall the description of $\fz_\fg$
obtained in \cite{FF,F}.

\medskip

Let $\fz_{\fg,x}$ be the fiber of $\fz_\fg$ at $x$; this is
a commutative algebra. We have the natural maps
\begin{equation*} \label{center and endomorphisms}
\fz_{\fg,x} \longrightarrow
(\BV_{\fg,crit})^{\fg(\wh{\CO}_x)} \overset{\sim}\longleftarrow
\on{End}_{\hat\fg_{crit}}(\BV_{\fg,crit}),
\end{equation*}
where the left arrow is obtained from the definition of the center of a
chiral algebra, and the right arrow assigns to $e \in
\on{End}_{\hat\fg_{crit}}(\BV_{\fg,crit})$ the vector $e \cdot v^0$, where $v^0$ is the
canonical generator of $\BV_{\fg,crit}$. 

The resulting homomorphism of algebras $\fz_{\fg,x}\to 
\on{End}_{\hat\fg_{crit}}(\BV_{\fg,crit})$
is an isomorphism. In fact, for any chiral algebra $\CA$, its 
center $\fz(\CA)$ identifies with the D-module of endomorphisms of $\CA$ 
regarded as a chiral $\CA$-module. At the level of fibers, we have a map in one
direction $\fz(\CA)_x\to \on{End}_{\CA\text{--}\on{mod}}(\CA_x)$. This map is an isomorphism
if a certain flatness condition is satisfied. This condition is always satisfied if
$\CA$ is ``universal'', i.e., comes from a quasi-conformal vertex algebra, which
is the case of $\CA_{\fg,crit}$.

\ssec{}

For a chiral algebra $\CA$, let $\hat\CA_x$ be the canonical topological
associative algebra attached to the point $x$, see \cite{BD}, Sect. 3.6.2.
By definition, the category $\CA\text{--}\on{mod}$ 
endowed with the tautological forgetful functor to the category of
vector spaces, is equivalent to the category of discrete continuous
$\hat\CA_x$-modules, denoted $\hat\CA_x\text{--}\on{mod}$.

For example, when $\CA=\CA_{\fg,\kappa}$, the corresponding algebra 
$\hat\CA_{\fg,\kappa,x}$ identifies with the completed universal
enveloping algebra of $\hat\fg_\kappa$ modulo the relation
$1={\mathbf 1}$. We denote this algebra by $U'(\hat\fg_\kappa)$.

\medskip

When $\CA=\CB$ is commutative, the algebra $\hat\CB_x$ is commutative as well,
see \cite{BD}, Sects. 3.6.6 and 2.4.8.
In fact, $\hat\CB_x$ can be naturally represented as 
$\underset{\longleftarrow}{lim}\, \CB_x^i$, where $\CB^i$ 
are subalgebras of $\CB$, such that $\CB^i|_{X-x}\simeq \CB|_{X-x}$.
In particular, we have a surjective homomorphism $\hat\CB_x\to \CB_x$;
the subcategory $\CB_x\text{--}\on{mod}\subset \hat\CB_x\text{--}\on{mod}$ 
is the full subcategory of 
$\CB\text{--}\on{mod}$, whose objects are central $\CB$-modules, supported
at $x\in X$. (Recall that a $\CB$-module $\CM$ is called central if 
the action map $j_*j^*(\CB\boxtimes \CM)\to \Delta_!(\CB)$ comes from a map
$\CB\otimes\CM\to \CM$, i.e., 
factors through
$j_*j^*(\CB\boxtimes \CM)\twoheadrightarrow \Delta_!(\CB\otimes \CM)$.)

We will view $\on{Spec}(\hat\CB_x)$ as an ind-scheme 
$\underset{\longrightarrow}{lim}\, \on{Spec}(\CB_x^i)$; we have a closed
embedding $\on{Spec}(\CB_x)\hookrightarrow \on{Spec}(\hat\CB_x)$.

By taking $\CB=\fz_\fg$, we obtain a topological commutative algebra 
$\hat\fz_{\fg,x}$, which we will also denote by $\FZ_{\fg,x}$. The
corresponding map $\on{Spec}(\fz_{\fg,x})\hookrightarrow 
\on{Spec}(\FZ_{\fg,x})$ will be denoted by $\imath$.

\medskip

For any chiral algebra $\CA$ we have a homomorphism 
$$\widehat{\fz(A)}_x\to Z(\hat\CA_x),$$
where $Z(\hat\CA_x)$ is the center of $\hat\CA_x$. We do not know
whether this map is always an isomorphism, but can show that it
is an isomorphism for $\CA=\CA_{\fg,crit}$, using the description of 
$\fz_{\fg}$, given by \thmref{FF isom}(1) below (see \cite{BD1}, Theorem 3.7.7).
In other words, $\FZ_{\fg,x}$ maps isomorphically
to the center of $U'(\hat\fg_{crit})$.

\ssec{}  \label{FF description}

Let us recall the explicit description of $\fz_{\fg}$ and $\FZ_{\fg,x}$ due to
\cite{FF,F}. Let $^L G$ be the algebraic group of adjoint type whose
Lie algebra is the Langlands dual to $\g$. Denote by 
$\on{Op}_{^L G}({\mc D}_x)$ the affine scheme of $^L G$-{\em opers} 
on the disc 
${\mc D}_x = \on{Spec} (\wh{\mc O}_x)$. These are triples $({\mc F},{\mc
F}_B,\nabla)$, where ${\mc F}$ is a $^L G$--torsor over ${\mc D}_x$,
${\mc F}_B$ is its reduction to a fixed Borel subgroup $^L B \subset
{} ^L G$ and $\nabla$ is a connection on ${\mc F}$ (automatically
flat) such that ${\mc F}_B$ and $\nabla$ are in a special relative
position (see, e.g., \cite{F} for details). 

There exists an affine
$D_X$-scheme $J(\on{Op}_{^L G}(X))$ of jets of opers on $X$, whose
fiber at $x \in X$ is $\on{Op}_{^L G}({\mc D}_x)$ (see \cite{BD1},
Sect. 3.3.3), and so the corresponding sheaf of algebras of functions
$\on{Fun}\left(J(\on{Op}_{^L G}(X))\right)$ on $X$ is a commutative chiral
algebra. (In what follows, $\on{Fun}(\Y)$ stands for the ring of regular 
functions on a scheme $\Y$.)

The canonical topological algebra
associated to $\on{Fun}\left(J(\on{Op}_{^L G}(X))\right)$ at the point $x$
is nothing but the topological algebra of functions on the ind-affine
space $\on{Op}_{^L G}({\mc D}_x^\times)$ of $^L G$-opers on the
punctured disc ${\mc D}_x^\times = \on{Spec} (\wh{\mc K}_x)$, where
$\wh{\mc K}_x$ is the field of fractions of $\wh{\mc O}_x$. 
The following was established in \cite{FF,F}:

\begin{thm}    \label{FF isom}

\hfill

\smallskip

\noindent{\em (1)}
There exists a canonical isomorphism of $D_X$-algebras $$\fz_{\fg}
\simeq \on{Fun}\left(J(\on{Op}_{^L G}(X))\right).$$ 
In particular, we have an isomorphism of commutative algebras 
$\fz_{\fg,x} \simeq \on{Fun}\left(\on{Op}_{^L G}({\mc D}_x)\right)$
and of commutative topological algebras
$\FZ_{\fg,x} \simeq \on{Fun} \left(\on{Op}_{^L G}({\mc D}_x^\times)\right)$.

\smallskip

\noindent{\em (2)}
On the associated graded level, we have a commutative diagram of isomorphisms:
$$
\CD
\on{gr} (\fz_{\fg,x}) @<<<
\on{gr}\bigl(\on{Fun}\left(\on{Op}_{^L G}({\mc D}_x)\right)\bigr) \\
@VVV    @VVV   \\   
\on{Fun}\left(\left(\fg^*\times_{\BG_m}\Gamma({\mc D}_x,\Omega_X)
\right)^{G(\wh\CO_x)}\right)
@<<< \on{Fun}\left((^L\fg/^L G) \times_{\BG_m}\Gamma({\mc D}_x,\Omega_X)\right),
\endCD
$$
where $^L\fg/^L G=\on{Spec}(\on{Fun}(^L \fg)^{^L G})$.
\end{thm}

Note that in the lower left corner of the above commutative diagram we have
used the identification 
$\on{gr}(\BV_{\fg,crit}) \simeq \on{Sym}\left(\fg\otimes
(\wh\CK_x/\wh\CO_x)\right) \simeq
\on{Fun}\left(\fg^*\times_{\BG_m}\Gamma({\mc D}_x,\Omega_X) \right)$, and
$$\fg^*/G\simeq \fh^*/W\simeq {}^L \fh/W\simeq {}^L \fg/^L G.$$

\ssec{}   \label{comm alg}

To proceed we need to recall some more material from \cite{BD}
about commutative D-algebras (which, according to our conventions,
we do not distinguish them from commutative chiral algebras).

If $\CB$ is a commutative $D_X$-algebra, consider the $\CB$-module 
$\Omega^1(\CB)$ of (relative with respect to $X$) differentials on $\CB$, i.e.,
$\Omega^1(\CB)\simeq I_B/I^2_B$, where $I_B$ is the kernel of the
product $\CB\underset{\CO_X}\otimes \CB\to \CB$. From now on we will assume
that $\CB$ is finitely generated as a $D_X$-algebra; in this case
$\Omega^1(\CB)$ is finitely generated as a $\CB\otimes D_X$-module.

Recall that geometric points of the scheme $\on{Spec}(\CB_x)$ (resp., of the 
ind-scheme $\on{Spec}(\hat\CB_x)$) are the same as horizontal sections
of $\on{Spec}(\CB)$ over the formal disc $\D_x$ (resp., the formal
punctured disc $\D_x^\times$), see \cite{BD}, Sect. 2.4.9.
Let us explain the geometric meaning of $\Omega^1(\CB)$ in terms of these
identifications.

\medskip

Let $z$ be a point of $\on{Spec}(\CB_x)$, corresponding
to a horizontal section $\phi_z:\widehat\CO_x\to \CB_x$. Evidently,
we have: $\phi_z^*(\Omega^1(\CB))_x\simeq T^*_z(\on{Spec}(\CB_x))$,
where $T^*_z$ denotes the cotangent space at $z$.

From the definition of $\hat\CB_x$ we obtain a map
\begin{equation}  \label{map of cotangent}
H^0_{DR}\left(\D^\times_x,\phi_z^*(\Omega^1(\CB))\right)\to
T^*_z(\on{Spec}(\hat\CB_x)).
\end{equation}
(Since the D-module $\phi_z^*(\Omega^1(\CB))$ on $\D_x$ is finitely generated,
its de Rham cohomology over the formal and formal punctures disc makes
obvious sense.)
One can show that the map of \eqref{map of cotangent} is actually 
an isomorphism.

From the short exact sequence
\begin{equation}  \label{cotangent sequence}
0\to H^0_{DR}(\D_x,\phi_z^*(\Omega^1(\CB)))\to
H^0_{DR}(\D^\times_x,\phi_z^*(\Omega^1(\CB)))\to \phi_z^*(\Omega^1(\CB))_x\to 0,
\end{equation}
we obtain also an identification
$$H^0_{DR}(\D_x,\phi_z^*(\Omega^1(\CB)))\simeq N^*_{z}(\CB_x),$$
where $N_z^*(\CB_x)$ denotes the conormal 
to $\on{Spec}(\CB_x)$ inside $\on{Spec}(\hat\CB_x)$ at the point $z$.

\medskip

Assume now that $\CB$ is smooth (see \cite{BD}, Sect. 2.3.15
for the definition of smoothness). In this case $\Omega^1(\CB)$
is a finitely generated projective $\CB\otimes D_X$-module.

Consider the dual of $\Omega^1(\CB)$, i.e.,
$$\Theta(\CB):= 
\on{Hom}_{\CB\otimes D_X}\left(\Omega^1(\CB),\CB\otimes D_X\right).$$ 
This is a central $\CB$-module, called the tangent module to $\CB$.
Moreover, $\Theta(\CB)$ carries a canonical structure of Lie-*
algebroid over $\CB$ (see below). Evidently, $\Theta(\CB)$ is also
projective and finitely generated as a $\CB\otimes D_X$-module.

By dualizing the members of the short exact sequence
\eqref{cotangent sequence}, we obtain the identifications
(cf. \cite{BD}, Sect. 2.5.21):
\begin{align*}
&H^0_{DR}(\D_x,\phi_z^*(\Theta(\CB)))\simeq T_z(\on{Spec}(\CB_x)),\,\,
H^0_{DR}(\D^\times_x,\phi_z^*(\Theta(\CB)))\simeq T_z(\on{Spec}(\hat\CB_x)), \\
& \text{ \hskip3cm and } \phi_z^*(\Theta(\CB))_x\simeq N_z(\CB_x).
\end{align*}

\medskip

The next definition will be needed in \secref{sections on Gr}.
Let $\CI$ denote the kernel $\hat\CB_x\to \CB_x$. The quotient
$\CI/\CI^2$ is a topological module over $\CB_x$, and the normal
bundle, $N(\CB_x)$, to $\on{Spec}(\CB_x)$ inside $\on{Spec}(\hat\CB_x)$
can always be defined as the group ind-scheme 
$\on{Spec}(\on{Sym}_{\CB_x}(\CI/\CI^2))$. Let now $\CE\subset N(\CB_x)$
be a group ind-subscheme, and let $\CE^\perp$
be its annihilator in $\CI/\CI^2$.

We introduce the subcategory $\hat\CB_x\text{--}\on{mod}_\CE$ inside
the category $\hat\CB_x\text{--}\on{mod}$ of all chiral $\CB$-modules
supported at $x$ by imposing the following two conditions:

\noindent (1) We require that a module $\CM$, viewed as a
quasicoherent sheaf on $\on{Spec}(\hat\CB_x)$, is supported on the
formal neighborhood of $\on{Spec}(\CB_x)$.  In particular, $\CM$
acquires a canonical increasing filtration $\CM=\underset{i\geq 1}\cup\, 
\CM_i$, where $\CM_i\subset \CM$ is the submodule consisting
of sections annihilated by $\CI^i$.

\noindent (2) We require that the natural map
$\CI/\CI^2\underset{\CB_x}\otimes \CM_{i+1}/\CM_i\to
\CM_i/\CM_{i-1}$
vanish on $\CE^\perp\subset \CI/\CI^2$.

Note that the category $\hat\CB_x\text{--}\on{mod}_\CE$ is in general
not abelian.

\ssec{}  \label{notion of algebroid}

Let us now recall the notion of Lie-* algebroid over a 
commutative $D_X$ algebra $\CB$ (cf. \cite{BD}, Sect. 2.5).

Let $L$ be a central $\CB$-module. A structure of a {\em Lie-* algebroid over}
$\CB$ on $L$ is the data of a Lie-* bracket $L\boxtimes L\to
\Delta_!(L)$ and an action map $L\boxtimes \CB\to \Delta_!(\CB)$,
which satisfy the natural compatibility conditions given in \cite{BD},
Sect. 1.4.11 and 2.5.16. 

\medskip

If $\CB$ is smooth, then $\Theta(\CB)$ is well-defined, and it
carries a canonical structure of Lie-* 
algebroid over $\CB$. It is universal in the sense that for any Lie-* algebroid
$L$, its action on $\CB$ factors through a canonical map of Lie-* algebroids 
$\varpi:L\to \Theta(\CB)$, called the anchor map.

\medskip

Recall now that a structure on $\CB$ of chiral-Poisson (or, coisson,
in the terminology of \cite{BD}) algebra is a Lie-* bracket (called chiral-Poisson bracket)
$\CB\boxtimes \CB\to \Delta_!(\CB)$,
satisfying the Leibniz rule with respect to the multiplication on $\CB$
(cf. \cite{BD}, Sect. 1.4.18 and 2.6.). 

If $\CB$ is a chiral-Poisson algebra, $\Omega^1(\CB)$
acquires a unique structure of Lie-* algebroid, such that
the de Rham differential $d:\CB\to \Omega^1(\CB)$
is a map of Lie-* algebras, and the composition 
$$\CB\boxtimes \CB\overset{d\times \on{id}}\longrightarrow 
\Omega^1(\CB)\boxtimes \CB\to \Delta_!(\CB)$$
coincides with the chiral-Poisson bracket.

Following \cite{BD}, Sect. 2.6.6, we call a chiral-Poisson structure on
$\CB$ elliptic if (a) $\CB$ is smooth, (b) the anchor map
$\varpi:\Omega^1(\CB)\to \Theta(\CB)$ is injective, and (c)
$\on{coker}(\varpi)$ is a projective $\CB$-module of finite rank.

\ssec{}    \label{coisson}

Finally, let us recall the definition of the chiral-Poisson structure on $\fz_\fg$. 
Consider the flat $\BC[[\hslash]]$-family of chiral algebras $\CA_{\fg,\hslash}$,
corresponding to the pairing
$\kappa_\hslash=\kappa_{crit}+\hslash\cdot \kappa_0$, where
$\kappa_0$ is an arbitrary fixed non-zero invariant inner product.

For two sections $a,b\in \fz_\fg$, consider two arbitrary sections
$a_\hslash,b_\hslash\in \CA_{\fg,\hslash}$, whose values modulo
$\hslash$ are $a$ and $b$ respectively, and consider
$[a_\hslash,b_\hslash]\in \Delta_!(\CA_{\fg,\hslash})$. By assumption,
the last expression vanishes modulo $\hslash$. Therefore the section
$\frac{1}{\hslash}[a_\hslash,b_\hslash]\in \Delta_!(\CA_{\fg,\hslash})$
is well-defined. Moreover, its value mod $\hslash$ does not depend
on the choice of $a_\hslash$ and $b_\hslash$. 

Therefore we obtain a map $$a,b\in \fz_\fg \; \mapsto \;
\frac{1}{\hslash}[a_\hslash,b_\hslash] \; \on{mod} \; \hslash \in
\Delta_!(\CA_{\fg,crit}),$$ and it is easy to see that its image
belongs to $\Delta_!(\fz_\fg)$. Furthermore, it is straightforward to
verify that the resulting map $\fz_\fg\boxtimes \fz_\fg\to
\Delta_!(\fz_\fg)$ satisfies the axioms of the chiral-Poisson bracket,
see \cite{BD}, Sect. 2.7.1.

\medskip

Let us now describe in terms of
\thmref{FF isom} above the Lie-* algebroid $\Omega^1(\fz_\fg)$,
resulting from the chiral-Poisson structure on $\fz_\fg$.

First, recall from \cite{BD}, Sect. 2.4.11, that if $\CM$ is a central
module over a commutative chiral algebra $\CB$, then we can form
a topological module, denoted $\hat{h}{}^{\CB}_x(\CM)$ over $\CB_x$.
Applying this construction to $\CB=\fz_\fg$ and $\CM=\Omega^1(\fz_\fg)$
we obtain a topological Lie-* algebroid 
$\CG_{crit}:=\hat{h}{}^{\fz_\fg}_x(\Omega^1(\fz_\fg))$
(see \cite{BD}, Sect. 2.5.18 for details).

Let now $\fF_x$ be the universal $^L G$-torsor over $\on{Op}_{^L G}({\mc D}_x)$,
whose fiber over a given oper $({\mc F},{\mc F}_B,\nabla)\in \on{Op}_{^L G}({\mc D}_x)$
is the $^L G$-torsor of horizontal sections of $\F$, or, equivalently, the
fiber of $\F$ at $x\in \D_x$. Let us denote by $^L \CG_{\on{Op}}$ the
corresponding Atiyah algebroid over $\on{Op}_{^L G}({\mc D}_x)$, which by definition
consists of $^L G$-invariant vector fields on the total space of $\fF_x$.
We have a short exact sequence
$$0\to (^L \fg)_{\fF_x}\to {}^L \CG_{\on{Op}}\to T\left(\on{Op}_{^L G}({\mc D}_x)\right)\to 0,$$
where $T\left(\on{Op}_{^L G}({\mc D}_x)\right)$ denotes the tangent algebroid, and
$(^L \fg)_{\fF_x}$, which is the kernel of the anchor map, is the twist of the
adjoint representation by the $^L G$-torsor $\fF_x$.

The following was established in \cite{BD1}, Theorem 3.6.7 (see also \cite{BD}, Sect. 2.6.8),
using the fact that the isomorphism of $D_X$-algebras, given by \thmref{FF isom}(1), respects the
chiral-Poisson structures on both sides, where $J(\on{Op}_{^L G}(X))$ acquires a chiral-Poisson
structure  by its realization via the Drinfeld-Sokolov reduction.

\begin{thm}   \label{BD descr of Gelfand-Dikii}

\hfill

\smallskip

\noindent {\em (1)}
The chiral-Poisson structure on $\fz_\fg$ is elliptic.

\smallskip

\noindent {\em (2)}
Under the isomorphism $\on{Spec}(\fz_{\fg,x}) \simeq \on{Op}_{^L G}({\mc D}_x)$,
the algebroid $\CG_{crit}$ corresponds to the algebroid $^L \CG_{\on{Op}}$.
\end{thm}

\section{The renormalized chiral algebra}  \label{renormalized algebra}

\ssec{}

We will now refine the structure of chiral-Poisson algebra on
$\fz_\fg$ and obtain a chiral version of the renormalized 
universal enveloping algebra at the critical level introduced in
\cite{BD1}.

\medskip

First, we introduce a Lie-* algebra $\CA^\sharp_\fg$, which fits in a
short exact sequence
$$0\to \CA_{\fg,crit}\to \CA^\sharp_\fg\to \fz_\fg\to 0.$$

Namely, in the family of chiral algebras $\CA_{\fg,\hslash}$
consider the following subspace $\CA^\sharp_{\fg,\hslash}$, which
contains $\CA_{\fg,\hslash}$ and is contained in
$\frac{1}{\hslash}\cdot \CA_{\fg,\hslash}$:
$$
\CA^\sharp_{\fg,\hslash} = \{ \frac{a}{\hslash} \, | \, a \in
\CA_{\fg,\hslash}, a \, \on{mod}
\, \hslash \in \fz_\fg \}.
$$
Define $\CA^\sharp_\fg$ as $\CA^\sharp_{\fg,\hslash}/\hslash\cdot
\CA_{\fg,\hslash}$.  By repeating the construction of the
chiral-Poisson structure on $\fz_\fg$ from \secref{coisson}, we obtain
a Lie-* algebra structure on $\CA^\sharp_\fg$.

Note that the composition
$$\fz_\fg\boxtimes \CA^\sharp_\fg\to \CA^\sharp_\fg\boxtimes
\CA^\sharp_\fg\to \Delta_!(\CA^\sharp_\fg)$$ factors as
$\fz_\fg\boxtimes \CA^\sharp_\fg\to \fz_\fg\boxtimes \fz_\fg\to
\Delta_!(\fz_\fg)$, where the last arrow is chiral-Poisson bracket on
$\fz_\fg$.

\medskip

\begin{prop}  \label{BD algebroid}
There exist a unique Lie-* algebroid $\CA^{\flat}_\fg$ over $\fz_\fg$, which
fits into the following commutative diagram:
$$
\CD
0 @>>> \CA_{\fg,crit}/\fz_\fg @>>> \CA^{\flat}_\fg @>>>
\Omega^1(\fz_\fg) @>>> 0 \\
@. @AAA @AAA @AAA @. \\
0 @>>> \CA_{\fg,crit} @>>> \CA^\sharp_\fg @>>> \fz_\fg @>>> 0.
\endCD
$$
In the above diagram the rows are exact, and the rightmost vertical
map is the de Rham differential $\fz_\fg\to \Omega^1(\fz_\fg)$. 
\end{prop}

\begin{proof}

Recall the following general construction. Let $L$ be a Lie-* algebra
acting on a commutative chiral algebra $\CB$. Then we can form a
central $\CB$-module $\ol{\on{Ind}}_\CB(L):=\CB\otimes L$, which
will carry a natural structure of Lie-* algebroid over $\CB$. This is
analogous to the usual construction in differential geometry, when we
have a Lie-* algebra $\mathfrak l$ acting on a manifold $\Y$ and we
form the algebroid $\CO_\Y\otimes {\mathfrak l}$.

By taking $\CB=\fz_\fg$ and $L=\CA^\sharp_\fg/\fz_\fg$, we thus obtain a 
Lie-* algebroid $\ol{\on{Ind}}_{\fz_\fg}(\CA^\sharp_\fg/\fz_\fg)$
on $\fz_\fg$. We have a short exact sequence
$$0\to \fz_\fg\otimes (\CA_{\fg,crit}/\fz_\fg)\to
\ol{\on{Ind}}_{\fz_\fg}(\CA^\sharp_\fg/\fz_\fg)\to \fz_\fg\otimes
\fz_\fg\to 0.$$ To obtain from
$\ol{\on{Ind}}_{\fz_\fg}(\CA^\sharp_\fg/\fz_\fg)$ the desired
extension $\CA^{\flat}_\fg$, we need to take the quotient by two kinds
of relations.  First, we must pass from $\fz_\fg\otimes
(\CA_{\fg,crit}/\fz_\fg)$ to just $\CA_{\fg,crit}/\fz_\fg$, using the
structure of $\fz_\fg$-module on $\CA_{\fg,crit}$. Secondly, we must
impose the Leibniz rule to pass from the free $\fz_\fg$-module
$\fz_\fg\otimes \fz_\fg$ to $\Omega^1(\fz_\fg)$.
We will impose these two relations simultaneously. 

Consider the following three maps 
$\CA^\sharp_\fg\otimes \CA^\sharp_\fg\to \ol{\on{Ind}}_{\fz_\fg}(\CA^\sharp_\fg/\fz_\fg)$:

\medskip

\noindent (1) 
The first map is the projection
$\CA^\sharp_\fg\otimes \CA^\sharp_\fg\to
\fz_\fg\otimes (\CA^\sharp_\fg/\fz_\fg)$.

\medskip

\noindent (2) 
The second map is the projection
$\CA^\sharp_\fg\otimes \CA^\sharp_\fg\to 
(\CA^\sharp_\fg/\fz_\fg)\otimes \fz_\fg\simeq
\fz_\fg\otimes (\CA^\sharp_\fg/\fz_\fg)$.

\medskip

\noindent (3) To define the third map,  
note that chiral bracket on $\CA_{\fg,\hslash}$,
multiplied by $\hslash$, induces a map $j_*j^*(\CA^\sharp_\fg\boxtimes
\CA^\sharp_\fg)\to \Delta_!(\CA^\sharp_\fg)$. Composing
the latter with the projection $\Delta_!(\CA^\sharp_\fg)\to 
\Delta_!(\CA^\sharp_\fg/\fz_\fg)$, we obtain a map that vanishes
on $\CA^\sharp_\fg\boxtimes\CA^\sharp_\fg\subset j_*j^*(\CA^\sharp_\fg\boxtimes
\CA^\sharp_\fg)$, thereby giving rise to a map
$\CA^\sharp_\fg\otimes \CA^\sharp_\fg\to \fz_\fg\otimes (\CA^\sharp_\fg/\fz_\fg)$.

\medskip

By taking the linear combination of these three maps, namely
(1)--(2)--(3), we obtain a new map 
$\CA^\sharp_\fg\otimes \CA^\sharp_\fg\to \fz_\fg\otimes (\CA^\sharp_\fg/\fz_\fg)$.
We define $\CA^{\flat}_\fg$ as the quotient of 
$\ol{\on{Ind}}_{\fz_\fg}(\CA^\sharp_\fg/\fz_\fg)$ by the $\fz_\fg$-module, generated
by the image of the latter map.

One checks in a straightforward way that the Lie-* bracket
on $\ol{\on{Ind}}_{\fz_\fg}(\CA^\sharp_\fg/\fz_\fg)$ descends
to a Lie-* bracket on $\CA^{\flat}_\fg$, so that it becomes
a Lie-* algebroid over $\fz_\fg$. Moreover, by construction,
we have a short exact sequence

$$0\to (\CA_{\fg,crit}/\fz_\fg)'\to \CA^{\flat}_\fg\to \Omega^1(\fz_\fg)\to 0,$$
where $(\CA_{\fg,crit}/\fz_\fg)'$ is a quotient of $\CA_{\fg,crit}/\fz_\fg$. 
Let us show that $\CA_{\fg,crit}/\fz_\fg\to (\CA_{\fg,crit}/\fz_\fg)'$ is an isomorphism.

\medskip

Observe that the Lie-* algebra $\CA^\sharp_\fg$ acts on
$\CA_{\fg,crit}$. This action gives rise to
an action of the Lie-* algebroid
$\ol{\on{Ind}}_{\fz_\fg}(\CA^\sharp_\fg/\fz_\fg)$ on $\CA_{\fg,crit}$,
which is compatible with the $\fz_\fg$-module structure on
$\CA_{\fg,crit}$, and, moreover, it descends to an action of the
algebroid $\CA^{\flat}_\fg$ on $\CA_{\fg,crit}$.

The resulting Lie-* action of $\CA_{\fg,crit}$ on $\CA_{\fg,crit}$
obtained via 
$$\CA_{\fg,crit}\twoheadrightarrow (\CA_{\fg,crit}/\fz_\fg)'\hookrightarrow
\CA^{\flat}_\fg$$ coincides
with the initial Lie-* action of $\CA_{\fg,crit}$ on itself. By the definition
of the center, the kernel of the latter action is exactly $\fz_\fg$.

\end{proof}

\ssec{}   \label{flat modules}
Let us now introduce a category of modules over $\CA^{\flat}_\fg$,
which will be of interest for us.

First, note that the action of $\CA^\flat_\fg$ on $\CA_{\fg,crit}$,
introduced in the course of the proof of \propref{BD algebroid},
is compatible in the natural sense with the chiral bracket on
$\CA_{\fg,crit}$.

We define $\CA_{\fg}^{\flat}\text{--}\on{mod}$ to have as objects 
$\CM\in\CA_{\fg,crit}\text{--}\on{mod}$, endowed
with an additional action of the Lie-* algebroid $\CA^{\flat}_\fg$ 
(see \cite{BD}, Sect. 2.5.16 and 1.4.12 for the definition of the latter),
such that

\smallskip

\noindent (a) As a chiral module over $\fz_\fg$
(via $\fz_\fg\hookrightarrow \CA_{\fg,crit}$), $\CM$ is central.

\smallskip

\noindent (b) The two induced Lie-* actions
$(\CA_{\fg,crit}/\fz_\fg)\boxtimes \CM\to \Delta_!(\CM)$ 
(one coming from the $\CA^{\flat}_\fg$-action, and the other from
the $\CA_{\fg,crit}$-action and point (a) above) coincide.

\smallskip

\noindent (c) The chiral action of $\CA_{\fg,crit}$ and the Lie-*
action of $\CA^{\flat}_\fg$ on $\CM$ are compatible with the Lie-*
action of $\CA_{\fg}^{\flat}$ on $\CA_{\fg,crit}$.

\medskip

One can show that the category $\CA_{\fg}^{\flat}\text{--}\on{mod}$
is tautologically equivalent to the category of (discrete) modules
over the renormalized universal enveloping algebra introduced
in \cite{BD1}, Sect. 5.6.1.

For example, it is easy to see that if 
$\CM_\hslash$ is a flat $\BC[[\hslash]]$-family of chiral 
$\CA_{\fg,\hslash}$-modules such that the chiral $\fz_\fg$-module
$\CM:=\CM/\hslash\CM$ is central, then this $\CM$ is naturally an 
object of $\CA_{\fg}^{\flat}\text{--}\on{mod}$.

\ssec{}

In addition to the notion of a Lie-* algebroid there is also the
notion of a chiral Lie algebroid over a commutative chiral algebra
$\CB$, see \cite{BD}, Sect. 3.9.6. A Lie-* algebra $L$ is called a
{\em chiral Lie algebroid over} $\CB$ if we are given:

\smallskip

\noindent (1)
An action $L\boxtimes \CB\to\Delta_!(\CB)$
of $L$ as a Lie-* algebra on the commutative
chiral algebra $\CB$,

\smallskip

\noindent (2)
A chiral action $j_*j^*(\CB\boxtimes L)\to\Delta_!(L)$,
compatible with the action of $L$ on $\CB$ and the bracket
on $L$.

\smallskip

\noindent (3)
A map $\eta:\CB\to L$, compatible with both the $\CB$- and
$L$-actions,

\smallskip

such that the following conditions are satisfied:

\smallskip

\noindent (a)
The action in (1) is $\CB$-linear, in the sense that the two
natural maps $j_*j^*(\CB\boxtimes L)\boxtimes \CB\to \Delta_!(\CB)$
on $X^3$ coincide,

\smallskip

\noindent (b) 
The map $L\boxtimes \CB\overset{\on{id}\times \eta}\longrightarrow 
L\boxtimes L \to\Delta_!(L)$ equals the negative of 
$\CB\boxtimes L\hookrightarrow j_*j^*(\CB\boxtimes L)\to
\Delta_!(L)$,

\medskip

Note that if for a chiral Lie algebroid $L$ as above, the data 
of $\eta$ is zero, we retrieve the notion of Lie-* algebroid. In
most examples, however, the map $\eta$ is an injection. In this
case, the data of (1) is completely determined by (2) and (3),
and condition (a) is superfluous.

\medskip

It would be interesting to find out whether there exists a chiral Lie
algebroid $\CA^{\ren}_\fg$ over $\fz_\fg$, which is an extension
$$0\to \CA_{\fg,crit}\to \CA^{\ren}_\fg \to \Omega^1(\fz_\fg)\to 0,$$
such that the map $\CA^\sharp_\fg\to \CA^{\flat}_\fg$ lifts to a map
$\CA^\sharp_\fg\to \CA^{\ren}_\fg$.

However, we do not know how to construct such an
object. Instead, we will construct another chiral Lie algebroid
$\CA^{\ren,d}_\fg$, which is, in some sense, a double of
$\CA^{\ren}_\fg$. The construction of $\CA^{\ren,d}_\fg$ below is 
in terms of generators and relations. In the next section we will 
give a natural construction of $\CA^{\ren,d}_\fg$ via chiral
differential operators on the group $G$.

\medskip

Consider the Lie-* algebra $\CA^{\sharp,d}_\fg$ equal to
$$(\CA^\sharp_\fg\times \CA^\sharp_\fg) \underset{\fz_\fg\times
\fz_\fg}\times \fz_\fg,$$ where the map $\fz_\fg\to \fz_\fg\times
\fz_\fg$ is the anti-diagonal, i.e., $(\on{id},-\on{id})$. It fits
into a short exact sequence
$$0\to \CA_{\fg,crit}\times \CA_{\fg,crit}\to \CA^{\sharp,d}_\fg\to
\fz_\fg\to 0.$$

\begin{prop}  \label{doubled algebroid}
There exist a unique chiral algebroid $\CA^{\ren,d}_\fg$ over
$\fz_\fg$, which fits into the following commutative diagram with
exact rows:
$$
\CD
0 @>>> (\CA_{\fg,crit}\times\CA_{\fg,crit})/\fz_\fg
@>>> \CA^{\ren,d}_\fg @>>> \Omega^1(\fz_\fg) @>>> 0 \\
@. @AAA @AAA @AAA @. \\
0 @>>> \CA_{\fg,crit}\times\CA_{\fg,crit}
@>>> \CA^{\sharp,d}_\fg @>>> \fz_\fg @>>> 0,
\endCD
$$
where $\fz_\fg\hookrightarrow \CA_{\fg,crit}\times\CA_{\fg,crit}$
is the anti-diagonal embedding. 
\end{prop}

\begin{proof}

Let $L$ be a Lie-* algebra acting on a commutative 
chiral algebra $\CB$, as in the proof of \propref{BD algebroid}.
Then, following \cite{BD}, Sect. 3.9.9, one constructs a chiral Lie 
algebroid $\on{Ind}_\CB(L)$, which fits into a short exact sequence 
\begin{equation*} \label{induced algebroid}
0\to \CB\to \on{Ind}_\CB(L)\to \ol{\on{Ind}}_\CB(L)\to 0.
\end{equation*}

Indeed, consider the D-modules $j_*j^*(\CB\boxtimes L)$ 
and $\Delta_!(\CB)$ on $X\times X$. We have the maps
$$j_*j^*(\CB\boxtimes L) \leftarrow \CB\boxtimes L \rightarrow
\Delta_!(\CB),$$ where the left arrow is the natural inclusion, and the right
arrow is the negative of the Lie-* action.  Then the quotient
$\bigl(j_*j^*(\CB\boxtimes L)\oplus \Delta_!(\CB)\bigr)/\CB\boxtimes
L$ is supported on the diagonal, and therefore corresponds to a
D-module on $X$, which is by definition our $\on{Ind}_\CB(L)$. By construction, we
have the inclusions $\eta:\CB\to \on{Ind}_\CB(L)$ and $L\to
\on{Ind}_\CB(L)$, and a chiral action $j_*j^*(\CB\boxtimes L)\to
\Delta_!(\on{Ind}_\CB(L))$.  It is a straightforward verification to
show that these data extend uniquely to a Lie-* algebra structure on
$\on{Ind}_\CB(L)$ and a chiral action of $\CB$ on $\on{Ind}_\CB(L)$,
which satisfy the conditions of chiral Lie algebroid.

\medskip

Let us view $\CA^{\sharp,d}_\fg$ as a Lie-* algebra, which acts on
$\fz_\fg$ via $\CA^{\sharp,d}_\fg\to \fz_\fg$ and the chiral-Poisson
bracket on $\fz_\fg$. Consider the chiral Lie algebroid
$\on{Ind}_{\fz_\fg}(\CA^{\sharp,d}_\fg)$.  As in the case of
$\CA^\flat_\fg$ (see the proof of \propref{BD algebroid}), to obtain
from $\on{Ind}_{\fz_\fg}(\CA^{\sharp,d}_\fg)$ the desired
chiral algebroid $\CA^{\ren,d}_\fg$, we must take the quotient by some
additional relations.

The first set of relations is that we must identify the three copies of
$\fz_\fg$ inside $\on{Ind}_{\fz_\fg}(\CA^{\sharp,d}_\fg)$. One
copy is the image of the canonical embedding $\fz_\fg\hookrightarrow
\on{Ind}_{\fz_\fg}(\CA^{\sharp,d}_\fg)$ coming from the definition of
the induced algebroid. The other two copies come from $\fz_\fg\times
\fz_\fg\subset \CA^{\sharp,d}_\fg\hookrightarrow
\on{Ind}_{\fz_\fg}(\CA^{\sharp,d}_\fg)$. When we identify them, we
obtain a new chiral algebroid over $\fz_\fg$ which we denote by
$\on{Ind}'_{\fz_\fg}(\CA^{\sharp,d}_\fg)$.

The second set of relations is similar to what we had in the case of
$\CA^\flat_\fg$: they amount to killing the chiral $\fz_\fg$-submodule
generated by the image of a certain map $\CA^{\sharp,d}_\fg\otimes
\CA^{\sharp,d}_\fg\to \on{Ind}'_{\fz_\fg}(\CA^{\sharp,d}_\fg)$.  To
construct this map, we consider three morphisms from the D-module
$j_*j^*(\CA^{\sharp,d}_\fg\boxtimes \CA^{\sharp,d}_\fg)$ to
$\Delta_!(\on{Ind}_{\fz_\fg}(\CA^{\sharp,d}_\fg))$:

\noindent (1) 
The first map is $j_*j^*(\CA^{\sharp,d}_\fg\boxtimes \CA^{\sharp,d}_\fg)\to
j_*j^*(\fz_\fg\boxtimes \CA^{\sharp,d}_\fg)\to 
\Delta_!(\on{Ind}_{\fz_\fg}(\CA^{\sharp,d}_\fg))$,
where the first arrow comes from the natural projection
$\CA^{\sharp,d}_\fg\to \fz_\fg$.

\noindent (2) 
The second map is obtained from the first one by interchanging the roles of 
the factors in $j_*j^*(\CA^{\sharp,d}_\fg\boxtimes \CA^{\sharp,d}_\fg)$.

\noindent (3) To construct the third map, note that the chiral
bracket on $\CA_{\fg,\hslash}$ gives rise to a map
$$\hslash\cdot(\{\cdot,\cdot\},-\{\cdot,\cdot\}):
j_*j^*(\CA^{\sharp,d}_\fg\boxtimes \CA^{\sharp,d}_\fg)\to\Delta_!(\CA^{\sharp,d}_\fg),$$
and we compose it with the canonical map
$\Delta_!(\CA^{\sharp,d}_\fg)\to
\Delta_!(\on{Ind}_{\fz_\fg}(\CA^{\sharp,d}_\fg))$.

\medskip

Consider the linear combination (1)-(2)-(3) of the above three maps
as a new map from $j_*j^*(\CA^{\sharp,d}_\fg\boxtimes \CA^{\sharp,d}_\fg)$
to $\Delta_!(\on{Ind}_{\fz_\fg}(\CA^{\sharp,d}_\fg))$.
It is easy to see that the composition
$$\CA^{\sharp,d}_\fg\boxtimes \CA^{\sharp,d}_\fg\hookrightarrow
j_*j^*(\CA^{\sharp,d}_\fg\boxtimes \CA^{\sharp,d}_\fg)\to
\Delta_!(\on{Ind}_{\fz_\fg}(\CA^{\sharp,d}_\fg))\to 
\on{Ind}'_{\fz_\fg}(\CA^{\sharp,d}_\fg)$$
vanishes. Thus, we obtain the desired map
$\CA^{\sharp,d}_\fg\otimes
\CA^{\sharp,d}_\fg\to \on{Ind}'_{\fz_\fg}(\CA^{\sharp,d}_\fg)$.
We define $\CA^{\ren,d}_\fg$ as the quotient of 
$\on{Ind}'_{\fz_\fg}(\CA^{\sharp,d}_\fg)$ by the chiral
$\fz_\fg$-module, generated by the image of this map.

By construction, $\CA^{\ren,d}_\fg$ is a chiral $\fz_\fg$-module.  One
readily checks that the Lie-* bracket on
$\on{Ind}'_{\fz_\fg}(\CA^{\sharp,d}_\fg)$ descends to a Lie-* bracket
on $\CA^{\ren,d}_\fg$, such that, together with the map
$\fz_\fg\to \on{Ind}'_{\fz_\fg}(\CA^{\sharp,d}_\fg)\to
\CA^{\ren,d}_\fg$, these data define on $\CA^{\ren,d}_\fg$ a structure
of chiral Lie algebroid over $\fz_\fg$.

As in the case of $\CA^\flat_\fg$, we have a short exact sequence
$$0\to ((\CA_{\fg,crit}\times \CA_{\fg,crit})/\fz_\fg)'\to 
\CA^{\ren,d}_\fg\to \Omega^1(\fz_\fg)\to 0,$$
where $((\CA_{\fg,crit}\times \CA_{\fg,crit})/\fz_\fg)'$ is a certain
quotient of $(\CA_{\fg,crit}\times \CA_{\fg,crit})/\fz_\fg$. Let
us show that 
\begin{equation} \label{inj}
(\CA_{\fg,crit}\times \CA_{\fg,crit})/\fz_\fg\to
((\CA_{\fg,crit}\times \CA_{\fg,crit})/\fz_\fg)'
\end{equation}
is in fact an isomorphism.

Let $\CA^{\flat,d}_\fg$ be the Lie-* algebroid over $\fz_\fg$ equal to
$\CA^{\flat,d}_\fg/\fz_\fg$. We have a surjection
\begin{equation} \label{inj of ker anch}
(\CA_{\fg,crit}/\fz_\fg)\times (\CA_{\fg,crit}/\fz_\fg)\twoheadrightarrow
\on{ker}\left(\CA^{\flat,d}_\fg\to \Omega^1_{\fz_\fg}\right).
\end{equation}

As in the case of $\CA^\flat_\fg$, we show that $\CA^{\flat,d}_\fg$
acts naturally on the chiral algebra 
$\CA_{\fg,crit}\underset{\fz_\fg}\otimes \CA_{\fg,crit}$.
This implies that the map of \eqref{inj of ker anch} is an isomorphism.

Thus, it remains to show that the canonical map
$\fz_\fg\to \CA^{\ren,d}_\fg$ is injective. If it were not so,
the ideal $\on{ker}(\fz_\fg\to \CA^{\ren,d}_\fg)$ would be stable 
under the chiral-Poisson bracket on $\fz_\fg$. However, this is
impossible, since the above chiral-Poisson structure is elliptic
by \thmref{BD descr of Gelfand-Dikii}(1).

We remark that the isomorphism of \eqref{inj} can be alternatively
deduced from \thmref{embedding of algebroids} below.

\end{proof}

\ssec{}  \label{categories of modules}

Let $\CA_{\fg}^{\flat,d}$ be the Lie-* algebroid introduced 
in the proof of \propref{doubled algebroid}.

We introduce the category
$\CA_{\fg}^{\flat,d}\text{--}\on{mod}$ in a way analogous to
$\CA_{\fg}^\flat\text{--}\on{mod}$. Namely, the objects of 
$\CA_{\fg}^{\flat,d}\text{--}\on{mod}$ are modules over the
chiral algebra $\CA_{\fg,crit}\underset{\fz_\fg}\otimes \CA_{\fg,crit}$
(supported at $x\in X$) equipped with an extra Lie-* action of the 
Lie-* algebroid $\CA_{\fg}^{\flat,d}$ such that the conditions, analogous
to (a), (b) and (c) in the definition of $\CA_{\fg}^\flat\text{--}\on{mod}$,
hold.

\medskip

Next, we will introduce an appropriate category of chiral
modules over $\CA_{\fg}^{\ren,d}$.
First, recall from \cite{BD}, Sect. 3.9.24, the notion of chiral module
over a chiral Lie algebroid. 

If $L$ is a chiral algebroid over a commutative $D_X$-algebra $\CB$,
there exists a canonical chiral algebra $U(\CB,L)$, such the category
of chiral modules over $L$ (regarded as a chiral algebroid)
is equivalent to the category of chiral modules over $U(\CB,L)$
as a chiral algebra.

\medskip

Now let us introduce the category $\CA_{\fg}^{\ren,d}\text{--}\on{mod}$.
By definition, its objects are, as before, D-modules $\CM$ on $X$ 
supported at $x$ equipped with

\smallskip

\noindent (1) An action of the chiral algebra 
$\CA_{\fg,crit}\underset{\fz_\fg}\otimes \CA_{\fg,crit}$, 

\smallskip

\noindent (2) An action of the chiral algebroid
$\CA_{\fg}^{\ren,d}$, 

\medskip

\noindent such that the two induced chiral brackets
$j_*j^*\left((\CA_{\fg,crit}\times\CA_{\fg,crit}/\fz_\fg)
\boxtimes \CM\right)\to \Delta_!(\CM)$ coincide.

\medskip

Observe that one can reformulate the definition of
$\CA_{\fg}^{\ren,d}\text{--}\on{mod}$ as modules (supported at $x\in X$)
over a certain chiral algebra. Namely, let $U^{\ren,d}(L_{\fg,crit})$
be the quotient of the chiral algebra $U(\fz_\fg,\CA_{\fg}^{\ren,d})$
by the following relation:

\smallskip

We have a map
$U((\CA_{\fg,crit}\times\CA_{\fg,crit})/\fz_\fg)\to
U(\fz_\fg,\CA_{\fg}^{\ren,d})$ coming from the embedding of Lie-* algebras
$(\CA_{\fg,crit}\times\CA_{\fg,crit})/\fz_\fg\to \CA_{\fg}^{\ren,d}$.
In addition, we have a map $U((\CA_{\fg,crit}\times\CA_{\fg,crit})/\fz_\fg)\to
\CA_{\fg,crit}\underset{\fz_\fg}\otimes \CA_{\fg,crit}$. We need to
kill the ideal in $U(\fz_\fg,\CA_{\fg}^{\ren,d})$ generated by the image of
the kernel of the latter map.

\medskip

We define a PBW-type filtration on $U^{\ren,d}(L_{\fg,crit})$,
by setting $F^0\left(U^{\ren,d}(L_{\fg,crit})\right)$ to be the 
image of $\CA_{\fg,crit}\underset{\fz_\fg}\otimes \CA_{\fg,crit}$,
and by requiring inductively that 
$F^{i+1}\left(U^{\ren,d}(L_{\fg,crit})\right)$ is the smallest
$D_X$-submodule such that
$$j_*j^*\left((\CA^{\ren,d}_\fg\boxtimes 
F^i\left(U^{\ren,d}(L_{\fg,crit})\right)\right)\to
\Delta_!\left(F^{i+1}\left(U^{\ren,d}(L_{\fg,crit})\right)\right) \text{ and }$$
$$j_*j^*\left((\CA_{\fg,crit}\underset{\fz_\fg}\otimes \CA_{\fg,crit})\boxtimes
F^{i+1}\left(U^{\ren,d}(L_{\fg,crit})\right)\right)\to
\Delta_!\left(F^{i+1}\left(U^{\ren,d}(L_{\fg,crit})\right)\right).$$
In this case we automatically have also:
$$\CA^{\ren,d}_\fg\boxtimes 
F^i\left(U^{\ren,d}(L_{\fg,crit})\right)\to
\Delta_!\left(F^{i+1}\left(U^{\ren,d}(L_{\fg,crit})\right)\right).$$

We have a natural surjection on the associated graded level: 

\begin{equation}  \label{ass graded}
(\CA_{\fg,crit}\underset{\fz_\fg}\otimes \CA_{\fg,crit})
\underset{\fz_\fg}\otimes \on{Sym}_{\fz_\fg}(\Omega^1(\fz_\fg))
\twoheadrightarrow \on{gr}\left(U^{\ren,d}(L_{\fg,crit})\right).
\end{equation}

From \cite{BD}, Theorem 3.9.12 it follows that this map
is an isomorphism.

\ssec{}  \label{twist by tau}

Let now $\tau$ be an automorphism of $\fz_\fg$ as a chiral-Poisson algebra.
We can form the Lie-* algebra 
$$\CA^{\sharp,\tau}_\fg:=(\CA^\sharp_\fg\times \CA^\sharp_\fg)
\underset{\fz_\fg\times \fz_\fg}\times 
\fz_\fg,$$
where the map $\fz_\fg\to \fz_\fg\times \fz_\fg$
is now $(\on{id},-\tau)$.

Repeating the construction of \propref{doubled algebroid}, we obtain a
chiral algebroid $\CA^{\ren,\tau}_\fg$, which fits in a short exact
sequence
$$0\to (\CA_{\fg,crit}\times\CA_{\fg,crit})/\fz_\fg\to 
\CA^{\ren,\tau}_\fg\to \Omega^1(\fz_\fg)\to 0,$$
where $\fz_\fg$ is embedded into $\CA_{\fg,crit}\times\CA_{\fg,crit}$
also via $(\on{id},-\tau)$.

We will denote by $\CA^{\flat,\tau}_\fg$ the Lie-* algebroid on
$\fz_\fg$ equal to the quotient $\CA^{\ren,\tau}_\fg/\fz_\fg$.
Finally, in a way similar to the above, we introduce the corresponding 
categories of modules, $\CA_{\fg}^{\flat,\tau}\text{--}\on{mod}$ and 
$\CA_{\fg}^{\ren,\tau}\text{--}\on{mod}$, and the chiral algebra
$U^{\ren,\tau}(L_{\fg,crit})$.

\medskip

Note, however, that according to \cite{BD1} (and which something that we 
will have to use later), every automorphism of $\fz_\fg$, respecting the 
chiral-Poisson structure, comes from an outer automorphism of $\fg$.
This implies that, as abstract algebroids, $\CA^{ren,\tau}_\fg$ and
$\CA^{ren,d}_\fg$ are, in fact, isomorphic. 

\section{Chiral differential operators at the critical level}  \label{diff op}

\ssec{}

Recall that $\fD_{G,\kappa}$ denotes the chiral algebra of
differential operators on the group $G$ at level $\kappa$ (see
\cite{AG}), and $\fl,\fr$ are the two embeddings
$$\CA_{\fg,\kappa}\overset{\fl} \to
\fD_{G,\kappa}\overset{\fr}\leftarrow
\CA_{\fg,2\kappa_{crit}-\kappa}.$$

\medskip

Recall that if $\CM$ is a Lie-* module over a Lie-* algebra $L$, 
the centralizer of $L$ is the maximal D-submodule $\CM'\subset \CM$
such that the Lie-* bracket $L\boxtimes \CM'\to \Delta_!(\CM)$ vanishes.

\begin{lem}  \label{two centralizers}
The centralizer of $\fl(\CA_{\fg,\kappa})$ in $\fD_{G,\kappa}$
equals $\fr(\CA_{\fg,2\kappa_{crit}-\kappa})$. Conversely,
the centralizer of $\fr(\CA_{\fg,2\kappa_{crit}-\kappa})$
equals $\fl(\CA_{\fg,\kappa})$.
\end{lem}

\begin{proof}

The inclusion of $\fl(\CA_{\fg,\kappa})$ into the centralizer of
$\fr(\CA_{\fg,2\kappa_{crit}-\kappa})$ is just the fact that the
images of $\fl$ and $\fr$ Lie-* commute with each other.  The fact
that this inclusion is an equality is established as follows. Let
$\hat\fg_\kappa$ be the affine Kac-Moody algebra corresponding to a
point $x\in X$. The fiber $\CA_{\fg,\kappa,x}$ of $\CA_{\fg,\kappa}$
at $x$ is a $\hat\fg_\kappa$-module, equal to the vacuum module
$\BV_{\fg,\kappa}$.

Denote by $\fD_{G,\kappa,x}$ the fiber of $\fD_{G,\kappa}$ at $x$.
This is a module over $\hat\fg_\kappa\times
\hat{\fg}_{2\kappa_{crit}-\kappa}$. Recall that as
$\hat\fg_\kappa$-module, $\fD_{G,\kappa,x}$ is the induced module
$$\on{Ind}_{\fg(\wh{\CO}_x)\oplus \BC {\mb 1}}^{\hat{\fg}_\kappa}
\left(\on{Fun}\left(G(\wh{\CO}_x)\right) \right).$$ 

Moreover, the commuting right action of
$\fg(\wh{\CO}_x) \subset \hat{\fg}_{2\kappa_{crit}-\kappa}$ comes by
transport of structure from the right action of $\fg(\wh{\CO}_x)$ on
$\on{Fun}\left(G(\wh{\CO}_x)\right)$. In other words, as a right
$\fg(\wh{\CO}_x)$-module,
$$\on{Ind}_{\fg(\wh{\CO}_x)\oplus \BC {\mb 1}}^{\hat{\fg}_\kappa}
\left(\on{Fun} \left(G(\wh{\CO}_x)\right) \right)\simeq U(\fg \otimes
t^{-1}\BC[t^{-1}]) \otimes \on{Fun} \left(G(\wh{\CO}_x)\right),$$ where
$\fg(\wh{\CO}_x)$ acts through the second factor and $t$ is a
uniformizer in $\wh{\CO}_x$. At the level of fibers, the embedding
$\fl$ is just the natural embedding
$$\BV_{\fg,\kappa}\simeq \on{Ind}_{\fg(\wh{\CO}_x)\oplus \BC {\mb
1}}^{\hat{\fg}_\kappa} (\BC)\to \on{Ind}_{\fg(\wh{\CO}_x)\oplus \BC
{\mb 1}}^{\hat{\fg}_\kappa} \left(\on{Fun} \left(G(\wh{\CO}_x)\right) \right)$$
corresponding to the unit $\BC\to \on{Fun} \left(G(\wh{\CO}_x)\right)$. We have to
show that $\BV_{\fg,\kappa}\subset \fD_{G,\kappa,x}$ equals
$\left(\fD_{G,\kappa,x}\right)^{\fg(\wh{\CO}_x)}$, for
$\fg(\wh{\CO}_x)\subset \hat\fg_{2\kappa_{crit}-\kappa}$. But this
immediately follows from the above description of $\fD_{G,\kappa,x}$
as a $\fg(\wh{\CO}_x)$-module.

To finish the proof, observe that the roles of $\fl$ and $\fr$ in the
definition of $\fD_{G,\kappa}$ are symmetric, and in particular,
$\fD_{G,\kappa,x}$ is isomorphic to $\on{Ind}_{\fg(\wh{\CO}_x)\oplus
\BC}^{\hat{\fg}_{2\kappa_{crit}-\kappa}} \left(\on{Fun} \left(G(\wh{\CO}_x)\right)
\right)$ as a $\hat\fg(\wh{\CO}_x)\times
\hat{\fg}_{2\kappa_{crit}-\kappa}$-module. Indeed, we have a map
from the latter to the former, by the definition of the induction, and
this is map is clearly an isomorphism at the level of associate
graded spaces, by the PBW theorem.

\end{proof}

\medskip

Now we specialize to $\kappa=\kappa_{crit}$. Then $\fl$ and $\fr$ are
two different embeddings of $\CA_{\fg,crit}$ into
$\fD_{G,crit}$. \lemref{two centralizers} implies the following:

\begin{cor} \label{centers}
$\fl(\fz_\fg)=\fl(\CA_{\fg,crit})\cap \fr(\CA_{\fg,crit})=\fr(\fz_\fg)$.
\end{cor}

Let $\tau$ be the involution of the Dynkin diagram of $\fg$, which
sends a weight $\lambda$ to $-w_0(\lambda)$. We lift $\tau$ to an
outer automorphism of $\fg$, and it gives rise to a canonically
defined involution of $\fz_\fg$, which we will also denote by $\tau$.

\begin{thm}   \label{embedding of algebroids}
The two compositions $\fz_\fg\hookrightarrow
\CA_{\fg,crit}\overset{\fl}\to \fD_{G,crit}$ and
$\fz_\fg\hookrightarrow \CA_{\fg,crit}\overset{\fr}\to \fD_{G,crit}$
are intertwined by the automorphism $\tau:\fz_\fg\to \fz_\fg$.

We have an embedding of the chiral algebroid $\CA^{\ren,\tau}_\fg$
into $\fD_{G,crit}$ such that the maps $\fl$ and $\fr$ are the
compositions of this embedding and the canonical maps
$$\CA_{\fg,crit}\rightrightarrows (\CA_{\fg,crit}\times
\CA_{\fg,crit})/\fz_\fg \to \CA^{\ren,\tau}_\fg.$$ This embedding
extends to a homomorphism of chiral algebras
$U^{\ren,\tau}(L_{\fg,crit})\to \fD_{G,crit}$.
\end{thm}

The rest of this section is devoted to the proof of this theorem.

\ssec{}

The first step will be  to construct a map
$$\psi:\Omega^1(\fz_\fg)\to \fD_{G,crit}/
\bigl(\fl(\CA_{\fg,crit})+\fr(\CA_{\fg,crit})\bigr).$$

Note that if $\CM$ is a central module over a commutative chiral
algebra $\CB$, we have a naturally defined notion of derivation
$\CB\to \CM$, which amounts to a map of $\CB$-modules
$\Omega^1(\CB)\to \CM$.
We take $\CB=\fz_\fg$, and $\CM$ to be the centralizer of 
$\fz_\fg$ in the chiral $\fz_\fg$-module 
$\fD_{G,crit}/\left(\fl(\CA_{\fg,crit})+\fr(\CA_{\fg,crit})\right)$.
Thus, we need to construct a map
$$\fz_\fg\to \fD_{G,crit}/\left(\fl(\CA_{\fg,crit})+\fr(\CA_{\fg,crit})\right),$$
whose image Lie-* commutes with $\fz_\fg$, and which satisfies the Leibniz
rule.

\medskip

By letting the level $\kappa$ vary in the $\BC[[\hslash]]$-family $\kappa_{\hslash}$,
we obtain a flat $\BC[[\hslash]]$-family of chiral algebras $\fD_{G,\hslash}$. 
Note that the map $\fl$ extends to a map $\fl_\hslash:\CA_{\fg,\hslash}\to \fD_{G,\hslash}$, 
whereas the map $\fr$ gives rise to a map $\fr_\hslash:\CA_{\fg,-\hslash}\to 
\fD_{G,\hslash}$ (the negative appears due to the sign inversion in 
$\kappa\mapsto 2\kappa_{crit}-\kappa$).

Let $a$ be an element of $\fz_\fg$, and choose 
elements $a'_\hslash\in \CA_{\fg,\hslash}$ and 
$a''_{-\hslash}\in \CA_{\fg,-\hslash}$, which map to $a$
mod $\hslash$. Consider the element
$\fl_\hslash(a'_\hslash)-\fr_\hslash(a''_{-\hslash})\in \fD_{G,\hslash}$.
By definition, it vanishes mod $\hslash$; hence we obtain an element
$$\frac{\fl_\hslash(a'_\hslash)-\fr_\hslash(a''_{-\hslash})}{\hslash}
\; \on{mod} \; \hslash\in \fD_{G,crit}$$ which is well-defined modulo
$\fl(\CA_{\fg,crit})+\fr(\CA_{\fg,crit})$. This defines the required
map.  The fact that it is a derivation is a straightforward
verification.

\medskip

Note that the chiral bracket on $\fD_{G,crit}$ gives rise
to a well-defined Lie-* bracket
$$\fz_\fg\boxtimes 
\bigl(\fD_{G,crit}/\left(\fl(\CA_{\fg,crit})+\fr(\CA_{\fg,crit})\right)\bigr)\to
\Delta_!\left(\fD_{G,crit}\right),$$
where $\fz_\fg$ is thought of as embedded into $\fD_{G,crit}$ via
$\fz_\fg\hookrightarrow \CA_{\fg,crit}\overset{\fl}\to \fD_{G,crit}$.

\begin{lem}  \label{compatibility with Poisson}
The composition
$$\fz_\fg\boxtimes \Omega^1(\fz_\fg)\overset{\fl\times \psi}\rightarrow 
\fz_\fg\boxtimes 
\left(\fD_{G,crit}/\left(\fl(\CA_{\fg,crit})+\fr(\CA_{\fg,crit})\right)\right)\to
\Delta_!\left(\fD_{G,crit}\right)$$
factors as
$\fz_\fg\boxtimes \Omega^1(\fz_\fg)\to \Delta_!(\fz_\fg)\to
\Delta_!\left(\fD_{G,crit}\right)$,
where the first arrow is the chiral-Poisson structure on $\fz_\fg$.
A similar assertion holds for $\fz_\fg$ mapping to $\fD_{G,crit}$ via $\fr$.
\end{lem}

\begin{proof}

For two sections $a,b\in \fz_\fg$, and $a'_\hslash,a''_{-\hslash}$ as above,
we have
$$[\fl_\hslash(a'_\hslash)-
\fr_\hslash(a''_{-\hslash}),\fl_\hslash(b_\hslash)]=
[\fl_\hslash(a'_\hslash),\fl_\hslash(b_\hslash)]=\fl_\hslash([a'_\hslash,b_\hslash]),$$
because the images of $\fl_\hslash$ and $\fr_\hslash$ Lie-* commute in
$\fD_{G,\hslash}$. Hence, the assertion follows from the definition of
the chiral-Poisson structure on $\fz_\fg$.

\end{proof}

\ssec{} \label{jet construction} 
Since the images of $\fz_\fg$ in
$\fD_{G,crit}$ under $\fl$ and $\fr$ coincide, we obtain that there
exists an automorphism $\tau'$ of $\fz_\fg$, as a commutative chiral
algebra such that $\fl|_{\fz_\fg}=\fr|_{\fz_\fg}\circ \tau'$.  Our
goal now is to show that $\tau'=\tau$.

\lemref{compatibility with Poisson} implies that $\tau'$ is in fact an
automorphism of $\fz_\fg$ as a chiral-Poisson algebra. According to
Proposition 3.5.13 and Theorem 3.6.7 of \cite{BD1}, the chiral-Poisson
structure on $\fz_\fg$ is rigid, i.e., its group of automorphisms
equals the group of automorphisms of the Dynkin diagram of $\fg$.

Therefore, in order to prove that $\tau'=\tau$, it suffices to show,
that the two automorphisms coincide at the associate graded
level. Recall that if $\C$ is a commutative $\CO_X$-algebra, $\CJ(\C)$
denotes the corresponding commutative chiral algebra, obtained by
the jet construction from $\C$ (see \cite{BD}, Sect. 2.3.2).
Recall that the chiral algebras $\CA_{\fg,\kappa}$ and $\fD_{G,crit}$
are naturally filtered (see \cite{BD}, Sect. 3.7.13 and 3.9.11), and
we have:
$$\on{gr}(\CA_{\fg,\kappa}) \simeq 
\CJ\left(\on{Sym}(\fg\otimes \omega^{\otimes -1}_X)\right)
\simeq \CJ\left(\on{Fun}\left( \fg^* \times_{\BG_m} \omega_X\right) \right)$$
and 
$$\on{gr}(\fD_{G,\kappa})\simeq \CJ(\CO_G\otimes \CO_X)\otimes
\CJ\left(\on{Sym}(\fg\otimes \omega^{\otimes -1}_X)\right)\simeq
\CJ\left(\on{Fun}\left( T^*G \times_{\BG_m} \omega_X \right)\right),$$ so
that the maps $\on{gr}(\fl)$ and $\on{gr}(\fr)$ come from the (moment)
maps $T^*G\rightrightarrows \fg^*$ corresponding to the action of
$\fg$ on $G$ by left and right translation, respectively.

Moreover,
$$\on{gr}(\fz_\fg)\hookrightarrow \CJ\left(\on{Sym}(\fg\otimes
\omega^{\otimes -1}_X)^G\right) \simeq \CJ\left(\on{Fun}\left( \fg^*/G
\times_{\BG_m} \omega_X\right) \right).$$ 
(This inclusion is, in fact, an equality, by \thmref{FF isom}(2).) 
Therefore, the required assertion follows from the fact that
the two maps $T^*G\rightrightarrows \fg^*\to \fg^*/G$ differ by the
automorphism $\tau$.

\ssec{}

To finish the proof of \thmref{embedding of algebroids}, we will
identify $\CA^{\ren,\tau}_\fg$ with $$\Omega^1(\fz_\fg)
\underset{\fD_{G,crit}/\left(\fl(\CA_{\fg,crit}) +
\fr(\CA_{\fg,crit})\right)}\times \fD_{G,crit}.$$

\medskip

Note that the construction of the map $\psi$ gives in fact a map
$\CA^{\sharp,\tau}_\fg/\fz_\fg\to \fD_{G,crit}$. Indeed, a section of
$\CA^{\sharp,\tau}_\fg$ has a form
$\frac{(a'_\hslash,a''_{-\hslash})}{\hslash}$ for $a'_\hslash\in
\CA_{\fg,\hslash}$, $a''_{-\hslash}\in \CA_{\fg,-\hslash}$, such that
$$a'\, mod\, \hslash=-\tau(a'') \; \on{mod} \; \hslash\in \fz_\fg.$$
We associate to it a section of $\fD_{G,crit}$ equal to
$\frac{\fl_\hslash(a'_\hslash)+\fr_\hslash(a''_{-\hslash})}{\hslash}$.

In addition, $\fD_{G,crit}$ is obviously a chiral $\fz_\fg$-module, so
we obtain a map $\on{Ind}_{\fz_\fg}(\CA^{\sharp,\tau}_\fg)\to
\fD_{G,crit}$, and it is straightforward to check that the relations,
defining $\CA^{\ren,\tau}_\fg$ as a quotient of
$\on{Ind}_{\fz_\fg}(\CA^{\sharp,\tau}_\fg)$, hold.

Finally, we obtain a homomorphism of chiral algebras
$U(\fz_\fg,\CA^{\ren,\tau}_\fg)\to \fD_{G,crit}$, and it is easy to see that
it annihilates the ideal defining $U^{\ren,\tau}(L_{\fg,crit})$ as a
quotient of $U(\fz_\fg,\CA^{\ren,\tau}_\fg)$.

\section{The functor of global
sections on the affine Grassmannian} \label{sections on Gr}

\ssec{}

Let $\fD_{G,crit}\text{--}\on{mod}^{G(\wh{\CO}_x)}$ be the category of
chiral $\fD_{G,crit}$-modules supported at the point $x\in X$,
which are $G(\wh{\CO}_x)$-integrable with respect to the
embedding $\fr:\CA_{\fg,crit}\to \fD_{G,crit}$.

Let $\F$ be a critically twisted D-module on $\Gr_G$, and
$\CM_\F$--the corresponding object of
$\fD_{G,crit}\text{--}\on{mod}^{G(\wh{\CO}_x)}$.  According to \thmref{AG},
\begin{equation}    \label{relate}
\Gamma(\Gr_G,\F)\simeq \on{Hom}_{\fg(\wh{\CO}_x)}(\BC,\CM_\F)\simeq
\on{Hom}_{\hat\fg_{crit}}(\BV_{\fg,crit},\CM_\F),
\end{equation}
where $\CM_\F$ is regarded as a $\hat\fg_{crit}$-module via $\fr$,
and $\BV_{\fg,crit}\simeq
\on{Ind}^{\hat\fg_{crit}}_{\fg(\wh{\CO}_x)\oplus \BC {\mb 1}}(\BC)$
is the vacuum module, i.e., the fiber $\CA_{\fg,crit,x}$ of
$\CA_{\fg,crit}$ at $x$.

\medskip

Recall that $\hat\fg_{crit}\text{--}\on{mod}$ denotes the category of all discrete
$\hat\fg_{crit}$-modules supported at $x\in X$, and let
$\hat\fg_{crit}\text{--}\on{mod}^{G(\wh{\CO}_x)}$ be the subcategory of
$G(\wh{\CO}_x)$-integrable modules. Obviously, $\BV_{\fg,crit}$
belongs to $\hat\fg_{crit}\text{--}\on{mod}^{G(\wh{\CO}_x)}$, but the main
difficulty in the proof of \thmref{main} is that, in contrast to the
negative or irrational level cases, $\BV_{\fg,crit}$ is not projective in this
category.

\medskip

Let now $\hat\fg_{crit}\text{--}\on{mod}_{\reg}$ denote the subcategory of
$\hat\fg_{crit}\text{--}\on{mod}$ consisting of modules, which are
central (cf. \cite{BD}, Sect. 3.3.7)
with respect to the action of $\fz_\fg$. Let us denote by
$\hat\fg_{crit}\text{--}\on{mod}^{G(\wh{\CO}_x)}_{\reg}$ the intersection
$\hat\fg_{crit}\text{--}\on{mod}_{\reg}\cap
\hat\fg_{crit}\text{--}\on{mod}^{G(\wh{\CO}_x)}$. The module
$\BV_{\fg,crit}$ belongs to 
$\hat\fg_{crit}\text{--}\on{mod}^{G(\wh{\CO}_x)}_{\reg}$, but 
the modules from $\fD_{G,crit}\text{--}\on{mod}^{G(\wh{\CO}_x)}$, regarded
as objects of $\hat\fg_{crit}\text{--}\on{mod}^{G(\wh{\CO}_x)}$, do not
belong there.

The following projectivity result is essentially due to \cite{BD1}
(see \secref{other results} for the proof).

\begin{thm} \label{drinfeld}
The module $\BV_{\fg,crit}$ is a projective generator of the category
$\hat\fg_{crit}\text{--}\on{mod}^{G(\wh{\CO}_x)}_{\reg}$. In particular,
the functor
$\hat\fg_{crit}-\on{mod}^{G(\wh{\CO}_x)}_{\reg} \to \on{Vect}$ given by 
$\CM \mapsto \on{Hom}_{\hat\fg_{crit}}(\BV_{\fg,crit},\CM)$,
is exact.
\end{thm}

\medskip

Consider the functors
\begin{align*}
\sF: \hat\fg_{crit}\text{--}\on{mod}^{G(\wh{\CO}_x)}_{\reg} &\to
\fz_{\fg,x}\text{--}\on{mod}, \qquad {\mc M} \mapsto
\Hom_{\hat\fg_{crit}}(\BV_{\fg,crit},\CM), \\ \sG:
\fz_{\fg,x}\text{--}\on{mod} &\to
\hat\fg_{crit}\text{--}\on{mod}^{G(\wh{\CO}_x)}_{\reg}, \qquad {\mc F}
\mapsto \BV_{\fg,crit} \underset{\fz_{\fg,x}}\otimes {\mc F}.
\end{align*}

Now \thmref{drinfeld} implies the following:

\begin{thm} \label{Drinfeld's equivalence of categories}
The functors $\sF$ and $\sG$ are mutually inverse equivalences of
categories.
\end{thm}

By combining this theorem with \thmref{FF isom},
we obtain that the category 
$\hat\fg_{crit}\text{--}\on{mod}^{G(\wh{\CO}_x)}_{\reg}$ is equivalent to
the category of quasicoherent sheaves on the scheme $\on{Op}_{^L G}({\mc D}_x)$.

\ssec{}  \label{functors}

Consider the functor $\imath^!:\FZ_{\fg,x}\text{--}\on{mod}\longrightarrow
\fz_{\fg,x}\text{--}\on{mod}$, which takes a $\FZ_{\fg,x}$-module to its
maximal submodule, scheme-theoretically supported on
$\on{Spec}(\fz_{\fg,x})$, i.e., for an object $\CM\in \FZ_{\fg,x}\text{--}\on{mod}$,
$\imath^!(\CM)$ consists of elements annihilated by $\on{ker}(\FZ_{\fg,x}\to \fz_{\fg,x})$.
We will denote by the same symbol $\imath^!$ the corresponding functors
$$
\hat\fg_{crit}\text{--}\on{mod} \to \hat\fg_{crit}\text{--}\on{mod}_{\reg} \qquad
\on{and} \qquad \hat\fg_{crit}\text{--}\on{mod}^{G(\wh{\CO}_x)} \to
\hat\fg_{crit}\text{--}\on{mod}^{G(\wh{\CO}_x)}_{\reg}.
$$

According to formula \eqref{relate}, the functor of global sections
$\Gamma: \on{D}_{crit}(\Gr_G)\text{--}\on{mod} \to \on{Vect}$ can be
viewed as a functor $\fD_{G,crit}\text{--}\on{mod}^{G(\wh{\CO}_x)}\to
\on{Vect}$ given by $\CM\mapsto
\on{Hom}_{\hat\fg_{crit}}(\BV_{\fg,crit},\CM)$. Since
$\BV_{\fg,crit}$ is supported on $\on{Spec}(\fz_{\fg,x})$, we obtain
that this functor factors as
$$\CM\mapsto \imath^!(\CM)\mapsto 
\on{Hom}_{\hat\fg_{crit}\text{--}\on{mod}^{G(\wh{\CO}_x)}_{\reg}}(\BV_{\fg,crit},\imath^!(\CM)).$$
But according to \thmref{drinfeld}, the second functor is exact.
Therefore \thmref{main} is equivalent to the following:

\begin{thm} \label{critical reformulation}
The composition
$$\fD_{G,crit}\text{--}\on{mod}^{G(\wh{\CO}_x)}\to 
\hat\fg_{crit}\text{--}\on{mod}^{G(\wh{\CO}_x)}
\overset{\imath^!}\longrightarrow \hat\fg_{crit}
\text{--}\on{mod}^{G(\wh{\CO}_x)}_{\reg},$$
where the first arrow is the forgetful functor corresponding to the
embedding $\fr$, is exact.
\end{thm}

\ssec{}

Let $\FZ_{\fg,x}\text{--}\on{mod}_{\nil}$ be the full subcategory of
$\FZ_{\fg,x}\text{--}\on{mod}$, whose objects are those modules, which are
set-theoretically supported on $\on{Spec}(\fz_{\fg,x})\subset
\on{Spec}(\FZ_{\fg,x})$, i.e., modules 
supported on the formal neighborhood of $\on{Spec}(\fz_{\fg,x})$. 
Let $\hat\fg_{crit}\text{--}\on{mod}_{\nil}$,
(resp., $\hat\fg_{crit}\text{--}\on{mod}_{\nil}^{G(\wh{\CO}_x)}$)
denote the corresponding full subcategory of $\hat\fg_{crit}\text{--}\on{mod}$ 
(resp., $\hat\fg_{crit}\text{--}\on{mod}^{G(\wh{\CO}_x)}$).

Let $\wt{\imath}{}^!:\FZ_{\fg,x}\text{--}\on{mod}\longrightarrow
\FZ_{\fg,x}\text{--}\on{mod}_{\nil}$ be the functor that attaches to a
$\FZ_{\fg,x}$-module its maximal submodule, which is supported on the
formal neighborhood of $\on{Spec}(\fz_{\fg,x})$.  In other words, for
$\CM\in \FZ_{\fg,x}\text{--}\on{mod}$, $\wt{\imath}{}^!(\CM)$ consists of
all sections annihilated by some power of the ideal of
$\on{Spec}(\fz_{\fg,x})$ in $\on{Spec}(\FZ_{\fg,x})$.  We will denote
in the same way the corresponding functors
$$
\hat\fg_{crit}\text{--}\on{mod} \to \hat\fg_{crit}\text{--}\on{mod}_{\nil} \qquad
\on{and} \qquad \hat\fg_{crit}\text{--}\on{mod}^{G(\wh{\CO}_x)} \to
\hat\fg_{crit}\text{--}\on{mod}_{\nil}^{G(\wh{\CO}_x)}.
$$
Clearly, $\imath^!\simeq \imath^!\circ \wt{\imath}{}^!$.

\begin{prop}   \label{restriction to formal}
Every object $\CM\in\hat\fg_{crit}\text{--}\on{mod}^{G(\wh{\CO}_x)}$
can be canonically decomposed as a direct sum
$\CM=\CM^{\nil}\oplus \CM^{\on{non-reg}}$, where
$\wt{\imath}{}^!(\CM^{\nil})\simeq \CM^{\nil}$ and
$\wt{\imath}{}^!(\CM^{\on{non-reg}})$ is supported away from
$\on{Spec}(\fz_{\fg,x})\subset \on{Spec}(\FZ_{\fg,x})$.
\end{prop}

\begin{cor}  \label{exactness of restr to formal}
The functor 
$\wt{\imath}{}^!:\hat\fg_{crit}\text{--}\on{mod}^{G(\wh{\CO}_x)} \to
\hat\fg_{crit}\text{--}\on{mod}^{G(\wh{\CO}_x)}_{\nil}$ is exact.
\end{cor}

\begin{proof}(of \propref{restriction to formal})

For an irreducible $\fg$-module $V^\lambda$ with a dominant highest
weight $\lambda \in \Lambda^+$, let $\BV^\lambda_{\fg,crit}$ be the
corresponding Weyl module in
$\hat\fg_{crit}\text{--}\on{mod}^{G(\wh{\CO}_x)}$, as defined in
\secref{Weyl}; in particular, $\BV_{\fg,crit}=\BV^0_{\fg,crit}$. Let
$\Y^\lambda\subset \on{Spec}(\FZ_{\fg,x})$ be the closed sub ind-scheme
corresponding to the annihilating ideal of $\BV^\lambda_{\fg,crit}$ in
$\FZ_{\fg,x}$. In particular, $\Y^0 = \on{Spec}(\fz_{\fg,x})$.

By definition, every object in the category
$\hat\fg_{crit}\text{--}\on{mod}^{G(\wh{\CO}_x)}$ has a filtration whose
successive quotients are generated by vectors, on which the subalgebra
$\fg\otimes t\BC[[t]]\subset \fg(\wh{\CO}_x)$ acts trivially. In particular,
such a subquotient is a quotient of $\BV^\lambda_{\fg,crit}$ for some
$\lambda$. Therefore, the support in $\on{Spec}(\FZ_{\fg,x})$ of
every object from $\hat\fg_{crit}\text{--}\on{mod}^{G(\wh{\CO}_x)}$ is
contained in the union of the formal neighborhoods of $\Y^\lambda$ for
$\lambda\in \Lambda^+$.

\begin{lem}   \label{Sugawara calculation}
For $\lambda\neq 0$, $\Y^\lambda \cap \on{Spec}(\fz_{\fg,x}) =
\emptyset$.
\end{lem}

This lemma implies the proposition. Indeed, for 
$\CM\in \hat\fg_{crit}\text{--}\on{mod}^{G(\wh{\CO}_x)}$
we define $\CM^{\nil}$ to be direct summand of $\CM$ supported
on the formal neighbourhood of $\Y^0$, and $\CM^{\on{non-reg}}$
to be the direct summand supported on the union of the formal 
neighborhoods of $\Y^\la$ with $\la \neq 0$.

\end{proof}

\begin{proof}(of \lemref{Sugawara calculation})

Recall the operator $S_0$ given by formula \eqref{S0}. At the critical
level this operator commutes with the action of $\hat\fg_{crit}$,
i.e., it belongs to $\FZ_{\fg,x}$. But according to formula
\eqref{S0}, $S_0$ acts on $V^\lambda \subset \BV^\lambda_{\fg,crit}$,
and hence on the entire $\BV^\lambda_{\fg,crit}$, by the scalar
$C_\fg(\la)$ equal to the value of the Casimir operator on
$V^\lambda$. This scalar is zero for $\lambda=0$ and non-zero for
$\lambda\neq 0$. This proves the lemma.

\end{proof}

\ssec{}

Recall from \secref{comm alg} that if $\CE$ is a group ind-subscheme
of the normal bundle $N(\fz_{\fg,x})$, we can introduce the subcategory
$\FZ_{\fg,x}\text{--}\on{mod}_\CE$, such that
$$\fz_{\fg,x}\text{--}\on{mod}\subset \FZ_{\fg,x}\text{--}\on{mod}_\CE\subset 
\FZ_{\fg,x}\text{--}\on{mod}.$$

We define $\CE$ as follows. By \thmref{FF isom}(1), the $D_X$-algebra $\fz_\fg$ 
is non-canonically isomorphic to a free algebra. 
This implies, in particular, that the fiber $\Theta(\fz_\fg)_x$
of $\Theta(\fz_\fg)$ at $x$ is locally free of countable rank over $\fz_{\fg,x}$,
and $N(\fz_{\fg,x})$ can be identified with the total space of the resulting
vector bundle.

Therefore, to specify a group ind-subscheme $\CE\subset N(\fz_{\fg,x})$
it would be sufficient to specify a $\fz_{\fg,x}$-submodule in 
$\Theta(\fz_\fg)_x$, which is locally a direct summand. Such a
submodule is given by the image of the anchor map
$\varpi: \Omega^1(\fz_\fg)_x\to \Theta(\fz_\fg)_x$; it is locally
a direct summand, as follows from \thmref{BD descr of Gelfand-Dikii}(1).

For $\CE$ defined in this way, let us
denote by $\hat\fg_{crit}\text{--}\on{mod}_{\CE}$ the subcategory of
$\hat\fg_{crit}\text{--}\on{mod}$ whose objects are the
$\hat\fg_{crit}$-modules, such that the underlying chiral
$\fz_\fg$-module belongs to $\FZ_{\fg,x}\text{--}\on{mod}_\CE$.  In
\secref{other results} we will prove the following

\begin{thm} \label{support on tangent bundle}
The category $\hat\fg_{crit}\text{--}\on{mod}^{G(\wh{\CO}_x)}_{\nil}$ is
contained in $\hat\fg_{crit}\text{--}\on{mod}_{\CE}$.
\end{thm}

In other words, this theorem says that the inclusion
$$\hat\fg_{crit}\text{--}\on{mod}^{G(\wh{\CO}_x)}_{\CE}:=
\hat\fg_{crit}\text{--}\on{mod}^{G(\wh{\CO}_x)}\cap\,
\hat\fg_{crit}\text{--}\on{mod}_{\CE}\subset
\hat\fg_{crit}\text{--}\on{mod}^{G(\wh{\CO}_x)}_{\nil}$$ is in fact an
equivalence.

\ssec{}   \label{Beilinson}
The following remark was suggested by A.~Beilinson:

Let $\on{Spec}(\FZ_{\fg,x,\on{m.f.}})$ be 
the smallest formal subscheme inside $\on{Spec}(\FZ_{\fg,x})$,
which contains $\on{Spec}(\fz_{\fg,x})$, and which is preserved
by the Poisson bracket. (The subscript "$\on{m.f.}$"
stands for "monodromy free".) Let $\FZ_{\fg,x}\text{--}\on{mod}_{\on{m.f.}}$
be the subcategory of $\FZ_{\fg,x}\text{--}\on{mod}$ consisiting of
modules supported on $\on{Spec}(\FZ_{\fg,x,\on{m.f.}})$, and let
$\hat\fg_{crit}\text{--}\on{mod}_{\on{m.f.}}$ be the corresponding subcategory
in $\hat\fg_{crit}\text{--}\on{mod}$. Beilinson has suggested that the following
strengthening of \thmref{support on tangent bundle} might be true:

\begin{conj}  \label{monodromy free}
The subcategory $\hat\fg_{crit}\text{--}\on{mod}^{G(\wh{\CO}_x)}_{\nil}$ is
contained in $\hat\fg_{crit}\text{--}\on{mod}_{\on{m.f.}}$. 
\end{conj}

In other words, this conjecture says that any $G(\CO_x)$-integrable $\hat\fg_{crit}$-module,
which is set-theoretically supported on $\on{Spec}(\fz_{\fg,x})$, is supported on 
the formal subscheme $\on{Spec}(\FZ_{\fg,x,\on{m.f.}})$.

If we could prove this conjecture, the
proof of \thmref{critical reformulation} would have been more elegant, 
since instead of the obscure condition (2) in the definition of 
$\FZ_{\fg,x}\text{--}\on{mod}_\CE$, we would work with a clearer
geometric concept of support on a subscheme.

\ssec{}    \label{some cat}

Recall now the category $\CA_{\fg}^{\ren,\tau}\text{--}\on{mod}$, introduced
in \secref{categories of modules} and \ref{twist by tau}. We have
a natural forgetful functor
$\CA_{\fg}^{\ren,\tau}\text{--}\on{mod}\to \hat\fg_{crit}\text{--}\on{mod}$
coming from the ``right'' copy of $\CA_{\fg,crit}$ in $\CA_{\fg}^{\ren,\tau}$.
Let $\CA_{\fg}^{\ren,\tau}\text{--}\on{mod}^{G(\wh{\CO}_x)}$
(resp., $\CA_{\fg}^{\ren,\tau}\text{--}\on{mod}_{\nil}$,
$\CA_{\fg}^{\ren,\tau}\text{--}\on{mod}_\CE$,
$\CA_{\fg}^{\ren,\tau}\text{--}\on{mod}^{G(\wh{\CO}_x)}_{\nil}$, etc.)
be the preimages of the corresponding subcategories of
$\hat\fg_{crit}\text{--}\on{mod}$ under the above forgetful functor.
Note, that by \thmref{support on tangent bundle}, the inclusion
$$\CA_{\fg}^{\ren,\tau}\text{--}\on{mod}^{G(\wh{\CO}_x)}_\CE\hookrightarrow
\CA_{\fg}^{\ren,\tau}\text{--}\on{mod}^{G(\wh{\CO}_x)}_{\nil}$$
is in fact an equivalence.

It is easy to see that the functor $\wt{\imath}{}^!:\FZ_{\fg,x}\text{--}\on{mod}\to
\FZ_{\fg,x}\text{--}\on{mod}_{\nil}$ gives rise to a well-defined functor
$\wt{\imath}{}^!: \CA_{\fg}^{\ren,\tau}\text{--}\on{mod}\to
\CA_{\fg}^{\ren,\tau}\text{--}\on{mod}_{\nil}$. In particular, the corresponding
functor
$$\CA_{\fg}^{\ren,\tau}\text{--}\on{mod}^{G(\wh{\CO}_x)} 
\overset{\wt{\imath}^!}\longrightarrow
\CA_{\fg}^{\ren,\tau}\text{--}\on{mod}^{G(\wh{\CO}_x)}_{\nil}$$
is exact, by \corref{exactness of restr to formal},.

\medskip

Recall now the Lie-* algebroid $\CA_{\fg}^{\flat,\tau}$ 
and the corresponding category $\CA_{\fg}^{\flat,\tau}\text{--}\on{mod}$
(see \secref{categories of modules}, \secref{twist by tau}).
We claim that the functor
$\imath^!:\FZ_{\fg,x}\text{--}\on{mod}\to \fz_{\fg,x}\text{--}\on{mod}$ gives
rise to a functor from $\CA_{\fg}^{\ren,\tau}\text{--}\on{mod}$ to
$\CA_{\fg}^{\flat,\tau}\text{--}\on{mod}$. 

Indeed, given an object ${\mc M}$ of the category
$\CA_{\fg}^{\ren,\tau}\text{--}\on{mod}$, we consider it as a
$\FZ_{\fg,x}$-module and take its maximal submodule $\imath^!(\CM)$
supported on $\on{Spec} (\fz_{\fg,x})$. We consider $\imath^!(\CM)$
as a Lie-* module over $\CA_{\fg}^{\ren,\tau}$. But now the Lie-* action of
the diagonal $\fz_{\fg} \subset
(\CA_{\fg,crit}\times\CA_{\fg,crit})/\fz_\fg \subset
\CA^{\ren,\tau}_\fg$ will be zero. Therefore, the Lie-* action of
$\CA_{\fg}^{\ren,\tau}$ on $\imath^!(\CM)$ will
factor through the action of the Lie-* algebra $\CA_{\fg}^{\flat,\tau}
= \CA_{\fg}^{\ren,\tau}/\fz_\fg$. Moreover, $\imath^!(\CM)$ is clearly
preserved by the chiral
action of $\CA_{\fg,crit}\underset{\fz_\fg}\otimes \CA_{\fg,crit}$,
and these two structures make $\imath^!(\CM)$ an object of
the category $\CA_{\fg}^{\flat,\tau}\text{--}\on{mod}$. By a slight abuse
of notation, we denote the resulting functor 
$\CA_{\fg}^{\ren,\tau}\text{--}\on{mod}\to \CA_{\fg}^{\flat,\tau}\text{--}\on{mod}$
also by $\imath^!$.

\medskip

In the next section we will prove the following theorem, which can be regarded
as a version of the Kashiwara theorem in the theory of D-modules.

\begin{thm} \label{Kashiwara}
The functor $\imath^!:\CA_{\fg}^{\ren,\tau}\text{--}\on{mod}_\CE\to
\CA^{\flat,\tau}_\fg\text{--}\on{mod}$ is an equivalence of categories. 
In particular, it is exact.
\end{thm}

If we denote by $\CA_{\fg}^{\flat,\tau}\text{--}\on{mod}^{G(\wh{\CO}_x)}$
the corresponding subcategory of $\CA_{\fg}^{\ren,\tau}$, we obtain that the functor
$$\CA_{\fg}^{ren,\tau}\text{--}\on{mod}^{G(\wh{\CO}_x)}_\CE\to
\CA_{\fg}^{\flat,\tau}\text{--}\on{mod}^{G(\wh{\CO}_x)}$$
is also exact (and, in fact, an equivalence).

\ssec{}

We are now able to finish the proof of \thmref{critical reformulation},
modulo Theorems \ref{support on tangent bundle} and \ref{Kashiwara}.

Recall from \thmref{embedding of algebroids} that we have a homomorphism
of chiral algebras $U^{\ren,\tau}(L_{\fg,crit})\to \fD_{G,crit}$.
Hence, the forgetful functor
$\fD_{G,crit}\text{--}\on{mod}\to \hat\fg_{crit}\text{--}\on{mod}$ factors as
$$\fD_{G,crit}\text{--}\on{mod}\to \CA_{\fg}^{\ren,\tau}\text{--}\on{mod}\to
\hat\fg_{crit}\text{--}\on{mod}.$$

We have a commutative diagram of functors
\begin{equation}   \label{comm diag funct} 
\CD
\CA_{\fg}^{ren,\tau}\text{--}\on{mod}^{G(\wh{\CO}_x)} @>{\imath^!}>>
\CA_{\fg}^{\flat,\tau}\text{--}\on{mod}^{G(\wh{\CO}_x)} \\
@VVV  @VVV  \\
\hat\fg_{crit}\text{--}\on{mod}^{G(\wh{\CO}_x)} @>{\imath^!}>>
\hat\fg_{crit}\text{--}\on{mod}^{G(\wh{\CO}_x)}_{\on{reg}},
\endCD
\end{equation}
where the vertical arrows are the forgetful functors.

Thus, to prove \thmref{critical reformulation}, it is sufficient to show
that the composition
$$\CA_{\fg}^{ren,\tau}\text{--}\on{mod}^{G(\wh{\CO}_x)}\overset{\wt{\imath}^!}\longrightarrow
\CA_{\fg}^{ren,\tau}\text{--}\on{mod}^{G(\wh{\CO}_x)}_{\nil}
\simeq
\CA_{\fg}^{ren,\tau}\text{--}\on{mod}^{G(\wh{\CO}_x)}_{\CE}\overset{\imath^!}\longrightarrow
\CA_{\fg}^{\flat,\tau}\text{--}\on{mod}^{G(\wh{\CO}_x)}$$
is exact. But, as we have just seen, all the above arrows are exact functors.

\medskip

Our plan now is as follows. In the next section we will prove
\thmref{Kashiwara} and hence complete the proof of Theorems
\ref{critical reformulation} and \ref{main}, modulo Theorems
\ref{drinfeld} and \ref{support on tangent bundle}. These theorems
will be proved simultaneously in \secref{other results}. Finally, 
in \secref{faithfulness} we will prove that the functor of global
sections, considered as a functor
$\on{D}_{crit}(\Gr_G)\text{--}\on{mod}\to
\CA_{\fg}^\flat\text{--}\on{mod}$, is fully faithful.

\section{Proof of \thmref{Kashiwara}}   \label{proof of Kashiwara}

\ssec{}  \label{Kash gen}

Let us recall the setting of the original Kashiwara theorem.  Let $\X$
be a smooth variety, and $\imath:\Y\hookrightarrow \X$ an embedding of
a smooth closed subvariety. Consider the category $D_\X\text{--}\on{mod}$
of right D-modules on $\X$, and its subcategory $D_\X\text{--}\on{mod}_\Y$
of right D-modules set-theoretically supported on $\Y$. Finally,
consider the category $D_\Y\text{--}\on{mod}$ of D-modules on $\Y$.

We have the functor $\imath^!:D_\X\text{--}\on{mod}_\Y\to
\CO_\Y\text{--}\on{mod}$ which sends a D-module $\CM$ to its maximal
$\CO_\X$-submodule consisting of sections annihilated by the ideal of
$\Y$. Then the resulting $\CO_\X$-module naturally acquires a right
action of the ring of differential operators on $\Y$, and so we obtain
a functor $\imath^!:D_\X\text{--}\on{mod}_\Y\to
D_\Y\text{--}\on{mod}$. Kashiwara's theorem asserts that this functor is an
equivalence of categories. 

Our \thmref{Kashiwara} should be regarded as a generalization of the
above theorem, when the ring of differential operators is replaced by
a certain algebroid (cf. \secref{sec Beil} below). In the proof we 
will use the same argument as in the proof of the original Kashiwara theorem. 

\ssec{}   \label{functor !}

We start by constructing a functor $\imath_!:
\CA^{\flat,\tau}_\fg\text{--}\on{mod}\to
\CA_{\fg}^{\ren,\tau}\text{--}\on{mod}_\CE$, which will be the left
adjoint of $\imath^!$.

Given an object $\CN$ in $\CA^{\flat,\tau}_\fg\text{--}\on{mod}$, we regard
it as a Lie-* module over $\CA_{\fg}^{\ren,\tau}$, considered as a
Lie-* algebra. Let $\on{Ind}(\CN)$ denote the induced 
chiral $\CA_{\fg}^{\ren,\tau}$-module (see \cite{BD}, Sect. 3.7.15),
where $\CA_{\fg}^{\ren,\tau}$ is again considered merely as a Lie-* 
algebra (and not as a chiral algebroid).

Thus, $\on{Ind}(\CN)$ is a chiral module over the chiral universal
enveloping algebra $U(\CA_{\fg}^{\ren,\tau})$. We have the surjections
$$U(\CA_{\fg}^{\ren,\tau})\twoheadrightarrow
U(\fz_\fg,\CA_{\fg}^{\ren,\tau})\twoheadrightarrow U^{\ren,\tau}(L_{\fg,crit}),$$
and we set $\imath_!(\CN)$ to be the (maximal) quotient of $\on{Ind}(\CN)$, on
which the action of $U(\CA_{\fg}^{\ren,\tau})$ factors through an action of
$U^{\ren,\tau}(L_{\fg,crit})$, and for which the two maps
$$j_*j^*\left((\CA_{\fg,crit}\underset{\fz_\fg}\otimes \CA_{\fg,crit})\boxtimes \CN\right)
\to \Delta_!\left(\imath_!(\CN)\right),$$
one coming from the inital chiral action of 
$\CA_{\fg,crit}\underset{\fz_\fg}\otimes \CA_{\fg,crit}$ on $\CN$, and the other
from the homomorphism 
$\CA_{\fg,crit}\underset{\fz_\fg}\otimes \CA_{\fg,crit}\to U^{\ren,\tau}(L_{\fg,crit})$,
coincide.

By definition, $\imath_!(\CN)$ is an object of
$\CA_{\fg}^{\ren,\tau}\text{--}\on{mod}$. It is easy to see that the
functor $\CN\mapsto
\imath_!(\CN):\CA_{\fg}^{\flat,\tau}\text{--}\on{mod}\to
\CA_{\fg}^{\ren,\tau}\text{--}\on{mod}$ is the left adjoint to
$\imath^!:\CA_{\fg}^{\ren,\tau}\text{--}\on{mod}\to
\CA_{\fg}^{\flat,\tau}\text{--}\on{mod}$. One readily checks that
for $\CN=\CA_{\fg,crit}\underset{\fz_\fg}\otimes \CA_{\fg,crit}$,
with the action of $\CA^{\flat,\tau}_\fg$ introduced in the proof of
\thmref{doubled algebroid}, the resulting
object $\imath_!(\CN)$ is isomorphic to $U^{\ren,\tau}(L_{\fg,crit})$ itself.

\ssec{}

Let us regard $\imath_!(\CN)$ as a chiral module over $\CA_{\fg,crit}\underset{\fz_\fg}\otimes
\CA_{\fg,crit}$. We have a canonical map $\CN\to \imath_!(\CN)$, and 
the PBW filtration on $U^{\ren,\tau}(L_{\fg,crit})$ induces an increasing filtration
$F^i(\imath_!(\CN)), i\geq 1$ on
$\imath_!(\CN)$ with the $F^1(\imath_!(\CN))$ term being the image of $\CN$. The terms of this
filtration are stable under the chiral action of
$\CA_{\fg,crit}\underset{\fz_\fg}\otimes \CA_{\fg,crit}$ and the
Lie-* action of the entire $\CA_{\fg}^{\ren,\tau}$;
the action of $\fz_{\fg,x}$ on $\on{gr}(\imath_!(\CN))$ is central.

The description of $\on{gr}(U^{\ren,\tau}(L_{\fg,crit}))$ given by \eqref{ass graded}
implies that we have an isomorphism
\begin{equation} \label{ass graded module}
\CN\underset{\fz_{\fg,x}}\otimes
\on{Sym}_{\fz_{\fg,x}}(\Omega^1(\fz_{\fg,x}))\simeq\on{gr}(\imath_!(\CN)).
\end{equation}

Evidently, as a module over $\FZ_{\fg,x}$, $\imath_!(\CN)$ is
supported on the formal neighborhood of $\on{Spec}(\fz_{\fg,x})$.
Recall the setting if \secref{comm alg}; in particular, let
$\CI$ be the ideal $\on{ker}(\FZ_{\fg,x}\to \fz_{\fg,x})$.

\begin{lem}
For $\CN$ as above, $F^i(\imath_!(\CN))$ equals the submodule of
$\imath_!(\CN)$ annihilated by $\CI^i$. Moreover, as a module
over $\FZ_{\fg,x}$, $\imath_!(\CN)$ belongs to the
category $\FZ_{\fg,x}\text{--}\on{mod}_\CE$.
\end{lem}

\begin{proof}

Since the action of $\fz_\fg$ on $\on{gr}(\imath_!(\CN))$ is central,
we have $\CI\cdot F^{i+1}(\imath_!(\CN))\subset F^i(\imath_!(\CN))$.
Therefore, by induction, every element of $F^i(\imath_!(\CN))$
is annihilated by $\CI^i$.

To prove the inclusion in the other direction, consider the map
$$\CI/\CI^2\underset{\fz_{\fg,x}}\otimes 
\left(F^{i+1}(\imath_!(\CN))/F^i(\imath_!(\CN))\right)\to
\left(F^{i}(\imath_!(\CN))/F^{i-1}(\imath_!(\CN))\right),$$
and the corresponding dual map
\begin{equation} \label{dual map}
\left(F^{i+1}(\imath_!(\CN))/F^i(\imath_!(\CN))\right)\to
\Theta(\fz_\fg)_x\underset{\fz_{\fg,x}}\otimes 
\left(F^{i}(\imath_!(\CN))/F^{i-1}(\imath_!(\CN))\right).
\end{equation}

We need to show that the map of \eqref{dual map} is injective.
But this follows by combining the isomorphism \eqref{ass graded module}
and \thmref{BD descr of Gelfand-Dikii}. Indeed, this map is obtained
by tensoring with $\CN$ from
$$\on{Sym}^i_{\fz_{\fg,x}}(\Omega^1(\fz_{\fg,x}))\to
\on{Sym}^{i-1}_{\fz_{\fg,x}}(\Omega^1(\fz_{\fg,x}))
\underset{\fz_{\fg,x}}\otimes \Omega^1(\fz_{\fg,x})\to
\on{Sym}^{i-1}_{\fz_{\fg,x}}(\Omega^1(\fz_{\fg,x}))
\underset{\fz_{\fg,x}}\otimes \Theta(\fz_\fg)_x.$$

This also proves the second assertion of the lemma, since $\CE$ is by definition
the image of $\Omega^1(\fz_{\fg,x})$ in $\Theta(\fz_\fg)_x$.

\end{proof}

Thus, we obtain that $\CN\mapsto \imath_!(\CN)$ is a functor
$\CA^{\flat,\tau}_\fg\text{--}\on{mod}\to \CA_{\fg}^{\ren,\tau}\text{--}\on{mod}_\CE$,
left adjoint to $\imath^!$, and the adjunction map
$\CN\mapsto \imath^!\circ \imath_!(\CN)$ is an isomorphism.

\ssec{}

To prove \thmref{Kashiwara}, it remains to show that for every 
$\CM\in \CA_{\fg}^{\ren,\tau}\text{--}\on{mod}_\CE$, the adjunction map
$\imath_!\circ \imath^!(\CM)\to \CM$ is surjective. Indeed, from the
fact that $\imath^!\circ \imath_!(\CN) \simeq \CN$, we know that the
map $\imath^! \circ \imath_!\circ \imath^!(\CM) \to \imath^!(\CM)$
is an isomorphism, and we conclude that 
$\imath^!\left(\on{Ker} (\imath_!\circ \imath^!(\CM)\to \CM)\right)=0$. 
However, the functor $\imath^!$ is evidently faithful, by condition (1) in 
the definition of $\fZ_{\fg,x}\text{--}\on{mod}_\CE$.

\medskip

Locally on $\on{Spec}(\fz_{\fg,x})$, let us choose a basis $\xi_k$ in
the space of sections of the vector bundle $\CE\subset
N(\fz_{\fg,x})$. (This is possible since $\CE\simeq \Omega(\fz_\fg)_x$,
and we know that $\Omega(\fz_\fg)$ is locally free over $\fz_\fg\otimes D_X$,
by \thmref{FF isom}(1).) Let us choose functions $f_k$ in the ideal $\CI$, 
such that under the pairing $\langle\cdot,\cdot\rangle:\CI/\CI^2\otimes \Theta(\fz_\fg)_x\to
\fz_{\fg,x}$, we have $\langle f_l,\xi_k\rangle=\delta_{k,l}$.

Let $\wh{\xi}_k$ be an arbitrary element in 
$H^0_{DR}(\D^\times_x,\CA_{\fg}^{\ren,\tau})$, which projects to
$\xi_k$ under
$$H^0_{DR}(\D^\times_x,\CA_{\fg}^{\ren,\tau})\to
H^0_{DR}(\D^\times_x,\Omega^1(\fz_\fg))\to (\Omega^1(\fz_\fg))_x\simeq
\CE.$$
Consider the natural action of $H^0_{DR}(\D^\times_x,\CA_{\fg}^{\ren,\tau})$ on
$\FZ_{\fg,x}$; we have:
$$\wh{\xi}_k(f_l)=\delta_{k,l} \; \on{mod} \; \CI.$$
 
\medskip

Let $\CM$ be an object of $\CA_{\fg}^{\ren,\tau}\text{--}\on{mod}_\CE$, and
let $\CM_i$ be the canonical filtration on it, as in \secref{comm alg}. 
We have to show that the subspace
$\CM_1$ generates $\CM$ under the action of $\CA_{\fg}^{\ren,\tau}$. We 
will argue by induction, and assume that the subspace $\CM_{i-1}$ can be
obtained from $\CM_1$ by the action of $\CA_{\fg}^{\ren,\tau}$.

Consider the action of $\CI/\CI^2$ on the extension
$$0\to \CM_{i-1}/\CM_{i-2}\to \CM_i/\CM_{i-2}\to \CM_i/\CM_{i-1}\to
0.$$ 
By definition, this action factors through $(\CI/\CI^2)/\CE^\perp$,
and the span of $f_k$'s is dense in the latter quotient.

If $m$ is an element of $\CM_i$, then $f_k\cdot m \in
\CM_{i-1}$, and for all but finitely many indices $k$ the element
$f_k\cdot m$ will belong to $\CM_{i-2}$. Therefore, the operator
$\delta:=\sum_k \wh{\xi}_k\cdot f_k$ is well-defined on
$\CM_i/\CM_{i-1}$. To perform the induction step, it would be enough
to show that $\delta$ is surjective.

We claim that $\delta$ is in fact a scalar operator that acts as
multiplication by $i-1$. To prove this statement, we assume by
induction, that $\delta$ acts as $j-1$ on $\CM_j/\CM_{j-1}$ for all
$j<i$.

Given an element $m\in \CM_i/\CM_{i-1}$, consider the finite sum
$\displaystyle \sum_{k=1,...,N} \wh{\xi}_k\cdot f_k\cdot m$, which
includes all the indices $k$ for which $f_k\cdot m\notin \CM_{i-2}$.
We must show that for any index $l$,
$$f_l\cdot \left(\sum_{k=1,...,N} \wh{\xi}_k\cdot f_k\cdot m-
(i-1)\cdot m\right)\in \CM_{i-2}.$$

Without loss of generality, we can assume that the initial finite
set of $k$'s included $l$, as well as the corresponding set of indices
for the element $f_l\cdot m\in \CM_{i-1}/\CM_{i-2}$. Then we have
\begin{align*}
&f_l\cdot \left(\sum_{k=1,...,N} \wh{\xi}_k\cdot f_k\cdot m-
(i-1)\cdot m\right)= \left(\sum_{k=1,...,N} \wh{\xi}_k\cdot f_k\cdot
f_l \cdot m -(i-2)\cdot f_l\cdot m\right)+ \\ &+\left(\sum_{l\neq
k=1,...,N} \wh{\xi}_k(f_l)\cdot f_k\cdot m\right)+
\left(\wh{\xi}_l(f_l)\cdot f_l\cdot m-f_l\cdot m\right).
\end{align*}

In the above expression, the first term belongs to $\CM_{i-2}$, by the
induction hypothesis on the action of $\delta$ on
$\CM_{i-1}/\CM_{i-2}$. The second term belongs to $\CM_{i-2}$, because
for $k\neq l$ we have $\wh{\xi}_k(f_l)\in \CI$, and the third term
also belongs to $\CM_{i-2}$, because $\wh{\xi}_l(f_l)=1\,mod\,\, \CI$.
This completes the proof of the induction step, and hence, of
\thmref{Kashiwara}.

\ssec{}   \label{sec Beil}

Recall that formal scheme $\on{Spec}(\FZ_{\fg,x,\on{m.f.}})$ introduced
in \secref{Beilinson}. From \thmref{Kashiwara} we obtain the following
corollary:

\begin{cor}
Every object of $\CA^{ren,\tau}_\fg\text{--}\on{mod}_\CE$ is supported
on $\on{Spec}(\FZ_{\fg,x,\on{m.f.}})$.
\end{cor}

\begin{proof}

Let us write $\CM\in \CA^{ren,\tau}_\fg\text{--}\on{mod}_\CE$ as
$\imath_!(\CN)$ for $\CN\in \CA^{\flat,\tau}_\fg\text{--}\on{mod}$,
and consider the filtration $F^i(\imath_!(\CN))$. 

Of course, $F^1(\imath_!(\CN))$ is supported on 
$\on{Spec}(\fz_{\fg,x})\subset \on{Spec}(\FZ_{\fg,x,\on{m.f.}})$, and
let us  assume by induction that $F^{i-1}(\imath_!(\CN))$ is supported on
$\on{Spec}(\FZ_{\fg,x,\on{m.f.}})$. 
However, since $F^{i-1}(\imath_!(\CN))$ is stable under the chiral
action of $\CA_{\fg,crit}\underset{\fz_\fg}\otimes \CA_{\fg,crit}$,
and 
$j_*j^*\left(\CA^{\ren,\tau}\boxtimes F^{i-1}(\imath_!(\CN))\right)$
maps surjectively onto $F^i(\imath_!(\CN))$, and taking into account
that $\CA^{\ren,\tau}/\left(\CA_{\fg,crit}\times \CA_{\fg,crit}/\fz_\fg\right)
\simeq \Omega^1(\fz_\fg)$, we obtain that $F^i(\imath_!(\CN))$ is also
supported on $\on{Spec}(\FZ_{\fg,x,\on{m.f.}})$.

\end{proof}

Following a suggestion of Beilinson, let us note that \thmref{Kashiwara} can be 
viewed in the following general framework. (We formulate it in the finite-dimensional 
situation, for simplicity.)

Let $\X$ and $\Y$ be as in \secref{Kash gen}. Let $L_X$ be a Lie algebroid
on $\X$, and let $L_\Y$ be its pull-back back to $\Y$. Let $\hat\Y$ be the formal
neighbourhood of $\Y$ in $\X$, and let $\hat\Y'\supset \Y$ be the smallest 
ind-subscheme of $\hat\Y$, stable under the action of $L_X$.
Let $L_X\text{--}\on{mod}$ be the category of all right $L_X$-modules,
$L_X\text{--}\on{mod}_\Y$ its subcategory of modules supported on $\hat\Y$,
and let $L_Y\text{--}\on{mod}$ be the category of right $L_Y$-modules.
Let also $L_X\text{--}\on{mod}'_\Y$ be the subcategory of 
$L_X\text{--}\on{mod}_\Y$, consisting of modules supported on $\hat\Y'$. 

We have the direct image functor $\imath_!:L_Y\text{--}\on{mod}\to
L_X\text{--}\on{mod}_\Y$, but it is easy to see that its image belongs
in fact to $L_X\text{--}\on{mod}'_\Y$. And we have the right adjoint of
$\imath_!$, denoted $\imath^!:L_X\text{--}\on{mod}_\Y\to L_Y\text{--}\on{mod}$.

Assume now that $L_X$ is $\CO_\X$-flat, and that $L_Y$ is a locally direct
summand in $N(\Y)$. Let $\CO_\X\text{--}\on{mod}_{L_Y}$ be the subcategory
of $\CO_\X\text{--}\on{mod}$ defined as in \secref{comm alg}. We have:

\begin{thm}
\hfill

\smallskip

\noindent{\em (1)}
An object $\CM\in L_X\text{--}\on{mod}_\Y$ belongs to
$L_X\text{--}\on{mod}'_\Y$ if and only if, as a $\CO_\X$-module,
it belongs to $\CO_\X\text{--}\on{mod}_{L_Y}$.

\smallskip

\noindent{\em (2)}
The functor $\imath_!$ is an equivalence of categories
$L_Y\text{--}\on{mod}\to L_X\text{--}\on{mod}'_\Y$ with
the quasi-inverse given by $\imath^!$.
\end{thm}

\section{Proof of Theorems \ref{support on tangent bundle}
and \ref{drinfeld}.}   \label{other results}

\ssec{}

Let us start by observing that we have a natural map
\begin{equation}   \label{deform}
\Omega^1(\fz_{\fg,x}) \to
\on{Ext}^1_{\hat\fg_{crit}\text{--}\on{mod}}(\BV_{\fg,crit},\BV_{\fg,crit}). 
\end{equation}

This map can be constructed in the following general framework. Let
$\CA_{\BC[\hslash]/\hslash^2}$ be a flat $\BC[\hslash]/\hslash^2$-family
of chiral algebras, and set $\CA=\left(\CA_{\BC[\hslash]/\hslash^2}\right)/\hslash$.
We claim that there is a canonical map
\begin{equation}  \label{gen deform}
\Omega^1(\fz(\CA)_x)\to \on{Ext}^1_{\CA\text{--}\on{mod}}(\CA_x,\CA_x).
\end{equation}

Indeed, consider $(\CA_{\BC[\hslash]/\hslash^2})_x$ as an extension
$$0\to \CA_x\to (\CA_{\BC[\hslash]/\hslash^2})_x\to \CA_x\to 0$$
in the category of $\CA_{\BC[\hslash]/\hslash^2}$-modules.
Let ${\mathbf e}$ denote its class in 
$\on{Ext}^1_{\CA_{\BC[\hslash]/\hslash^2}\text{--}\on{mod}}(\CA_x,\CA_x)$.

For an element $a\in \fz(\CA)_x$, viewed as an endomorphism of 
$\CA_x$ as a $\CA_{\BC[\hslash]/\hslash^2}$-module, we can produce two more elements
of $\on{Ext}^1_{\CA_{\BC[\hslash]/\hslash^2}\text{--}\on{mod}}(\CA_x,\CA_x)$,
namely, $a\cdot {\mathbf e}$ and ${\mathbf e}\cdot a$. However, it is
easy to see that their difference already belongs to
$\on{Ext}^1_{\CA\text{--}\on{mod}}(\CA_x,\CA_x)$.
Moreover, one readily checks that the resulting map
$\fz(\CA)_x\to \on{Ext}^1_{\CA\text{--}\on{mod}}(\CA_x,\CA_x)$ is a derivation,
i.e., gives rise to a map in \eqref{gen deform}.

\medskip

Explicitly, for $\CA=\CA_\fg$
the map of \eqref{deform} looks as follows. Note that
$$\on{Ext}^1_{\hat\fg_{crit}\text{--}\on{mod}}(\BV_{\fg,crit},\BV_{\fg,crit})\simeq
\on{H}^1(\fg(\wh{\CO}_x),\BV_{\fg,crit}).$$ Given an element $da\in
\Omega^1(\fz_{\fg,x})$, where $a \in \fz_{\fg,x} \subset
\BV_{\fg,crit}$, consider its deformation $a_\hslash\in
\BV_{\fg,\hslash}$ and define the corresponding 1-cocycle on
$\fg(\wh{\CO}_x)$ by
\begin{equation}    \label{cocycle}
g\in \fg(\wh{\CO}_x)\mapsto \left. \frac{g\cdot
a_\hslash}{\hslash}\right|_{\hslash=0}.
\end{equation}
This gives the desired map.

\medskip

The next proposition states that the map of
\eqref{deform} is an isomorphism. This is a particular
case of the following general theorem established in \cite{FT}:

\begin{thm} \label{FT}
We have a canonical isomorphism
between $\Omega^i(\fz_{\fg,x})$ and the relative cohomology
$\on{H}^i(\fg(\wh{\CO}_x),\fg,\BV_{\fg,crit})$.
\end{thm}

For $i=1$ we have
$\on{H}^1(\fg(\wh{\CO}_x),\fg,\BV_{\fg,crit})\simeq
\on{H}^1(\fg(\wh{\CO}_x),\BV_{\fg,crit})$, and we obtain 
that \eqref{deform} is indeed an isomorphism.

Here we will give a different proof of this fact, using 
some results from \cite{BD1}.

\begin{prop}  \label{Ext^1}
The above map $\Omega^1(\fz_{\fg,x}) \to
\on{Ext}^1_{\fg_{crit}}(\BV_{\fg,crit},\BV_{\fg,crit})$ is an isomorphism.
\end{prop}

\begin{proof}

Recall the setting of \secref{diff op}. Consider the short exact
sequence
$$0\to \BV_{\fg,crit}\overset{\fr}\longrightarrow \fD_{G,crit,x}\to
\fD_{G,crit,x}/\fr(\BV_{\fg,crit})\to 0.$$ We know that
$\on{H}^0(\fg(\wh{\CO}_x),\fD_{G,crit,x})\simeq
\fl(\BV_{\fg,crit})$, and by a similar argument we obtain that
$\on{H}^1(\fg(\wh{\CO}_x),\fD_{G,crit,x})=0$, since
$\on{H}^1\left(\fg(\wh{\CO}_x),\on{Fun}\left(G(\wh{\CO}_x)\right)\right)=0$.  
Therefore from the
long exact sequence we obtain an isomorphism
$$\on{H}^1(\fg(\wh{\CO}_x),\BV_{\fg,crit})\simeq
\bigl(\fD_{G,crit,x}/
\fr(\BV_{\fg,crit})\bigr)^{\fr(\fg(\wh{\CO}_x))}/\fl(\BV_{\fg,crit}).$$

By letting the point $x$ move, we obtain from the subspace 
$$\bigl(\fD_{G,crit,x}/
\fr(\BV_{\fg,crit})\bigr)^{\fr(\fg(\wh{\CO}_x))}/\fl(\BV_{\fg,crit})
\subset \fD_{G,crit,x}/(\fr(\BV_{\fg,crit})+\fl(\BV_{\fg,crit}))$$ a
D-submodule, which we will denote by $\wt{\Omega}^1(\fz_\fg) \subset
\fD_{G,crit}/(\fr(\CA_{\fg,crit})+\fl(\CA_{\fg,crit}))$. We can form the
Cartesian squares
\begin{align*}
&\wt{\CA}^\flat_\fg:=
\wt{\Omega}^1(\fz_\fg)\underset{\fD_{G,crit}/(\fr(\CA_{\fg,crit}) +
\fl(\CA_{\fg,crit}))} \times \fD_{G,crit}/\fr(\CA_{\fg,crit}), \\
&\wt{\CA}^{\ren,\tau}_\fg :=
\wt{\Omega}^1(\fz_\fg)\underset{\fD_{G,crit}/(\fr(\CA_{\fg,crit}) +
\fl(\CA_{\fg,crit}))} \times \fD_{G,crit}.
\end{align*}

Let us show that the chiral bracket on $\fD_{G,crit}$ induces
on $\wt{\Omega}^1(\fz_\fg)$ and $\wt{\CA}^\flat_\fg$ a structure of
Lie-* algebroids over $\fz_\fg$, and on $\wt{\CA}^{\ren,\tau}_\fg$ a
structure of chiral algebroid. 

For that, let us note that $\wt{\CA}^{\ren,\tau}_\fg$ is by definition
the normalizer on $\fr(\CA_{\fg,crit})$ in $\fD_{G,crit}$, i.e., the
maximal D-submodule of $\fD_{G,crit}$, for which the Lie-* bracket
sends $\fr(\CA_{\fg,crit})$ to itself. Since $\fl(\CA_{\fg,crit})$
is the centralizer of $\fr(\CA_{\fg,crit})$ (by \lemref{two centralizers}),
we obtain that $\wt{\CA}^{\ren,\tau}_\fg$ normalizes $\fl(\CA_{\fg,crit})$
as well. (By symmetry, we immediately obtain that $\wt{\CA}^{\ren,\tau}_\fg$
is in fact the entire normalizer of $\fl(\CA_{\fg,crit})$ in
$\fD_{G,crit}$.) Now, \corref{centers} implies that $\wt{\CA}^{\ren,\tau}_\fg$
normalizes also $\fz_\fg$. This implies the above assertion about the
algebroid structures.

Note also, that if $\CM$ is a chiral module over $\fD_{G,crit}$, we obtain
that the Lie-* algebroid $\wt{\CA}^\flat_\fg$ acts naturally on the subspace
$\CM^{\fr(\fg(\wh{\CO}_x))}$, and $\wt{\Omega}^1(\fz_\fg)$ acts on the
subspace $\CM^{\fl(\fg(\wh{\CO}_x))}\cap \CM^{\fr(\fg(\wh{\CO}_x))}$.

\medskip

From the construction of the map 
$\Omega^1(\fz_\fg)\to
\fD_{G,crit}/\fr(\CA_{\fg,crit})+\fl(\CA_{\fg,crit})$ 
in \secref{diff op}, it is clear that its image is in $\wt{\Omega}^1(\fz_\fg)$. 
(Of course, we are about to prove that $\Omega^1(\fz_\fg)$ is in fact
isomorphic to $\wt{\Omega}^1(\fz_\fg)$.) On the
level of fibers, the above map coincides with the map
$$\Omega^1(\fz_\fg)\underset{\fz_\fg}\otimes \fz_{\fg,x} =
\Omega^1(\fz_{\fg,x}) \to \on{H}^1(\fg(\wh{\CO}_x),\BV_{\fg,crit}),$$
given by formula \eqref{cocycle}.  We obtain also morphisms of
algebroids $\CA^\flat_\fg\to \wt{\CA}^\flat_\fg$ and
$\wt{\CA}^{\ren,\tau}_\fg\to \CA^{\ren,\tau}_\fg$.

\medskip

Thus, we have a sequence of morphisms of algebroids
$\Omega^1(\fz_\fg)\to \wt{\Omega}^1(\fz_\fg)\to \Theta(\fz_\fg)$, and we
want to prove that the first arrow is an isomorphism. At the level of
fibers, the second map can be viewed as follows. For an extension
$$0\to \BV_{\fg,crit}\to \CM\to \BV_{\fg,crit}\to 0,$$ the action of
$\CI/\CI^2\simeq N^*(\fz_{\fg,x})$ on $\CM$ defines an endomorphism
of $\BV_{\fg,crit}$, i.e., an element of $\fz_{\fg,x}$. 

According to
Proposition 6.2.4 of \cite{BD1}, there are no non-trivial extensions
of $\BV_{\fg,crit}$ by itself in the category
$\hat\fg_{crit}\text{--}\on{mod}_{\reg}$. This implies that
$\wt{\Omega}^1(\fz_\fg)\to \Theta(\fz_\fg)$ is an injection.
Indeed, otherwise, we would obtain an extension $\CM$ as above,
which belongs to $\hat\fg_{crit}\text{--}\on{mod}_{\reg}$.

\medskip

Thus, $\wt{\Omega}^1(\fz_\fg)$ is ``squeezed'' between
$\Omega^1(\fz_\fg)$ and $\Theta(\fz_\fg)$. To prove that
$\wt{\Omega}^1(\fz_\fg)$ in fact coincides with
$\Omega^1(\fz_{\fg,x})$, we will use Theorem 5.5.3 of \cite{BD1}.
This theorem asserts that $\Omega^1(\fz_\fg)$ coincides with the Atiyah
algebroid corresponding to a certain principal $^L G$-bundle over
$\on{Spec}(\fz_{\fg})$. 

This bundle is constructed as follows. We have an equivalence
between the category $\on{Rep}({}^L G)$ of finite-dimensional
representations of $^L G$ and the category of
$G(\wh{\CO}_x)$-equivariant objects in
$\on{D}_{crit}(\Gr_G)\text{--}\on{mod}$; for $V\in \on{Rep}({}^L G)$, let
$\F^V\in \on{D}_{crit}(\Gr_G)\text{--}\on{mod}$ be the corresponding
D-module. Set
$$\CV_{\fz_{\fg,x}}:=
\on{Hom}_{\hat\fg_{crit}}(\BV_{\fg,crit},\Gamma(\Gr_G,\F^V)).$$
Theorem 5.4.8 and Sect. 5.5.1 of \cite{BD1} imply that
$V\mapsto \CV_{\fz_{\fg,x}}$ is a tensor functor from $\on{Rep}({}^L G)$ 
to the category of locally free finitely generated $\fz_{\fg,x}$-modules.
By letting the point $x$ move along the curve, for every $V$ as above,
we obtain a $\fz_\fg\otimes D_X$-module, denoted $\CV_{\fz_\fg}$,
and the assignment $V\mapsto \CV_{\fz_\fg}$ is the sought-for
principal $^L G$-bundle on $\on{Spec}(\fz_{\fg})$.

The assertion of Theorem 5.5.3 of \cite{BD1},
combined with \thmref{BD descr of Gelfand-Dikii}, implies that
the Lie-* algebroid $\Omega^1(\fz_\fg)$ is the universal Lie-* 
algebroid over $\fz_\fg$,
whose action lifts to the above $^L G$-bundle. However, as we have
seen above, the Lie-* algebroid $\wt{\Omega}^1(\fz_\fg)$ acts on
every $\fz_{\fg,x}$-module of the form 
$\on{Hom}_{\hat\fg_{crit}}(\BV_{\fg,crit},\Gamma(\Gr_G,\F))$, for
$\F\in \on{D}_{crit}(\Gr_G)\text{--}\on{mod}$. Again, globally over $X$,
we obtain that $\wt{\Omega}^1(\fz_\fg)$ acts on all $\CV_{\fz_\fg}$.
This implies that we have a splitting $\wt{\Omega}^1(\fz_\fg)\to
\Omega^1(\fz_\fg)$, which is automatically an isomorphism.

\end{proof}

\ssec{}  \label{flatness}

Next we will prove that 
$\on{Ext}^i_{\hat\fg_{crit}\text{--}\on{mod}}(\BV_{\fg,crit},\BV_{\fg,crit})$
is flat as a $\fz_{\fg,x}$-module for any $i$. This statement can be formally
deduced from \thmref{FT}, but we will give a different proof.

First, we claim that the topological Lie algebroid $\hat{h}^{\fz_\fg}_x(\Omega^1(\fz_\fg))$
over $\on{Spec}(\fz_{\fg,x})$
(defined as in \cite{BD}, Sect. 2.5.18) acts on every such $\on{Ext}^i$.
Indeed, consider the Lie algebra $H^0_{DR}(\D_x,\Omega^1(\fz_\fg))$.
The Lie-* action of $\CA^\flat_{\fg}$ on $\CA_{\fg,crit}$ yields an action
of $H^0_{DR}(\D_x,\Omega^1(\fz_\fg))$ on the associative algebra
$\hat\CA_{\fg,crit,x}$ by outer derivations. Since the $\CA_{\fg,crit}$-action
on $\BV_{\fg,crit}$ lifts to an action of $\CA^\flat_{\fg}$, we obtain that
$H^0_{DR}(\D_x,\Omega^1(\fz_\fg))$ indeed acts on every
$\on{Ext}^i_{\hat\fg_{crit}\text{--}\on{mod}}(\BV_{\fg,crit},\BV_{\fg,crit})$.
Since $\BV_{\fg,crit}$ is a $\fz_\fg$-module, so is $\on{Ext}^i$, and the
above $H^0_{DR}(\D_x,\Omega^1(\fz_\fg))$-action extends to an action
of its completion $\hat{h}^{\fz_\fg}_x(\Omega^1(\fz_\fg))$.

\medskip

By identifying ${\mathcal K}_x$ with $\BC((t))$, we endow $\hat\fg_{crit}$
with a $\BZ$-grading by letting $t$ have degree $-1$. In this case
$\BV_{\fg,crit}$ is a non-negatively graded $\hat\fg_{crit}$-module. Moreover,
the terms of the standard complex computing the cohomology
$\on{Ext}^i_{\hat\fg_{crit}\text{--}\on{mod}}(\BV_{\fg,crit},\BV_{\fg,crit})\simeq
\on{H}^i(\fg(\wh{\CO}_x),\BV_{\fg,crit})$ are also non-negatively graded.

By applying Lemma 6.2.2 of \cite{BD1}, we conclude that the above $\on{Ext}^i$
is free over $\fz_{\fg,x}$.

\begin{lem}   \label{freeness}
The module $\BV_{\fg,crit}$ is flat over $\fz_{\fg,x}$.
\end{lem}

\begin{proof}

The lemma is proved by passing to the associate graded.  Recall
from \thmref{FF isom}(1) that $\BV_{\fg,crit}$ and $\fz_{\fg,x}$ are
naturally filtered, and
\begin{align*}
&\on{gr}(\BV_{\fg,crit})\simeq
\on{Fun}\left(\g^*\times_{\BG_m} \Gamma({\mc D}_x,\Omega_X)\right)\simeq
\CJ\left(\on{Fun}(\fg^*\times_{\BG_m}\omega_X)\right)_x \\
&\on{gr} (\fz_{\fg,x})\simeq\on{Fun}\left((\g^*/G) \times_{\BG_m}
\Gamma({\mc D}_x,\Omega_X)\right) \simeq
\CJ\left(\on{Fun}(\fg^*/G\times_{\BG_m}\omega_X)\right)_x.
\end{align*}
Now we apply Theorem A.4 of \cite{EF}, which exactly asserts that 
$\CJ\left(\on{Fun}(\fg^*\times_{\BG_m}\omega_X)\right)$ is flat
over $\CJ\left(\on{Fun}(\fg^*/G\times_{\BG_m}\omega_X)\right)$.
\end{proof}

Finally, we are ready to prove the following:

\begin{prop}  \label{all Exts}
For any $\fz_{\fg,x}$-module $\CL$, the natural map
$$\on{Ext}^i_{\hat\fg_{crit}\text{--}\on{mod}}(\BV_{\fg,crit},
\BV_{\fg,crit})
\underset{\fz_{\fg,x}}\otimes \CL\to
\on{Ext}^i_{\hat\fg_{crit}\text{--}\on{mod}}(\BV_{\fg,crit},
\BV_{\fg,crit}\underset{\fz_{\fg,x}}\otimes \CL)$$
is an isomorphism.
\end{prop}

\begin{proof}

From the identification
$\on{Ext}^i_{\hat\fg_{crit}\text{--}\on{mod}}(\BV_{\fg,crit},
\BV_{\fg,crit}\underset{\fz_{\fg,x}}\otimes \CL)\simeq
\on{H}^i(\fg(\wh{\CO}_x),\BV_{\fg,crit}\underset{\fz_{\fg,x}}\otimes \CL)$
and the standard complex, computing cohomology of the Lie algebra
$\fg(\wh{\CO}_x)$, we obtain that the functor
$\CL\mapsto \on{Ext}^i_{\hat\fg_{crit}\text{--}\on{mod}}(\BV_{\fg,crit},
\BV_{\fg,crit}\underset{\fz_{\fg,x}}\otimes \CL)$
commutes with direct limits. Therefore, to prove the proposition,
we can suppose that $\CL$ is finitely presented.

Since $\fz_{\fg,x}$ is isomorphic to a polynomial algebra (by \thmref{FF isom}(1)),
any finitely presented module $\CL$ admits a finite resolution by 
projective modules:
$$0\to \CP_n\to...\to \CP_1\to \CP_0\to \CL\to 0.$$
By \lemref{freeness}, the complex
$$0\to \BV_{\fg,crit}\underset{\fz_{\fg,x}}\otimes\CP_n\to...\to
\BV_{\fg,crit}\underset{\fz_{\fg,x}}\otimes\CP_1\to
\BV_{\fg,crit}\underset{\fz_{\fg,x}}\otimes\CP_0\to 
\BV_{\fg,crit}\underset{\fz_{\fg,x}}\otimes\CL\to 0$$ 
is also exact.

Thus, we have a spectral sequence, converging to 
$\on{Ext}^i_{\hat\fg_{crit}\text{--}\on{mod}}(\BV_{\fg,crit},
\BV_{\fg,crit}\underset{\fz_{\fg,x}}\otimes \CL)$, with the $E^{i,-j}_1$-term 
isomorphic to
$\on{Ext}^i_{\hat\fg_{crit}\text{--}\on{mod}}(\BV_{\fg,crit},
\BV_{\fg,crit}\underset{\fz_{\fg,x}}\otimes \CP^j)$.

Since each $\CP^j$ is projective, we evidently have 
$$\on{Ext}^i_{\hat\fg_{crit}\text{--}\on{mod}}(\BV_{\fg,crit},
\BV_{\fg,crit}\underset{\fz_{\fg,x}}\otimes \CP^j)\simeq
\on{Ext}^i_{\hat\fg_{crit}\text{--}\on{mod}}(\BV_{\fg,crit},
\BV_{\fg,crit}) \underset{\fz_{\fg,x}}\otimes \CP^j.$$
But since all $\on{Ext}^i_{\hat\fg_{crit}\text{--}\on{mod}}(\BV_{\fg,crit},
\BV_{\fg,crit})$ are $\fz_{\fg,x}$-flat, this spectral sequence
degenerates at $E_2$, implying the assertion of the proposition.

\end{proof}

\begin{cor}  \label{vanishing of ext}
$\on{Ext}^1_{\hat\fg_{crit}\text{--}\on{mod}_{\reg}}(\BV_{\fg,crit},
\BV_{\fg,crit}\underset{\fz_{\fg,x}}\otimes \CL)=0$.
\end{cor}

\begin{proof}

We have a map
$$\CI/\CI^2\underset{\fz_{\fg,x}}\otimes
\on{Ext}^1_{\hat\fg_{crit}\text{--}\on{mod}}(\BV_{\fg,crit},
\BV_{\fg,crit}\underset{\fz_{\fg,x}}\otimes \CL)\to
\on{Hom}_{\hat\fg_{crit}\text{--}\on{mod}}(\BV_{\fg,crit},
\BV_{\fg,crit}\underset{\fz_{\fg,x}}\otimes \CL),$$
and its adjoint 
\begin{equation}  \label{dmap}
\on{Ext}^1_{\hat\fg_{crit}\text{--}\on{mod}}(\BV_{\fg,crit},
\BV_{\fg,crit}\underset{\fz_{\fg,x}}\otimes \CL)\to
\Theta(\fz_\fg)_x\underset{\fz_{\fg,x}}\otimes
\on{Hom}_{\hat\fg_{crit}\text{--}\on{mod}}(\BV_{\fg,crit},
\BV_{\fg,crit}\underset{\fz_{\fg,x}}\otimes \CL).
\end{equation}
It is easy to see that
$$\on{Ext}^1_{\hat\fg_{crit}\text{--}\on{mod}_{\reg}}(\BV_{\fg,crit},
\BV_{\fg,crit}\underset{\fz_{\fg,x}}\otimes \CL)\subset
\on{Ext}^1_{\hat\fg_{crit}\text{--}\on{mod}}(\BV_{\fg,crit},
\BV_{\fg,crit}\underset{\fz_{\fg,x}}\otimes \CL)$$
is exactly the kernel of the latter map.

However, by \propref{all Exts} applied to $i=0$ and $1$,
we can identify both sides in \eqref{dmap} with
$$\Omega^1(\fz_{\fg,x})\underset{\fz_{\fg,x}}\otimes \CL\to
\Theta(\fz_{\fg,x})\underset{\fz_{\fg,x}}\otimes \CL,$$
and the latter map is injective, since
$\on{coker}(\Omega^1(\fz_\fg)\to \Theta(\fz_{\fg,x}))$ is
flat as a $\fz_\fg$-module, by \thmref{BD descr of Gelfand-Dikii}.

\end{proof}

\ssec{}

Recall the functor $\sF:\hat\fg_{crit} \text{--}
\on{mod}_{\reg}^{G(\wh{\CO}_x)}\to \fz_{\fg,x}\text{--}\on{mod}$ and its
left adjoint $\sG:\fz_{\fg,x}\text{--}\on{mod}\to
\hat\fg_{crit}\text{--}\on{mod}_{\reg}^{G(\wh{\CO}_x)}$ defined in
\secref{functors}.

Note that the functor $\sF$ is faithful. Indeed, a module $\CM\in
\hat\fg_{crit}\text{--}\on{mod}_{\reg}^{G(\wh{\CO}_x)}$ necessarily
contains a non-zero vector, which is annihilated by $\fg\otimes t\BC[[t]]\subset 
\hat\fg_{crit}$. Therefore we have a non-zero map
$\BV^\lambda_{\fg,crit} \to \CM$, and by \lemref{Sugawara
calculation}, $\lambda$ must be equal to $0$. 

Note that the assertion of \propref{all Exts} for $i=0$ implies that
the adjunction morphism $$\CL\to \sF\circ \sG(\CL)$$ is an isomorphism. 
We claim that this, combined with \corref{vanishing of ext}, formally implies 
\thmref{drinfeld}.

\medskip

We have to show that for $\CM\in \hat\fg_{crit}\text{--}\on{mod}_{\reg}^{G(\wh{\CO}_x)}$
the adjunction map
$$\sG\circ \sF(\CM)\to \CM$$
is an isomorphism.

This map is injective. Indeed, if $\CM'$ is the kernel of
$\sG\circ \sF(\CM)\to \CM$, by the left exactness of $\sF$, we would
obtain that 
$$\sF(\CM')=\on{ker}\bigl(\sF\circ \sG\circ \sF(\CM)\to \sF(\CM)\bigr)\simeq
\on{ker}\bigl(\sF(\CM)\to \sF(\CM)\bigr)=0.$$
But we know that the functor $\sF$ is faithful, so $\CM'=0$.

Let us prove that $\sG\circ \sF\to \on{Id}$ is surjective.  Let
$\CM''$ be the cokernel of $\sG\circ \sF(\CM)\to \CM$. We have the
long exact sequence
$$0\to \sF\circ \sG\circ \sF(\CM) \to \sF(\CM)\to \sF(\CM'')\to
R^1\sF(\sG\circ \sF(\CM))\to...$$ 
However, \corref{vanishing of ext} implies that 
$R^1\sF(\sG(\CL))=0$ for any $\fz_{\fg,x}$-module $\CL$.
Therefore, the above portion of the long exact sequence amounts to
a short exact sequence
$$0\to \sF\circ \sG\circ \sF(\CM) \to \sF(\CM)\to \sF(\CM'')\to 0.$$
But the first arrow is an isomorphism, which implies that
$\sF(\CM'')=0$ and hence $\CM''=0$.
Thus, \thmref{Drinfeld's equivalence of categories} is proved.

\medskip

As a corollary, we obtain the following result. Let $\sigma\in
\on{Spec}(\fz_{\fg,x})$ be a $\BC$-point, and consider the subcategory
$\hat\fg_{crit}\text{--}\on{mod}_{\sigma}^{G(\wh{\CO}_x)}$ of
$\hat\fg_{crit}\text{--}\on{mod}_{\reg}^{G(\wh{\CO}_x)}$ whose
objects are the $\hat{\fg}_{crit}$-modules with central character
equal to $\sigma$. \thmref{Drinfeld's equivalence of
categories} implies that the category
$\hat\fg_{crit}\text{--}\on{mod}_{\sigma}^{G(\wh{\CO}_x)}$ is
equivalent to the category of vector spaces. In particular, the module
$$\BV_{\fg,\sigma}:=\BV_{\fg,crit}\underset{\fz_{\fg,x}}\otimes
\CC_\sigma$$ is irreducible.

\ssec{}

Finally, let us prove \thmref{support on tangent bundle}.
Let $\CM$ be an object of 
$\hat\fg_{crit}\text{--}\on{mod}^{G(\wh{\CO}_x)}$, and let
$\CM_i$ be the filtration as in \secref{comm alg}.
We need to show that the action of $\CI/\CI^2$ on
$\CM_{i+1}/\CM_{i-1}$, that maps $\CM_{i+1}/\CM_i$ to
$\CM_i/\CM_{i-1}$, factors through $(\CI/\CI^2)/\CE^\perp$.

We claim that for any extension in $\hat\fg_{crit}\text{--}\on{mod}^{G(\wh{\CO}_x)}$
$$0\to \CM^1\to \CM^2\to\CM^3\to 0$$
with $\CM^1,\CM^3\in \hat\fg_{crit}\text{--}\on{mod}^{G(\wh{\CO}_x)}_{\reg}$,
the map $\CI/\CI^2\underset{\fz_{\fg,x}}\otimes \CM_3\to \CM_1$ factors
through $(\CI/\CI^2)/\CE^\perp$.

Indeed, by \thmref{drinfeld}, we can map surjectively onto $\CM^3$ a module of the 
form $\BV_{\fg,crit}\underset{\fz_{\fg,x}}\otimes \CL$, where $\CL$ is a free
$\fz_{\fg,x}$-module. Hence, we can replace $\CM^3$ by $\BV_{\fg,crit}$.
Again, by \thmref{drinfeld}, $\CM^1$ has the form
$\BV_{\fg,crit}\underset{\fz_{\fg,x}}\otimes \CL^1$ for some $\fz_{\fg,x}$-module
$\CL^1$.

Now, our assertion follows from \propref{all Exts} for $i=1$.

\section{Faithfulness}   \label{faithfulness}

\ssec{}

Recall the category $\CA_{\fg}^{\flat,\tau}\text{--}\on{mod}^{G(\wh{\CO}_x)}$
introduced in \secref{some cat}. Observe that the functor 
$\sF:\hat\fg_{crit}\text{--}\on{mod}^{G(\wh{\CO}_x)}_{\reg}\to
\fz_{\fg,x}\text{--}\on{mod}$ extends to a functor
\begin{equation} \label{dr with tau}
\CA_{\fg}^{\flat,\tau}\text{--}\on{mod}^{G(\wh{\CO}_x)}\to 
\CA_{\fg}^\flat\text{--}\on{mod}.
\end{equation}
Moreover, \thmref{Drinfeld's equivalence of categories} implies that the latter is also
an equivalence of categories, with the quasi-inverse being 
$\CM\to \CM\underset{\fz_{\fg,x}}\otimes \BV_{\fg,crit}$.

\medskip

We obtain that the functor
$\Gamma:\on{D}_{crit}(\Gr_G)\text{--}\on{mod}\to
\hat\fg_{crit}\text{--}\on{mod}$ naturally factors as
$$\on{D}_{crit}(\Gr_G)\text{--}\on{mod}\to \CA_{\fg}^\flat\text{--}\on{mod}\to
\hat\fg_{crit}\text{--}\on{mod}_{\reg}\hookrightarrow 
\hat\fg_{crit}\text{--}\on{mod},$$
where the second arrow is the tautological forgetful functor.
We will denote the resulting functor
$\on{D}_{crit}(\Gr_G)\text{--}\on{mod}\to \CA_{\fg}^\flat\text{--}\on{mod}$
by $\Gamma^\flat$.

\medskip

\noindent{\it Remark.}
Suppose that $\F_\hslash$ is a $\BC[[\hslash]]$-flat family of 
$\kappa_\hslash$-twisted D-modules on $\Gr_G$. By taking global
sections, we obtain a $\BC[[\hslash]]$-family of 
$\hat\fg_{\kappa_\hslash}$-modules. \thmref{main} implies
that this family is flat as well. 

Set $\F_0=\F/\hslash$. We obtain that 
the $\hat\fg_{crit}$-module $\Gamma(\Gr_G,\F_0)$ has two
(a priori different) structures of object of
$\CA_{\fg}^\flat\text{--}\on{mod}$: one such structure has been
described above, and another is as in \secref{flat modules}. However,
it is easy to see that these structures in fact coincide.

\medskip

The main result of this section is the following

\begin{thm}  \label{fully faithfulness}
The above functor $\Gamma^\flat:\on{D}_{crit}(\Gr_G)\text{--}\on{mod}\to 
\CA_{\fg}^\flat\text{--}\on{mod}$ is fully faithful.
\end{thm}

\ssec{}

Recall the category $\CA_{\fg}^{\ren,\tau}\text{--}\on{mod}^{G(\wh{\CO}_x)}_{\nil}$.
By combining Theorems
\ref{Kashiwara} and \ref{support on tangent bundle} we obtain:

\begin{cor}   \label{total equivalence}
We have the following sequence of equivalences:
$$\CA_{\fg}^{\ren,\tau}\text{--}\on{mod}^{G(\wh{\CO}_x)}_{\nil}
\overset{\imath^!}\longrightarrow 
\CA_{\fg}^{\flat,\tau}\text{--}\on{mod}^{G(\wh{\CO}_x)}\to
\CA_{\fg}^\flat\text{--}\on{mod},$$
where the last functor is as in \eqref{dr with tau}.
\end{cor}

Recall now the functor
$\wt{\imath}^!:\hat\fg_{crit}\text{--}\on{mod}^{G(\wh{\CO}_x)}\to
\hat\fg_{crit}\text{--}\on{mod}^{G(\wh{\CO}_x)}_{\nil}$, and the corresponding functor
$$\wt{\imath}^!:\CA_{\fg}^{\ren,\tau}\text{--}\on{mod}^{G(\wh{\CO}_x)}\to
\CA_{\fg}^{\ren,\tau}\text{--}\on{mod}^{G(\wh{\CO}_x)}_{\nil}.$$

Let us denote by $\wt{\Gamma}$ the functor
$\on{D}_{crit}(\Gr_G)\text{--}\on{mod}\to 
\CA_{\fg}^{\ren,\tau}\text{--}\on{mod}^{G(\wh{\CO}_x)}_{\nil}$ equal to the
composition
$$\on{D}_{crit}(\Gr_G)\text{--}\on{mod}\simeq \fD_{G,crit}\text{--}\on{mod}^{G(\wh{\CO}_x)}
\to \CA_{\fg}^{\ren,\tau}\text{--}\on{mod}^{G(\wh{\CO}_x)}\overset{\wt{\imath}^!}
\longrightarrow \CA_{\fg}^{\ren,\tau}\text{--}\on{mod}^{G(\wh{\CO}_x)}_{\nil}.$$

The functor $\Gamma^\flat$ is the composition of $\wt{\Gamma}$,
followed by the $\CA_{\fg}^{\ren,\tau}\text{--}\on{mod}^{G(\wh{\CO}_x)}_{\nil}\to
\CA_{\fg}^\flat\text{--}\on{mod}$ of \corref{total equivalence}. So we are reduced to proving

\begin{thm} \label{sectors}
The functor $\wt{\Gamma}:\fD_{G,crit}\text{--}\on{mod}^{G(\wh{\CO}_x)}\to 
\CA_{\fg}^{\ren,\tau}\text{--}\on{mod}^{G(\wh{\CO}_x)}_{\nil}$
is fully faithful.
\end{thm}

\ssec{}

Recall from \propref{restriction to formal} that for
$\CM\in \hat\fg_{crit}\text{--}\on{mod}^{G(\wh{\CO}_x)}$, the object
$\wt{\imath}^!(\CM)$ is in fact a direct summand of $\CM$, denoted
$\CM^{\nil}$.

Let us first analyze the decomposition $\CM\simeq
\CM^{\nil}\oplus \CM^{\on{non-reg}}$ when the module $\CM$ equals $\fD_{G,crit,x}$
itself. Letting $x$ move, we obtain a direct sum decomposition of
D-modules $\fD_{G,crit}\simeq \fD^{\nil}_{G,crit}\oplus \fD^{\on{non-reg}}_{G,crit}$.

\begin{lem}   \label{neutral component}
The homomorphism $U^{\ren,\tau}(L_{\fg,crit})\to \fD_{G,crit}$
is an isomorphism onto $\fD^{\nil}_{G,crit}$. In particular,
$\fD^{\nil}_{G,crit}$ is a chiral subalgebra of $\fD_{G,crit}$.
\end{lem}

\begin{proof}

It is enough to show that the homomorphism
$U^{\ren,\tau}(L_{\fg,crit})\to \fD_{G,crit}$ induces an isomorphism
at the level of fibers. The fiber of $U^{\ren,\tau}(L_{\fg,crit})$,
viewed as an object of $\CA_{\fg}^{\ren,\tau}\text{--}\on{mod}$,
corresponds, under the equivalence of categories given by
\corref{total equivalence}, to $\BV_{\fg,crit}\in
\CA_{\fg}^\flat\text{--}\on{mod}$. Hence it remains to show that 
$(\fD_{G,crit,x})^{\fg(\wh{\CO}_x)}\simeq \BV_{\fg,crit}$, but this is 
the content of \lemref{two centralizers}.

\end{proof}

This lemma implies, in particular, that for any $\CM\in 
\fD_{G,crit}\text{--}\on{mod}^{G(\wh{\CO}_x)}$, the chiral action of
$\fD^{\nil}_{G,crit}$ maps $\CM^{\nil}$ to $\CM^{\nil}$ and $\CM^{\on{non-reg}}$ to
$\CM^{\on{non-reg}}$.

\begin{lem}  \label{orthogonality}
For $\CM\in \fD_{G,crit}\text{--}\on{mod}^{G(\wh{\CO}_x)}$, the chiral action of
$\fD^{\on{non-reg}}_{G,crit}$ maps $\CM^{\nil}$ to $\CM^{\on{non-reg}}$.
\end{lem}

\begin{proof}

According to \lemref{Sugawara calculation}, we can find
an element of $\fZ_{\fg,x}$, such that its action  
is nilpotent on $\fD^{0}_{G,crit,x}$ and invertible on
$\fD^{\on{non-reg}}_{G,crit,x}$. We can assume that this element
comes from a local section $a\in \fz_\fg$. For example,
$a$ can be taken to be the section corresponding to the
Segal-Sugawara $S_0$ operator.

Moreover, we can find a section $a'$ of $\fz_\fg\boxtimes \CO_X$
such that the $\CO$-module endomorphism of $\fD_{G,crit}$ 
given by 
\begin{equation} \label{action of a'}
b\mapsto (h\boxtimes \on{id})[a'\otimes b]
\end{equation}
is nilpotent on $\fD^{0}_{G,crit}$, and invertible on
$\fD^{\on{non-reg}}_{G,crit}$. In the above formula
$a'\otimes b$ is viewed as an element of the D-module 
$\fD_{G,crit}\boxtimes \fD_{G,crit}$ on $X\times X$, $[\cdot,\cdot]$
denotes the chiral bracket, and $(h\boxtimes \on{id})$ denotes the
De Rham projection $\Delta_!(\fD_{G,crit})\to \fD_{G,crit}$.

The chiral action gives rise to a map of D-modules on $X$
$$\varphi:j_x{}_*j_x^*(\fD^{\on{non-reg}}_{G,crit})\otimes \CM^{\nil}\to i_x{}_!(\CM),$$
where $i_x$ (resp., $j_x$) is the embedding of the point $x$ (resp., of
its complement). For a section $b\in j_x{}_*j_x^*(\fD^{\on{non-reg}}_{G,crit})$
and an element $m\in \CM^{\nil}$, consider the section
$$a'\otimes b\otimes m\in  
j_{2,x}{}_*j_{2,x}^*(\fz_\fg\boxtimes\fD^{\on{non-reg}}_{G,crit})\otimes \CM^{\nil},$$
where $j_{2,x}$ is the embedding of the complement to
$\Delta_X\cup X\times x$ into $X\times X$. 

By applying the Jacobi identity to the above section,
we obtain that the action of $a'$ on the 
image of $\varphi$, given by the same formula as \eqref{action of a'}, is invertible.
But this means that the subspace of $\CM$ corresponding to the
D-submodule $\on{Im}(\varphi)\subset i_x{}_!(\CM)$ is supported off
$\on{Spec}(\fz_{\fg,x})$. Therefore, this subspace belongs to $\CM^{\on{non-reg}}$.

\end{proof}

\ssec{}

Let us assume for a moment that for any non-zero 
$\CM \in \fD_{G,crit}\text{--}\on{mod}^{G(\wh{\CO}_x)}$, 
the component $\CM^{\nil}\simeq \wt\Gamma(\CM)$ is necessarily non-zero.
Let us show that the functor $\wt\Gamma$ is then full.

Let $\wh\fD_{G,crit,x}$ be the canonical associative algebra 
corresponding to the chiral algebra $\fD_{G,crit}$ and the point 
$x\in X$, see \secref{functors}. We have a decomposition
$\wh\fD_{G,crit,x}=\wh\fD^{\nil}_{G,crit,x}\oplus \wh\fD^{\on{non-reg}}_{G,crit,x}$,
where the first summand is a subalgebra. For $\CM \in 
\fD_{G,crit}\text{--}\on{mod}^{G(\wh{\CO}_x)}$, the action of
$\wh\fD^{\nil}_{G,crit,x}$ preserves the decomposition $\CM=\CM^{\nil}\oplus \CM^{\on{non-reg}}$;
moreover, by \lemref{orthogonality}, the action of $\wh\fD^{\on{non-reg}}_{G,crit,x}$
sends $\CM^{\nil}$ to $\CM^{\on{non-reg}}$.

Observe that for
$\CM\in \fD_{G,crit}\text{--}\on{mod}^{G(\wh{\CO}_x)}$, the map
\begin{equation} \label{generating map}
\wh\fD_{G,crit,x}\otimes \CM^{\nil}\to \CM
\end{equation}
is automatically surjective. Indeed, its image is a
$\fD_{G,crit}$-submodule $\CM_1\subset \CM$, which satisfies
$\CM_1^{\nil}=\CM^{\nil}$. But then, for the quotient module 
$\CM_2:=\CM/\CM_1$, we have: $\CM_2^{\nil}=0$, which implies $\CM_2=0$.

Similarly, for any element $m\in \CM$ we can always find 
a section $b$ of $\wh\fD_{G,crit,x}$, such that $b\cdot m$ is a 
non-zero element of $\CM^{\nil}$.

\medskip

Let now $\CM$ and $\CN$ be two objects of 
$\fD_{G,crit}\text{--}\on{mod}^{G(\wh{\CO}_x)}$, and let
$\phi:\CM^{\nil}\to \CN^{\nil}$ be a map in $\CA^{\ren,\tau}_\fg\text{--}\on{mod}$.
By \lemref{neutral component}, $\phi$ is a homomorphism
of $\wh\fD^{\nil}_{G,crit,x}$-modules.
We have to show that this map extends uniquely to a map of
$\wh\fD_{G,crit,x}$-modules $\CM\to\CN$. 

The uniqueness statement is clear from the surjectivity of
\eqref{generating map}. To prove the existence, let us
suppose by contradiction that the required extension does not exist.
This means that there exist elements $a_i\in \wh\fD_{G,crit,x}$,
and $m_i\in \CM^{\nil}$, such that
$\underset{i}\sum\, a_i\cdot m_i=0\in \CM$, but
$\underset{i}\sum\, a_i\cdot \phi(m_i)=n\neq 0$ in $\CN$. Let $b\in \wh\fD_{G,crit,x}$
be an element such that $0\neq b\cdot n\in \CN^{\nil}$. Let us write
$b\cdot a_i=:c_i=c'_i+c''_i$, where $c'_i\in \wh\fD^{\nil}_{G,crit,x}$,
and $c''_i\in \wh\fD^{\on{non-reg}}_{G,crit,x}$.

By \lemref{orthogonality}, we have $\underset{i}\sum\, c'_i\cdot \phi(m_i)\neq 0$.
However, for the same reason, $\underset{i}\sum\, c'_i\cdot m_i=0$, which contradicts
the fact that $\phi$ was a morphism of $\wh\fD^{\nil}_{G,crit,x}$-modules.

\ssec{}

Finally, let us show that $0\neq \F\in
\on{D}_{crit}(\Gr_G)\text{--}\on{mod}$, implies that $\CM_\F^{\nil}\neq 0,$
where $\CM_\F$ is the corresponding object of
$\fD_{G,crit}\text{--}\on{mod}^{G(\wh{\CO}_x)}$. By \thmref{main}
and \thmref{total equivalence}, this is equivalent to the fact that
$$0\neq \F\in
\on{D}_{crit}(\Gr_G)\text{--}\on{mod}\Rightarrow \Gamma(\Gr_G,\F)\neq 0.$$
Note that the same argument works also in the negative and
irrational level cases:

For a congruence subgroup $K\subset G(\wh{\CO}_x)$, let 
$\on{D}_{\kappa}(\Gr_G)\text{--}\on{mod}^K$ be the subcategory of 
(strongly) $K$-equivariant D-modules, and 
let $\hat\fg_{\kappa}\text{--}\on{mod}^K$ be the subcategory of
$K$-integrable modules. The functor $\Gamma$ of global sections
evidently maps $\on{D}_{\kappa}(\Gr_G)\text{--}\on{mod}^K$ to
$\hat\fg_{\kappa}\text{--}\on{mod}^K$.

Recall now the setting for the Harish-Chandra action of \cite{BD1}, Sect. 7.14.
Namely, let $G((t))$ be the loop group corresponding to the point $x\in X$,
and $G((t))/K$-the corresponding ind-scheme.
We have the convolution functor
$$\star:\on{D}_{\kappa}(G((t))/K)\text{--}\on{mod}\times
\on{D}_{\kappa}(\Gr_G)\text{--}\on{mod}^K\longrightarrow 
D^b(\on{D}_{\kappa}(\Gr_G)\text{--}\on{mod}),$$
where $D^b(\cdot )$ stands for the bounded derived category. 
In addition, we have the functor
$$\star:\on{D}_{\kappa}(G((t))/K)\text{--}\on{mod}\times
\hat\fg_{\kappa}\text{--}\on{mod}^K\longrightarrow
D^b(\hat\fg_{\kappa}\text{--}\on{mod}).$$
Moreover, the (derived) functor of global sections
$$R\Gamma:
D^b(\on{D}_{\kappa}(\Gr_G)\text{--}\on{mod}^K)\to 
D^b(\hat\fg_{\kappa}\text{--}\on{mod}^K)$$ 
intertwines the two actions. 

\begin{lem}  \label{convolution}
For any non-zero object $\F\in \on{D}_{\kappa}(\Gr_G)\text{--}\on{mod}^K$,
there exists a $G(\wh{\CO}_x)$-equivariant object 
$\F'\in \on{D}_{\kappa}(G((t))/K)\text{--}\on{mod}$, such that
$\F'\star \F\in D^b(\on{D}_{\kappa}(\Gr_G)\text{--}\on{mod}^{G(\wh{\CO}_x)})$
is non-zero.
\end{lem}

\begin{proof}

We have an equivalence of categories $$\F\mapsto \F^*:
\on{D}_{\kappa}(\Gr_G)\text{--}\on{mod}^K\to 
\on{D}_{-\kappa}(G((t))/K)\text{--}\on{mod}^{G(\wh{\CO}_x)},$$
corresponding to the involution $g\mapsto g^{-1}$ on $G((t))$.
For $\F\in \on{D}_{\kappa}(\Gr_G)\text{--}\on{mod}^K$ and
$\F'\in \on{D}_{\kappa}(G((t))/K)\text{--}\on{mod}^{G(\wh{\CO}_x)}$,
the fiber at $1\in \Gr_G$ of the convolution $\F'\star \F$ is
canonically isomorphic to $\on{H}_{DR}(G((t))/K,\F'\otimes \F^*)$. (Note that
$\F'\otimes \F^*$ is an object of the derived category of
usual (i.e. non-twisted) right D-modules on $G((t))/K$, therefore,
global cohomology makes sense.)

In particular, this global cohomology is non-zero for 
$\F'$ being (the direct image of) the constant D-module on a 
$G(\wh{\CO}_x)$-orbit $G(\wh{\CO}_x)\cdot g \subset G((t))/K$,
for some $g\in G((t))/K$, such that the fiber $(\F^*)_g$ is non-zero.

\end{proof}

Using this lemma, our non-vanishing assertion reduces to the fact that
for a non-zero $\F\in D^b(\on{D}_{\kappa}(\Gr_G)\text{--}\on{mod}^{G(\wh{\CO}_x)})$,
the object $R\Gamma(\Gr_G,\F)$ is non-zero either. 

To prove it, note that since the functor $\Gamma$ is exact,
we can assume that $\F$ belongs to the abelian category of D-modules.
By the semi-smallness result \cite{BD1}, Sect. 5.3.6, the convolution 
$\star$ is exact on $\on{D}_{\kappa}(\Gr_G)\text{--}\on{mod}^{G(\wh{\CO}_x)}$,
i.e., $\on{D}_{\kappa}(\Gr_G)\text{--}\on{mod}^{G(\wh{\CO}_x)}$ acquires
a structure of monoidal category. (Note that for $\kappa$ integral,
the Satake equivalence identifies 
$\on{D}_{\kappa}(\Gr_G)\text{--}\on{mod}^{G(\wh{\CO}_x)}$
with the category of representation of the Langlands dual group $^L G$.)

For an object $\F\in \on{D}_{\kappa}(\Gr_G)\text{--}\on{mod}^{G(\wh{\CO}_x)}$
(which we can assume to be finitely generated), let 
$\F^*\in \on{D}_{-\kappa}(\Gr_G)\text{--}\on{mod}^{G(\wh{\CO}_x)}$
be the object as in the proof of \lemref{convolution},
and take $\F'\in \on{D}_{\kappa}(\Gr_G)\text{--}\on{mod}^{G(\wh{\CO}_x)}$ be 
the Verdier dual of $\F^*$. 

Let $\delta_1$ be the delta-function twisted D-module, corresponding to the 
unit point $1\in \Gr_G$. By adjunction, we obtain a non-zero map
$\delta_1\to \F'\star \F$. This map is necessarily an injection, because
$\delta_1$ is irreducible in $\on{D}_{\kappa}(\Gr_G)\text{--}\on{mod}$.
Hence, by the exactness of $\Gamma$, we obtain
$$\BV_{\fg,\kappa}\simeq \Gamma(\Gr_G,\delta_1)\neq 0\Rightarrow 
\Gamma(\Gr_G,\F'\star \F)\neq 0,$$
which, in turn, implies that $\Gamma(\Gr_G,\F)\neq 0$.


\begin{thebibliography}{199}

\bibitem[AG]{AG} 
S. ~Arkhipov and D.~Gaitsgory, {\em Differential operators on the loop group 
via chiral algebras},  IMRN 2002, no. 4, 165--210.

\bibitem[BB]{BB} 
A.~Beilinson and J.~Bernstein, {\em Localisation de $g$-modules}, 
C.R.Acad.Sci.Paris Ser. I Math. 292 (1981), no. 1, 15--18.

\bibitem[CHA]{BD} A.~Beilinson and V.~Drinfeld, {\em Chiral algebras},
Preprint, available at www.math.uchicago.edu/$\sim$benzvi.

\bibitem[BD]{BD1} 
A.~Beilinson and V.~Drinfeld, 
{\em Quantization of Hitchin's integrable system and Hecke eigensheaves},
Preprint, available at www.math.uchicago.edu/$\sim$benzvi.

\bibitem[EF]{EF} 
D.~Eisenbud and E.~Frenkel, Appendix to
{\em Jet schemes of locally complete intersection canonical singularities},
by M.Mustata, Inv. Math. {\bf 145} (2001) 397--424.

\bibitem[FB]{FB} E.~Frenkel and D.~Ben-Zvi, {\em Vertex Algebras and
Algebraic Curves}, Mathematical Surveys and
Monographs {\bf 88}, AMS, 2001.

\bibitem[FF]{FF} B.~Feigin and E.~Frenkel, {\em Affine Kac-Moody
algebras at the critical level and Gelfand-Dikii algebras}, in {\em
Infinite Analysis}, eds. A.  Tsuchiya, T. Eguchi, M. Jimbo, Adv.
Ser. in Math. Phys. {\bf 16}, 197--215, Singapore: World Scientific,
1992.

\bibitem[F]{F} 
E.~Frenkel,
{\em Lectures on Wakimoto modules, opers and the center at the
critical level}, math.QA/0210029.

\bibitem[FT]{FT} 
E.~Frenkel and K.~Teleman, in preparation.

\bibitem[KK]{KK}
V.~Kac and D.~Kazhdan, {\em Structure of representations with
highest weight of infinite-dimensional Lie algebras}, Adv. in Math.
{\bf 34} (1979), 97--108.

\end{thebibliography}
\end{document}